3/26/2022

# Sequences of Nahm pole solutions to the SU(2) Kapustin-Witten equations

Clifford Henry Taubes[†]

Department of Mathematics
Harvard University
Cambridge, MA 02138 USA

chtaubes@math.harvard.edu

ABSTRACT: This paper describes the behavior of sequences of solutions to the Kapustin-Witten equations with Nahm pole asymptotics on the product of $(0,\infty)$ with a compact, oriented, Riemannian 3-manifold. (Each solution consists of a connection on a principle SU(2) bundle on the product manifold and a Lie algebra valued 1-form. These have prescribed singular asymptotics near 0 in the $(0,\infty)$ factor and prescribed behavior near $\infty$ in this same factor.) In short: These sequences have subsequences that either converge to another solution after acting termwise by an automorphism of the principle bundle, or they converge after renormalization to a (weak) $\mathbb{Z}/2$ harmonic 1-form from the 3-manifold (it is independent of the $(0,\infty)$ coordinate in the product structure).

[†] Supported in part by the National Science Foundation

## 1. The Kapustin-Witten equations with Nahm pole boundary conditions

Let X denote an oriented 4-dimensional manifold with a Riemannian metric; and let P denote a principal G bundle over X with G being a compact Lie group. The Kapustin-Witten equations [KW] are equations for a pair (to be denoted by $(A, a)$) of connection on P and 1-form on X with values in the vector bundle ad(P). This is the associated bundle via G's adjoint representation whose fiber is the Lie algebra of G. The Kapustin-Witten equationa asserts that A's curvature 2-form (denoted by $F_A$) and the exterior covariant derivative of $a$ as defined by A (denoted by $D_A a$) obey

- $F_A - a \wedge a = *_X D_A a$ .
- $D_A *_X a = 0$ .

$$(1.1)$$

In this equation, $*_X$ denotes the metric's Hodge star. Note that $F_A$ and $D_A a$ and $a \wedge a$ are all ad(P) valued 2-forms. Note also that $*_X$ in dimension 4 maps 2-forms to 2-forms and has square equal to the identity (and it maps 1 forms to 3 forms, so $D_A *_X a$ is a 4-form.)

Of principle interest for what follows is the case when X is the product manifold $(0, \infty) \times Y$ with the metric being the product of the Euclidean metric on the half-line and a Riemannian metric on Y. Also: This paper only considers the case when G is SU(2). (The equations in (1.1) in the case of any compact group G are not so interesting when X is compact: In this case, Kapustin and Witten [KW] show that the only solutions are those where $F_A - a \wedge a$ and $D_A a$ are all zero; which is to say that $A + ia$ defines a flat $G_{\mathbb{C}}$ connection on X.) Witten, in a series of papers and lectures (see [W1], [W2]) following on work of Giaotto and Witten [GW], suggested that the space of solutions to (1.1) should be particularly interesting in the $(0, \infty) \times Y$ case when $(A, a)$ obey certain constraints when the coordinate on $(0, \infty)$ (which is denoted by t here) limits to 0 and when it limits to $\infty$:

- *The pair* $(A, a)$ *is asymptotic to the Nahm pole as* $t \to 0$.
- *The restrictions of the* $SL(2; \mathbb{C})$ *connection* $A + ia$ *to the various* $t \in (0, \infty)$ *versions of* $\{t\} \times Y$ *converges to a flat* $SL(2; \mathbb{C})$ *connection on* Y *as* $t \to \infty$.

$$(1.2)$$

The Nahm pole asymptotics for the case G = SU(2) are described in the upcoming Definition 1.1. (See Witten [W1], [W2] and Mazzeo-Witten [MW] for the definition when G is not SU(2).) With regards to the $t \to \infty$ constraint: The upcoming Definition 1.2 makes a constraint on the $t \to \infty$ limit of $(A, \mathfrak{a})$ that is weaker than (1.2)'s second bullet constraint.

This article contributes to the Nahm pole story with the upcoming Theorem A and Theorem B that describe (in part) the behavior of sequences of solutions to (1.1) that obey Definition 1.1's interpretation of the top bullet in (1.2) and Definition 1.2's weaker version of the lower bullet in (1.2).

### a) The asymptotic constraints

To set the stage for the definitions: The Lie algebra of SU(2) is identified in what follows with the vector space of $2 \times 2$ anti-Hermitian, trace zero matrices. The inner product of two matrices in this Lie algebra is then defined to be $-\frac{1}{2}$ times their trace. This inner product induces an inner product on the bundle ad(P) and, by extension, a homomorphism from $\otimes_2(\mathrm{ad}(P) \otimes T^*Y)$ to $\otimes_2 T^*Y$. In particular, supposing that $\mathfrak{a}$ is an ad(P) valued 1-form on Y, then the image of $\mathfrak{a} \otimes \mathfrak{a}$ by the latter homomorphism is a symmetric bilinear form which is denoted by $\langle \mathfrak{a} \otimes \mathfrak{a} \rangle$.

To continue with the stage setting: If A is a connection on P, then it induces a covariant derivative on ad(P) and, with the metric's Levi-Civita connection, a covariant derivative on sections of ad(P) valued tensor bundles. For example, on $\mathrm{ad}(P) \otimes T^*Y$. All of these covariant derivatives are denoted by $\nabla_A$. These covariant derivatives along just the Y-factor of $(0, \infty) \times Y$ are all denoted by $\nabla_A{}^{\perp}$.

Here one last bit of stage setting: An ad(P) valued 1-form on $(0,\infty)\times Y$, call it $a$, will be written as $a = a_0 dt + \mathfrak{a}$ with $\mathfrak{a}$ annihilating tangents to the $(0,\infty)$ factor of $(0,\infty)\times Y$.

**Definition 1.1**: *Let* P *denote a principal SU(2) bundle over* $(0, \infty) \times Y$ *and let* $(A, a_0 dt + \mathfrak{a})$ *denote a pair of connection on* P *and* ad(P)*-valued 1-form. This pair is said to be asympotic to the Nahm pole when the following conditions are met* (*the limits in the first two bullets are uniform over* Y*)*:

- $\lim_{t\to 0} |a_0| = 0$ .
- $\lim_{t\to 0} 4t^2\langle \mathfrak{a} \oplus \mathfrak{a} \rangle$ *is the Riemannian metric on* Y.
- $\int_{(0,1]\times Y} (|F_A|^2 + |\nabla_A^{\perp}\mathfrak{a}|^2)$ *is finite*.

This definition of the Nahm pole assymptotics for a solution to (1.1) for the group SU(2) is equivalent to the definition of the Nahm pole boundary condition given by Mazzeo and Witten [MW]. The equivalence is exhibited in the appendix to this article.

What follows is the definition that is used here for the $t \to \infty$ requirement in (1.2).

**Definition 1.2**: *Let* P *denote a principal SU(2) bundle over* $[0, \infty) \times Y$ *and let* $(A, a_0 dt + \mathfrak{a})$ *denote a pair of connection on* P *and* ad(P)*-valued 1-form. This pair is said to be asympotically like a flat* $\mathrm{Sl}(2;\mathbb{C})$ *connections on* Y *as* $t \to \infty$ *when*

- $\lim_{t\to\infty} |a_0| = 0$.
- $\int_{[1,\infty)\times Y} (|D_A a|^2 + |F_A - a\wedge a|^2)$ *is finite.*

With regards to this definition: Any SL(2;ℂ) connection on $(0, \infty) \times Y$ can be written as $A + ia$ with A being a connection on P and $a$ being an ad(P) valued 1-form. The connection is flat if and only if $D_A a = 0$ and $F_A = a\wedge a$. Therefore, if $(A, a)$ has a $t \to \infty$ limit which defines a flat SL(2;ℂ) connection on Y; and if this limit is approached sufficiently fast as $t \to \infty$, then the integral in the second bullet will be finite. However, if the integral in the second bullet is finite, then the various $\{t\}\times Y$ restrictions of $A + ia$ need not converge as $t \to \infty$ to a flat SL(2;ℂ) connection on Y except in the event that the Ricci curvature of Y is positive. Thus, the condition in the definition may be weaker than what is suggested by (1.2). (See [T1] for what will happen if there is no converging subsequence.)

Supposing that $(A, a)$ is a pair of connection on P and ad(P) valued 1-form over $(0,\infty)\times Y$. Use these to define a function of t (denoted by $\mathfrak{cs}$) by the rule whereby

$$\mathfrak{cs}(t) = \int_{\{t\}\times Y} (\langle a\wedge F_A\rangle - \tfrac{1}{3}\langle a\wedge a\wedge a\rangle) \tag{1.3}$$

As explained in Section 2d of [T2], if $(A, a)$ is a solution to (1.1) on $(0, \infty)$, then this function is non-decreasing. Moreover, if the condition in the second bullet of Definition 1.2 holds, then $\mathfrak{cs}(\cdot)$ has a well defined $t \to \infty$ limit. (If $(A, a)$ obeys (1.1), then $\mathfrak{cs}(s) - \mathfrak{cs}(t)$ for $s > t$ is half of the $[t,s] \times Y$ integral of $|D_A a|^2 + |F_A - a\wedge a|^2$.) The $t \to \infty$ limit of $\mathfrak{cs}$ is denoted by $\mathfrak{cs}_\infty$.

With Definitions 1.2 and 1.3 in hand, here is a third definition:

**Definition 1.3**: *A pair* $(A, a_0 dt + \mathfrak{a})$ *that obeys (1.1) and is described by Definitions 1.1 and 1.2 will be said to be a* *<u>Nahm pole solution</u>*.

As observed by [He], the top bullets in (1.1) and (1.2) require that the $a_0$ part of a Nahm pole solution is identically zero (see the upcoming Lemma 2.2). Keep this in mind.

**b) ℤ/2 harmonic 1-forms**

A ℤ/2 harmonic 1-form on Y consists of a data set $(Z, \mathcal{I}, \nu)$ with Z, $\mathcal{I}$ and $\nu$ as follows: First, Z is a closed, nowhere dense subset of Y. Second, $\mathcal{I}$ is a real line bundle on Y−Z. (This means that $\mathcal{I}$ is the associated $\mathbb{R}$-bundle to a principal {±1} bundle.) Third, $\nu$ is an $\mathcal{I}$-valued 1-form on Y−Z that is both closed and coclosed. (The bundle $\mathcal{I}$ has a canonical metric and metric compatible covariant derivative because it is associated to a principle ±1 bundle.) In addition, the norm of $\nu$ extends to the whole of Y as a Holder continuous function whose zero locus contains Z. Zhang [Zh] proved that Z is necessarily 2-rectifiable; in particular, it has finite 2-dimensional Hausdorff measure. It is also known (see [T3]) that Z contains a (relatively) open, dense subset that has the structure of a $C^1$ submanifold in Y.

By way of a local example: Use $(x_1, x_2, x_3)$ for the Euclidean coordinates of $\mathbb{R}^3$. Take Z to be the $x_3$ axis, and let $\mathcal{I}$ denote the line bundle on the complement of the $x_3$ axis with holonomy -1 around the unit circle in the $x_1$-$x_2$ plane. Then, set $f$ to be the real part of the $\mathcal{I}$-valued function $(x_1 + ix_2)^{3/2}$ and set $\nu$ to be $df$.

A $\mathbb{Z}/2$ harmonic form data set $(Z, \mathcal{I}, \nu)$ on Y defines an analogous data set on $(0,\infty)\times$ Y via pull-back by the projection to Y. This alternate view of $(Z, \mathcal{I}, \nu)$ is implicitly used below.

Now suppose that $(Z, \mathcal{I}, \nu)$ is a $\mathbb{Z}/2$ harmonic 1-form data set on Y. Let $\mathcal{T}$ denote the symmetric, bilinear form $\nu\otimes\nu$ (it is an honest section of $\otimes_2 T^*Y$.) This is smooth on the complement of Z where it is characterized by the algebraic and differential conditions stated momentarily in (1.4).

The conditions in the upcoming (1.4) refer to two differential operators on Y that act on sections of $\otimes_2 T^*Y$. The first is denoted by div. It is the formal, $L^2$ adjoint of the map from 1-forms to 2-tensors that sends any given 1-form to its covariant derivative (which is a priori a section of $\otimes_2 T^*Y$). When written using a local orthonormal frame $\{e^i\}_{i=1,2,3}$ for $T^*Y$, the 1-form $\mathrm{div}(\mathcal{T})$ is $\nabla_k \mathcal{T}_{ki} e^i$ where $\nabla_k$ denotes the Levi-Civita connection's directional covariant derivative along the vector field dual to $e^k$. (Repeated tensor indices are summed implicitly here and in what follows.) The second differential operator is denoted by curl. This sends sections of $\otimes_2 T^*Y$ to sections of $\otimes_2 T^*Y$. When written using the orthonormal frame, $\mathrm{curl}(\mathcal{T})$ has $e^i \otimes e^j$ component $\varepsilon^{imn}\nabla_m \mathcal{T}_{nj}$ with $\varepsilon$ denoting the volume 3-form. Two final bits of notation: If $\mathcal{T}$ is a section of $\otimes_2 T^*Y$, then $\mathrm{tr}(\mathcal{T})$ denotes its inner product with the metric (its trace). And, if $\mathcal{V}$ is another section of $\otimes_2 T^*Y$, then $\mathcal{T}\bullet\mathcal{V}$ is the section with $e^i \otimes e^j$ component given by $\mathcal{T}_{im}\mathcal{V}_{mj}$.

The following conditions characterize a $\mathbb{Z}/2$ harmonic 1-form using $\mathcal{T} = \nu\otimes\nu$:

- $\mathrm{tr}(\mathcal{T}\bullet\mathcal{T}) - \mathrm{tr}(\mathcal{T})^2 = 0$.
- $\mathrm{div}(\mathcal{T}) - \frac{1}{2}\, d\,\mathrm{tr}(\mathcal{T}) = 0$.
- $\mathcal{T}\bullet\mathrm{curl}(\mathcal{T}) = 0$.

(1.4)

Here, d is the exterior derivative on Y. (The top bullet says that the symmetric tensor $\mathcal{T}$ has rank 1; and the lower two bullets say that if $\mathcal{T}$ is written locally as $\nu \otimes \nu$ with $\nu$ being a 1-form, then $\nu$ must obey $d\nu = 0$ and $d{*}\nu = 0$. With regards to the third condition, it is sufficient that $\mathrm{tr}(\mathcal{T} \bullet \mathrm{curl}(\mathcal{T}))$ is zero.)

**Definition 1.4**: *A an almost everywhere bounded, symmetric section of* $\otimes_2 T^*Y$ *(to be denoted by $\mathcal{T}$) is said to be a weak, square-harmonic section of* $\otimes_2 T^*Y$ *when* $|\nabla \mathcal{T}|^2$ *and* $|\mathcal{T}|^2$ *are integrable on* Y *and when it obeys (1.3) in the following weak sense:*

- $\int_Y |\mathrm{tr}(\mathcal{T} \bullet \mathcal{T}) - (\mathrm{tr}\mathcal{T})^2| = 0$ .
- $\int_Y |\mathrm{div}\mathcal{T} - \frac{1}{2} d\, \mathrm{tr}(\mathcal{T})|^2 = 0$.
- $\int_Y |\mathcal{T} \bullet \mathrm{curl}(\mathcal{T})| = 0$.

There is no regularity theorem at present asserting in effect that a weak square-harmonic section of $\otimes_2 T^*Y$ comes from an honest $\mathbb{Z}/2$ harmonic 1-form data set. (The issue is, for the most part, whether $\mathcal{T}$ is necessarily continuous.)

## c) The convergence theorems

The following theorem is the first contribution of this paper.

**Theorem A**: *Let* $\{(A_n, a_n = \mathfrak{a}_n)\}_{n=1,2,\ldots}$ *denote a sequence of Nahm pole solutions on* $(0,\infty)\times Y$ *whose corresponding sequence of* $\mathfrak{cs}_\infty$ *values is bounded from above. For each* $n \in \{1, 2, \ldots\}$, *let* $\kappa_n(1)$ *denote the* $L^2$ *norm of* $\mathfrak{a}_n$ *on* $\{1\}\times Y$.

- *Assume that* $\{\kappa_n(1)\}_{n=1,2,\ldots}$ *has a bounded subsequence. There exists a subsequence* $\Lambda \subset \{1, 2, \ldots\}$, *a corresponding sequence* $\{g_n\}_{n\in\Lambda}$ *of automorphisms of* P *over* $(0, \infty)\times Y$, *and there exists a Nahm pole solution to (1.1) to be denoted by* $(A, a = \mathfrak{a})$; *and these are such that the sequence* $\{(g_n{}^*A_n, g_n{}^*\mathfrak{a}_n)\}_{n\in\Lambda}$ *converges to* $(A, \mathfrak{a})$ *in the following sense*:
  - a) *The associated sequence* $\{(g_n{}^*A_n - A,\ g_n{}^*\mathfrak{a}_n - \mathfrak{a})\}_{n\in\Lambda}$ *converges to zero in the* $C^\infty$-*topology on compact subsets of* $(0, \infty)\times Y$
  - b) $\lim_{n\in\Lambda} \int_0^1 (|\nabla_A(g_n^* A_n - A)|^2 + |g_n^* A_n - A|^2) = 0$ .
  - c) $\lim_{n\in\Lambda} \int_0^1 (|\nabla_A(g_n^* \mathfrak{a}_n - \mathfrak{a})|^2 + |g_n^* \mathfrak{a}_n - \mathfrak{a}|^2) = 0$ .
  - d) *The sequence of symmetric, bilinear forms* $\{t^2\langle \mathfrak{a}_n \otimes \mathfrak{a}_n\rangle\}_{n\in\Lambda}$ *converges uniformly in the* $C^0$ *topology on* $[0,1]\times Y$ *to* $t^2\langle \mathfrak{a} \otimes \mathfrak{a}\rangle$.

- *Assume that* $\{K_n(1)\}_{n=1,2,\ldots}$ *lacks bounded subsequences. Then, there exists a non-zero weak square-harmonic section of* $\otimes_2 T^*Y$ *on* Y *(to be denoted by* $\mathcal{T}$*), and there exists a subsequence* $\Lambda \subset \{1, 2, \ldots\}$*; and these are such that the sequence*

$$\{ \tfrac{1}{K_n(1)^2} \langle a_n \otimes a_n \rangle \}_{n\in\Lambda}$$

*converges to* $\mathcal{T}$ *weakly in the Sobolev* $L^2_1$ *topology on compact subsets of* $(0, \infty)\times Y$.

There are very good reasons to think that $\{(A_n, \frac{1}{K_n(1)} \mathfrak{a}_n)\}_{n\in\Lambda}$ converges in the manner of Theorem 1.2 in [T2] so as to define an honest $\mathbb{Z}/2$ harmonic 1-form data set $\{Z, \mathcal{I}, \nu\}$ on Y (which is viewed as a $\mathbb{Z}/2$ harmonic 1-form data set on $(0, \infty)\times Y$ that is independent of the $(0,\infty)$ factor.) The author hopes to address this issue in a sequel to this article.

Theorem A can be strengthened in the case when Y has positive Ricci curvature:

**Theorem B**: *Suppose that* Y *has positive Ricci curvature. Let* $\{(A_n, a_n = \mathfrak{a}_n)\}_{n=1,2,\ldots}$ *denote a sequence of Nahm pole solutions on* $(0,\infty)\times Y$ *whose corresponding* $\mathfrak{cs}_\infty$ *sequence is bounded from above. There exists a subsequence* $\Lambda \subset \{1, 2, \ldots\}$*, a corresponding sequence* $\{g_n\}_{n\in\Lambda}$ *of automorphisms of* P *over* $(0, \infty)\times Y$*, and there exists a Nahm pole solution to (1.1) which is denoted by* $(A, a = \mathfrak{a})$*; and these are such that the sequence* $\{(g_n{}^*A_n, g_n{}^*\mathfrak{a}_n)\}_{n\in\Lambda}$ *converges to* $(A, \mathfrak{a})$ *in the following sense*:

a) *The associated sequence* $\{(g_n{}^*A_n - A, g_n{}^*\mathfrak{a}_n - \mathfrak{a})\}_{n\in\Lambda}$ *converges to zero in the* $C^\infty$*-topology on compact subsets of* $(0, \infty) \times Y$

b) $\lim_{n\to\infty} \int_0^1 (|\nabla_A (g_n^* A_n - A)|^2 + | g_n^* A_n - A|^2 ) = 0$ .

*c)* $\lim_{n\to\infty} \int_0^1 (|\nabla_A (g_n^* \mathfrak{a}_n - \mathfrak{a})|^2 + | g_n^* \mathfrak{a}_n - \mathfrak{a}|^2 ) = 0$ .

d) *The sequence of symmetric, bilinear forms* $\{t^2\langle \mathfrak{a}_n \otimes \mathfrak{a}_n\rangle\}_{n\in\Lambda}$ *converges uniformly in the* $C^0$ *topology on* $[0,1]\times Y$ *to* $t^2\langle \mathfrak{a} \otimes \mathfrak{a}\rangle$.

This stronger theorem for the positive Ricci curvature case was foreshadowed in the recent preprint of Leung and Takahashi [LT].

## d) The remainder of this article

This article has ten sections and an appendix. Looking ahead, the proofs of Theorems A and B are in Section 10. The appendix explains why Definition 1.1 is equivalent to the Mazzeo-Witten definition of Nahm pole asymptotics

Here is a table of contents:

**e) Looking ahead**

The first task is to prove (or disprove) that the second bullet case in Theorem A leads to an honest $\mathbb{Z}/2$ harmonic 1-form. To be precise, the second bullet of Theorem A should say this:

*Assume that* $\{K_n(1)\}_{n=1,2,\ldots}$ *lacks bounded subsequences. Then, there exists a* $\mathbb{Z}/2$ *harmonic 1-form data set* $(\mathcal{I}, Z, \nu)$ *on* $Y$, *and there exists a subsequence* $\Lambda \subset \{1, 2, \ldots\}$, *a corresponding sequence of isometric isomorphism* $\{\tau_n\colon \mathcal{I} \to \mathrm{adP}\}$ *which is defined over* $(0, \infty)\times(Y-Z)$*; and these are such that* $\{\frac{1}{K_n(1)} a_n - \tau_n \otimes \nu\}_{n\in\Lambda}$ *converges to zero in the* $C^\infty$ *topology on compact subsets of* $(0, \infty)\times Y-Z$. *Meanwhile, both* $\{F_{A_n}\}_{n\in\Lambda}$ *and* $\{\nabla_{A_n}\tau_n\}_{n\in\Lambda}$ *converge to* $0$ *in the* $C^\infty$ *topology on compact subsets of* $(0, \infty)\times Y-Z$.

(1.5)

(Remember that the assumptions are that $\{(A_n, \mathfrak{a}_n)\}_{n=1,2,\ldots}$ is a sequence of Nahm pole solutions on $(0,\infty)\times Y$ whose corresponding $\mathfrak{cs}_\infty$ sequence is bounded from above. For each $n \in \{1, 2, \ldots\}$, what is denoted by $\kappa_n(1)$ is $L^2$ norm of $\mathfrak{a}_n$ on $\{1\}\times Y$.) There are compelling hints that (1.5) is true. A sequel to this article should either prove it, or explain why it is false. (The convergence in (1.5) is the sort of convergence that is described in Theorem 1.2 of [T2]. But, more work is needed to justify the use of that theorem here.)

In a different direction: Witten in [W1], [W2] and Mazzeo and Witten in [MW2] describe an analog of the Nahm pole asymptotics that is dictated by a given knot in the boundary $\{0\}\times Y$. There are probably analogs of Theorems A and B (and maybe (1.5)) for this knot version of Nahm pole asymptotics. In any event, much of the analysis in Sections 4-9 and the appendix of this paper should carry over to the knot version of the Nahm pole asymptotics (with some changes).

On a final note: The $\mathbb{Z}/2$ harmonic 1-forms on Y appear in the context of the Nahm pole boundary conditions (probably with knots too); and they also appear (see [T1]) as renormalized limits of non-convergent sequence of flat $SL(2;\mathbb{C})$ connections on Y (modulo bundle automorphisms). These $\mathbb{Z}/2$ harmonic 1-forms are still mysterious objects. They should be studied. (See [Tak] for a first step.)

**f) Conventions**

Various conventions will be employed for the most part without further comment. Here is a list:

- It is assumed in what follows that Y has volume 1. As long as Y has finite volume, the precise volume has no substantive bearing on what transpires. And, if the volume of Y is not 1, then various factors of the volume of Y clutter subsequent inequalities and identities.

- The Riemannian metric defines a corresponding inner product on all tensor bundles; and the Riemannian metric with the metric on adP define inner products on their tensor product with adP. All of these inner products are denoted by $\langle\, ,\rangle$.

- What is denoted by $c_0$ in what follows is a number that is greater than 1. Unless told otherwise, it is independent of any particular point in $[0,\infty)\times Y$ and independent of any particular solution to (1.1). Its value can be assumed to increase between successive appearances. (The 'otherwise' appears after Section 5.)

- What is denoted by $\chi$ is a smooth, non-increasing function on $\mathbb{R}$ that is one on $(-\infty, \frac{1}{4}]$ and zero on $[\frac{3}{4}, \infty)$. This is your favorite, standardized cut-off function. All 'bump' functions and cut-off functions are (implicitly) made from $\chi$ by rescalings (using $\chi(rt)$ in lieu of $\chi$ with $r$ a positive number) and taking suitable products of rescalings. This insures that their derivatives to any given order have uniform norm bounds after accounting for the rescalings.

- The Euclidean coordinate for the $(0, \infty)$ factor $(0, \infty) \times Y$ is denoted by t for the most part, but sometimes by s or x. (Note that [MW] use y, but this can be confused with a point in Y.)

By way of a final convention: The Lie algebra of SU(2) (which is the vector space of 2×2 anti-Hermitian matrices) is denoted by $\mathfrak{su}(2)$. An oriented orthonormal frame for $\mathfrak{su}(2)$ when needed is denoted by $\{\sigma_1, \sigma_2, \sigma_3\}$. These matrices square to -1 times the identity; and

$$\sigma_1\sigma_2 = -\sigma_3 \ . \tag{1.4}$$

Note in particular the minus sign on the right hand side. This is not the convention taken in [MW] and [He], nor is it (evidently) the convention taken by most earthlings.†

†Please forgive this. I have used this convention since back when George Washington and I were kids together. A change now will increase the odds of a sign error in any given identity to nearly 1 (but the odds may be near to 1 anyway, so check signs as you read).

**g) Acknowledgements**

Many of the basic identities and inequalities in this article can also be found in the fundamental papers by Kapustin and Witten [KW],Witten [W1], [W2], and Mazzeo-Witten [MW]; and in the important subsequent work of Siqi He [He], and Leung and Takahashi [LT].

## 2. Fundamental identities

Various Bochner-Weitzenboch identities play a central role in subsequent arguments. These identities (or closely related ones) can be found in [KW], [MW] and/or [LT]. (See also Sections 2 and 3 of [T2].) This section states some of these identities and derives some of their first consequences.

### a) Second order equations for $a$ and $F_A$

The ad(P)-valued 1-form $a$ obeys a second order equation that can be written schematically as:

$$\nabla_A^{\dagger}\nabla_A a + [a_*,[a, a_*]] + \mathrm{Ric}(a) = 0 \ . \tag{2.1}$$

The notation is as follows: The symbol $\nabla_A^{\dagger}$ denotes the formal, $L^2$ adjoint of the covariant derivative $\nabla_A$. Meanwhile, what is denoted by $[a_*,[a, a_*]]$ is best described by writing $a$ locally using an oriented, orthonormal frame for $T^*Y$: $\{e^1, e^2, e^3\}$. Let $e^0$ denote dt so the tetrad $\{e^{\nu}\}_{\nu=0,1,2,3}$ is an orthonormal frame for $T^*(\mathbb{R}\times Y)$. Write $a$ using this frame as $a_{\nu} e^{\nu}$ where repeated indices are implicitly summed. (Each $a_{\nu}$ is a section of ad(P).) Let $b = b_{\nu}e^{\nu}$ denote either $a$ or some other ad(P) valued 1-form. The component along $e^{\nu}$ of $[a_*,[b, a_*]]$ is the double commutator $[a_{\mu}, [b_{\nu}, a_{\mu}]]$ with it understood again that repeated indices are summed. (The commutator endomorphism from $\mathfrak{su}(2)\otimes\mathfrak{su}(2)$ to $\mathfrak{su}(2)$ is equivariant with respect to the adjoint action of SU(2) so it defines, fiberwise, an endomorphism from ad(P)$\otimes$ ad(P) to ad(P).) A useful formula for $[a_*,[b, a_*]]$ is this:

$$[a_{\mu}, [b, a_{\mu}]] = 4\,(|a|^2 b - \langle a_{\mu} b\rangle\, a_{\mu}) \ . \tag{2.2}$$

The last bit of notation from (2.1) concerns the term Ric($a$). What is denoted by Ric signifies here the endomorphism on $T^*Y$ (extended to ad(P)$\otimes T^*Y$) that is obtained from the Ricci tensor of the metric on Y by using the metric to contract indices. (Using the orthonormal frame, $\mathrm{Ric}_{\nu}(a) = 0$ when $\nu = 0$; and $\mathrm{Ric}_{\nu}(a) = \sum_{i=1,2,3} \mathrm{Ric}_{\nu i}\, \mathfrak{a}_i$ otherwise.)

By way of a derivation: The identity depicted in (2.1) is obtained from the equations in (1.1) by noting first that $D_A{*}D_A a = 0$ because $D_A F_A = 0$ (the latter is the Bianchi identity). Given this, then there is the identity $*_X D_A *_X D_A a - D_A *_X D_A *_X a = 0$; and writing the latter in terms of the covariant Laplacian $\nabla_A^{\dagger}\nabla_A$ leads (with the help again of the top bullet of (1.1)) to what is asserted by (2.1).

A second order equation is also obeyed by the curvature 2-form $F_A$. To write this equation, fix for the moment an oriented, orthonormal frame for $T^*Y$ to give, with dt, the frame $\{e^\nu\}_{\nu=0,1,2,3}$ for $T^*((0,\infty)\times Y)$. Any given ad(P)-valued 2-form (call it $\Theta$) is written as $\frac{1}{2}\Theta_{\mu\nu}\, e^\mu \wedge e^\nu$ (repeated index summation is implicit). Meanwhile, the covariant derivative of $a$ along the vector dual to $e^\mu$ is written as $(\nabla_{A\mu} a)_\nu\, e^\nu$. What follows is the equation for $F_A$.

$$(\nabla_A{}^\dagger \nabla_A F_A)_{\mu\nu} + [a_*, [F_{A\mu\nu}, a_*]] + 2[F_{A\mu\gamma}, F_{A\nu\gamma}] + \mathcal{R}^1{}_{\mu\nu\alpha\beta} F_{A\alpha\beta} = \mathcal{R}^2{}_{\mu\nu\alpha\beta}[a_\alpha, a_\beta] - 2[(\nabla_{A\mu} a)_\gamma, (\nabla_{A\nu} a)_\gamma] \tag{2.3}$$

where $\mathcal{R}^1$ and $\mathcal{R}^2$ are endomorphisms of $\wedge^2 T^*Y$ whose components are linear combinations of components of the metric's Riemann curvature tensor. (Repeated indices are summed in (2.3) also.) The identity in (2.3) follows from the equation $\nabla_{A\mu} F_{A\mu\nu} = [a_\mu, \nabla_\nu a_\mu]$ which follows in turn from (1.1).

**b) Local regularity**

An important point to keep in mind is that the equations in (1.1) are elliptic when the gauge invariance is accounted for; and that the non-linearities are essentially quadratic. The proposition that follows states one consequence of these facts and the identities in (2.1) and (2.3). By way of notation, the third bullet of the proposition uses $\theta_0$ to denote the product connection on the product principle SU(2) bundle.

**Proposition 2.1**: *There exists* $\kappa > 1$ *with the following significance: Fix* $r \in (0, \kappa^{-1})$ *and let* $B \subset (0,\infty)\times Y$ *denote a ball of radius* r*; and let* $\partial B$ *denote its boundary. Assume that* t *is strictly positive on* $\partial B$. *Let* $(A, a)$ *denote a solution to (1.1) on a neighborhood of* B, *and let* $\hat{\kappa}$ *denote* $r^{-3/2}$ *times the* $L^2$ *norm of* $|a|$ *on* $\partial B$.

- *The pointwise bound* $|a| \leq \kappa\, \hat{\kappa}$ *holds on the concentric, radius* $\frac{15}{16}$r *ball.*
- *If* $\hat{\kappa} \leq \kappa^{-2} r^{-1}$, *then, on the concentric, radius* $\frac{3}{4}$r *ball:*
  - a) $|\nabla_A a| \leq \kappa r^{-2}$ .
  - b) $|F_A| \leq \kappa r^{-2}$.
  - c) *For each integer* $k \geq 1$, *there exists* $c_k > 1$ *which is independent of* B, r *and* $(A, a)$ *and such that the pointwise norm of the covariant derivatives to order* k *of* $\nabla_A a$ *and* $F_A$ *are bounded by* $c_k r^{-k+2}$.
- *If* $\hat{\kappa} \leq \kappa^{-2} r^{-1}$, *then, on the concentric, radius* $\frac{3}{4}$r *ball, there is an isomorphism (to be denoted by* u*) from the product principle SU(2) bundle to P such that for each integer* $k \geq 0$, *the* $\theta_0$*-covariant derivatives of* $\hat{a} = u^*A - \theta_0$ *and* $\hat{a} = u^*a$ *to order* k *are bounded by* $c_k r^{-k-1}$ *with* $c_k$ *being independent of* B, r, *and* $(A, a)$.

***Proof of Proposition 2.1***: These are proved using arguments that are along the same lines as those employed in Section 3 of [T2]. The first bullet is proved by taking the inner product of (2.1) with $a$ to obtain an equation reading:

$$\tfrac{1}{2}\nabla^{\dagger}\nabla|a|^2 + |\nabla_A a|^2 + |[a_*, a]|^2 + \langle a, \mathrm{Ric}(a)\rangle = 0 \tag{2.4}$$

where $|[a_*, a]|^2$ is $\sum_{0\leq\alpha,\beta\leq 3}|[a_\alpha, a_\beta]|^2$. This is also $\langle a, [a_*, [a, a_*]]\rangle$ and also $2|a\wedge \mathfrak{a}|^2$. (Here and in what follows, $\nabla$ denotes the Levi-Civita covariant derivative.) The equation in (2.4) is used (after some preliminaries) with the Dirichelet Green's function on B for $\nabla^{\dagger}\nabla$ to obtain the bound in the top bullet. (See what is done in Sections 3a and Section 3b of [T2] for proving Proposition 3.1 in [T2].)

A preliminary observation is needed to prove the second bullet: To this end, let B´ denote the concentric ball with radius $\frac{15}{16}$r. Integrating (2.4) against a suitable bump function supported in B and equal to 1 on B´ leads to a $c_0\,\hat{\kappa}^2\,r^2$ bound for the integral of $|\nabla_A a|^2$ on B´ (the argument for this invokes (3.11) of [T2]). Meanwhile, the bound in the first bullet leads to a $c_0\,\hat{\kappa}^4\,r^4$ bound for the integral of $|[a_*, a]|$ on B´. Therefore, given in advance $c > 1$, and supposing that $\hat{\kappa} \leq c_0^{-1}c^{-1}r^{-1}$, then both of these integrals will be smaller than $c^{-2}$. As a consequence of (1.1), so will the integral over B´ of $|F_A|^2$. (In particular, Equations (3.1) and (3.2) in [T2] are obeyed on B´.)

With the preceding understood, one proves Item a) of the second bullet using an argument that differs only minimally from that given in Part 2 of Section 3e) of [T2] for proving the second bullet of [T2]'s Proposition 3.3.) In short, the equation in (2.1) is differentiated to obtain a second order equation for $\nabla_A a$ whose inner product with $\nabla_A a$ leads to an inhomogeneous Laplace equation for $|\nabla_A a|^2$ that is analogous to (2.4). The Dirichelet Green's function on B´ is then used to get an a priori bound for $|\nabla_A a|$ from the latter equation. Item b) of the second bullet is proved by taking the inner product of (2.3) with $F_A$ and then using similar manipulations. Item c) is proved by differentiating (2.1) and (2.3) the appropriate number of times to obtain Laplace equations for the norms of the derivatives of $\nabla_A a$ and $F_A$.

The third bullet follows from the second using Karen Uhlenbeck's justly famous theorem in [Uh] to obtain a priori bounds on $u^*A - \theta_0$ for a suitable isomorphism $u$. Note in this regard that Uhlenbeck's theorem requires a $c_0^{-1}$ bound for the integral of $|F_A|^2$ on B; and such a bound follows from Item b) of the second bullet of the lemma if $\kappa > c_0$.

### c) The vanishing of $a_0$

The identity in (2.1) has the following immediate corollary (see [He]):

**Lemma 2.2**: *If* $(A, a_0\,dt + \mathfrak{a})$ *is a Nahm pole solution, then* $a_0$ *is identically zero.*

***Proof of Lemma 2.2***: The dt component of (2.1) says that $\nabla_A{}^\dagger\nabla_A a_0 + [a_*,[a_0, a_*]] = 0$. Taking the inner product with $a_0$ leads to the equation

$$\tfrac{1}{2}\nabla^\dagger\nabla|a_0|^2 + |\nabla_A a_0|^2 + |[a, a_0]|^2 = 0 \tag{2.5}$$

This identity implies (via the maximum principle) that $|a_0|$ has no local maxima unless it is constant. But, if it has no local maxima, it is zero because its $t \to 0$ and $t \to \infty$ limits are zero.

Note that if the Ricci curvature of Y is non-negative, then the (2.5) implies that $d^\dagger d|\mathfrak{a}|^2 \le 0$. Granted this, then the maximum principle says that $|\mathfrak{a}|^2$ has no local maxima. It also implies that the function of t given by the integral of $|\mathfrak{a}|^2$ on any $\{t\} \times Y$ is a non-increasing function of t. More will be said momentarily about this function of t (with no assumption on the Ricci curvature.)

### d) The 1 + 3 notation

Because $a_0 = 0$, it is convenient to write Nahm pole solutions as $(A, \mathfrak{a})$ with $\mathfrak{a}$ viewed as t-dependent, Ad(P) valued 1-form on Y. In fact, it proves useful to dispense entirely (almost) with 4-dimensional notation and distinguish separately the $(0,\infty)$ and Y components of not just $a$ but the curvature of A and the covariant derivative $\nabla_A$. The curvature 2-form of A (which is an ad(P) valued 2-form on $(0,\infty)\times Y$) will be written as

$$F_A = dt\wedge E_A + *B_A \tag{2.6}$$

with $E_A$ and $B_A$ denoting t-dependent, ad(P) valued 1-forms on Y, and with $*$ denoting here and in what follows the metric Hodge dual on Y. Meanwhile, the covariant derivative $\nabla_A$ will be written as $dt\,\nabla_{At} + \nabla_A{}^\perp$ with $\nabla_{At}$ differentiating along the $(0, \infty)$ factor of $(0,\infty) \times Y$ and $\nabla_A{}^\perp$ differentating along the Y factor. The exterior covariant derivative along the Y factor is written below as $d_A$. The equations in (1.1) can be written using this new notation as three equations:

- $E_A = *d_A\mathfrak{a}$ .
- $\nabla_{At}\mathfrak{a} = B_A - *(\mathfrak{a}\wedge\mathfrak{a})$ .

- $d_A{*}\mathfrak{a} = 0$ .

(2.7)

These equations will be the starting point of the subsequent analysis.

**e) The functions K and N**

Introduce the function $\kappa$ on $(0, \infty)$ whose square is at time t is the integral of $|\mathfrak{a}|^2$ on $\{t\} \times Y$:

$$\kappa^2(t) = \int_{\{t\}\times Y} |\mathfrak{a}|^2 \ .$$

(2.8)

The Nahm pole condition implies that $\kappa$ is approximately $\frac{\sqrt{3}}{2t}$ for small t.

For later reference, the derivative of $\kappa^2$ can be written as

- $\frac{1}{2}\frac{d}{dt}\kappa^2 = \int_{\{t\}\times Y} \langle \mathfrak{a}, \nabla_{A_t}\mathfrak{a}\rangle$ .
- $\frac{1}{2}\frac{d}{dt}\kappa^2 = \int_{\{t\}\times Y} \langle \mathfrak{a}, B_A\rangle - \int_{\{t\}\times Y} \langle \mathfrak{a}\wedge\mathfrak{a}\wedge\mathfrak{a}\rangle$ .

(2.9)

with the second version coming from the first via the middle bullet of (2.7). The second derivative of $\kappa^2$ can be written using (2.1) as

$$\frac{1}{2}\frac{d^2}{dt^2}\kappa^2 = \int_{\{t\}\times Y} (|\nabla_A\mathfrak{a}|^2 + |[\mathfrak{a}_*, \mathfrak{a}]|^2 + \langle \mathrm{Ric}, \langle \mathfrak{a}\otimes\mathfrak{a}\rangle\rangle) \ .$$

(2.10)

By way of a reminder: What is denoted here by $|[\mathfrak{a}_*, \mathfrak{a}]|^2$ is $\sum_{1\le i,j\le 3}|[\mathfrak{a}_i, \mathfrak{a}_j]|^2$. This is also $\langle \mathfrak{a}, [\mathfrak{a}_*, [\mathfrak{a}, \mathfrak{a}_*]]\rangle$ and also $2|\mathfrak{a}\wedge\mathfrak{a}|^2$.

The function $\kappa$ is never zero. Indeed, were $\kappa(t_*)$ to vanish for $t_* \in (0, \infty)$, then $\mathfrak{a}$ would vanish along the whole of $\{t_*\} \times Y$. In this event, there would be two solutions to (1.1) on $(0, \infty)$ that agreed on $[t_*, \infty)$ and disagreed on $(0, t_*)$. The first is the original Nahm pole solution and the second is defined for $t \in (0, t_*)$ by the rule whereby $(A, \mathfrak{a})$ at $t < t_*$ is $(A|_{2t_*-t}, -\mathfrak{a}|_{2t_*-t})$. A version of Aronzjain's [Ar] unique continuation principle could then be invoked to see that this nonsensical.

Because $\kappa > 0$, its time derivative can be written as

$$\frac{d}{dt}\kappa = -\frac{N}{t}\kappa$$

(2.11)

with N being a smooth function on $(0, \infty)$. By virtue of (2.10): If $s > t$, then

$$\mathrm{N}(t)\,\frac{\mathrm{K}(t)^2}{t} - \mathrm{N}(s)\,\frac{\mathrm{K}(s)^2}{s} = \int_{[t,s]\times Y} \left(|\nabla_A \mathfrak{a}|^2 + |[\mathfrak{a}_*, \mathfrak{a}]|^2 + \langle \mathrm{Ric}, \langle \mathfrak{a}\otimes\mathfrak{a}\rangle\rangle\right) . \tag{2.12}$$

The identity in (2.11) turns out to be a useful rewriting of (2.9).

## f) The tensor $\mathfrak{t}$

The traceless part of $\langle \mathfrak{a}\otimes\mathfrak{a}\rangle$ is denoted here by $\mathfrak{t}$:

$$\mathfrak{t} = \langle \mathfrak{a}\otimes\mathfrak{a}\rangle - \tfrac{1}{3}|\mathfrak{a}|^2\,\mathfrak{g} \tag{2.13}$$

with $\mathfrak{g}$ denoting the Riemannian metric. Thus, $\mathfrak{t}$ is at any given $t\in(0,\infty)$ is a symmetric, traceless section of $\otimes_2 T^*Y$.

The size of $|\mathfrak{t}|$ relative to $|\mathfrak{a}|^2$ effectively determines the size of the $*(\mathfrak{a}\wedge\mathfrak{a})$ term that appears in the second bullet of (2.7). To say more about this, define a function $\mathcal{Q}$ on $(0,\infty)\times Y$ by the formula

$$\mathcal{Q} = \tfrac{1}{3}*\langle \mathfrak{a}\wedge\mathfrak{a}\wedge\mathfrak{a}\rangle . \tag{2.14}$$

This is also $\frac{1}{3}\langle \mathfrak{a}, *(\mathfrak{a}\wedge\mathfrak{a})\rangle$.

Use the metric to view $\mathfrak{t}$ as a section of Hom(TY; TY) so as to define $\mathfrak{t}^2$ and $\mathfrak{t}^3$. Use the metric to define their traces also. The identity

$$|\mathfrak{a}|^6 = \tfrac{27}{4}\,\mathcal{Q}^2 + \tfrac{9}{2}|\mathfrak{a}|^2\,\mathrm{trace}(\mathfrak{t}^2) - 9\,\mathrm{trace}(\mathfrak{t}^3) \tag{2.15}$$

follows because $\mathcal{Q}^2 = 4\det(\langle \mathfrak{a}\otimes\mathfrak{a}\rangle)$ and $|\mathfrak{a}|^2$ is the trace of $\langle \mathfrak{a}\otimes\mathfrak{a}\rangle$ and $\langle \mathfrak{a}\otimes\mathfrak{a}\rangle = \mathfrak{t} + \frac{1}{3}|\mathfrak{a}|^2\mathfrak{g}$. (To prove (2.15): Just check that it holds when $\langle \mathfrak{a}\otimes\mathfrak{a}\rangle$ and hence $\mathfrak{t}$ are diagonal with respect to the metric $\mathfrak{g}$.)

As explained directly, (2.15) leads to upper and lower bounds for $\mathrm{K}^2$:

- $\mathrm{K}^2 \le \frac{3}{4^{1/3}} \int_{\{t\}\times Y} |\mathcal{Q}|^{2/3} + c_0 \int_{\{t\}\times Y} |\mathfrak{t}|$ .
- $\mathrm{K}^2 \ge \frac{3}{4^{1/3}} \int_{\{t\}\times Y} |\mathcal{Q}|^{2/3}$

$$\tag{2.16}$$

To prove the top bullet: Divide both sides of (2.15) by $|\mathfrak{a}|^2$ and then take the square root of both sides to see that $|\mathfrak{a}|^2 \le \frac{3\sqrt{3}}{2}\,\frac{1}{|\mathfrak{a}|}\,|\mathcal{Q}| + c_0|\mathfrak{t}|$. (Keep in mind here that $|\mathfrak{t}| \le c_0|\mathfrak{a}|^2$.) The preceding bound leads to the inequality $|\mathfrak{a}|^2 \le \frac{3}{4^{1/3}}|\mathcal{Q}|^{2/3} + c_0|\mathfrak{t}|$ because $|\mathcal{Q}|^{1/3} \le \frac{2^{1/3}}{\sqrt{3}}\,|\mathfrak{a}|$. The

top bullet of (2.16) then follows by integrating over Y. The lower bullet follows from the pointwise bound $|\mathcal{Q}| \le \frac{2}{3}|\mathfrak{a}|^3$.

The inequalities in (2.16) say in effect that the function of t given by the rule

$$t \to \frac{3}{4^{1/3}} \int_{\{t\}\times Y} |\mathcal{Q}|^{2/3} \tag{2.17}$$

is a proxy for $\kappa^2$ at times t when the integral of $|\mathfrak{t}|$ is small relative to $\kappa^2$. Looking ahead, this observation is useful because of the appearance of the integral of $\mathcal{Q}$ in the formula in the second bullet of (2.9) for the time derivative of $\kappa^2$.

An important point in what follows with regards to $\mathfrak{t}$ is that its derivative with respect to t sees only the $B_A$ term in the second bullet (2.7):

$$\frac{\partial}{\partial t}\mathfrak{t} = \langle \mathfrak{a} \otimes B_A\rangle + \langle B_A \otimes \mathfrak{a}\rangle - \frac{2}{3}\langle \mathfrak{a}, B_A\rangle \mathfrak{g} \ . \tag{2.18}$$

(Keep in mind that $\frac{\partial}{\partial t}$ is the product metric's covariant derivative in the t-direction.)

### g) The integral of $\langle \mathfrak{a}, B_A\rangle$

The function of t given by the integral of $\langle \mathfrak{a}, B_A\rangle$ on the slice $\{t\} \times Y$ also appears in the second bullet of (2.9). It also plays a central role in subsequent arguments (as it did in [MK] and [LT]). Of particular importance with regards to this function is the formula that follows for its derivative:

$$\frac{d}{dt}\int_{\{t\}\times Y} \langle \mathfrak{a}, B_A\rangle = \int_{\{t\}\times Y} \left(|B_A|^2 + \tfrac{1}{2}|E_A|^2 + \tfrac{1}{2}|\nabla_A^{\perp}\mathfrak{a}|^2 + \tfrac{1}{2}\langle \mathrm{Ric}, \langle \mathfrak{a}\otimes\mathfrak{a}\rangle\rangle\right) . \tag{2.19}$$

(Note that if the Ricci curvature is non-negative, then the integral of $\langle \mathfrak{a}, B_A\rangle$ is non-decreasing.) This equation is derived from (2.7) with part of the Bianchi identity:

$$\nabla_{At} B_A = *d_A E_A \ ; \tag{2.20}$$

and the identity

$$\int_{\{t\}\times Y} \left(|\nabla_A^{\perp}\mathfrak{a}|^2 + 2\langle B_A \wedge \mathfrak{a}\wedge\mathfrak{a}\rangle + \langle \mathrm{Ric}, \langle \mathfrak{a}\otimes\mathfrak{a}\rangle\rangle\right) = \int_{\{t\}\times Y} \left(|d_A\mathfrak{a}|^2 + |d_A *\mathfrak{a}|^2\right) \tag{2.21}$$

which holds whether or not $(A, \mathfrak{a})$ obeys the equations in (2.7). (The rest of the Bianchi identity says that $d_A{*}B_A = 0$.)

**h) The function $\mathfrak{cs}$**

The function $\mathfrak{cs}$ from (1.3) can be written using the notation in this section as

$$\mathfrak{cs}(t) = \int_{\{t\}\times Y} \langle \mathfrak{a}, B_A\rangle - \tfrac{1}{3} \int_{\{t\}\times Y} \langle \mathfrak{a}\wedge\mathfrak{a}\wedge\mathfrak{a}\rangle \ . \tag{2.22}$$

For future purposes, the most convenient depiction of its derivative is

$$\tfrac{d}{dt}\, \mathfrak{cs} = \tfrac{1}{2} \int_{\{t\}\times Y} (|\nabla_t \mathfrak{a}|^2 + |d_A\mathfrak{a}|^2 + |B_A - *(\mathfrak{a}\wedge\mathfrak{a})|^2 + |E_A|^2)\, ; \tag{2.23}$$

but by virtue of (2.7), the derivative is also equal to the $\{t\}\times Y$ integral of $|\nabla_t\mathfrak{a}|^2 + |E_A|^2$; and it is also equal to the $\{t\}\times Y$ integral of $|B_A - *(\mathfrak{a}\wedge\mathfrak{a})|^2 + |d_A\mathfrak{a}|^2$. Note that up to the factor of $\frac{1}{2}$, the $[t, \infty)\times Y$ integral of what is depicted on the right hand side of (2.23) is the integral in the second bullet of (1.2).

**j) Some basic analysis**

The purpose of this section is to provide the lemma that follows which supplies an analytic tool for use in later sections.

**Lemma 2.3**: *There exists* $\kappa > 1$ *with the following significance: Let* A *denote a connection on the principle bundle* $P \to Y$ *and let* $\mathfrak{q}$ *denote a section of* $T^*Y \otimes ad(P)$. *Supposing that* $\varepsilon > 0$, *then* $\int_Y |\nabla_A^{\perp}\mathfrak{q}|^2 + \kappa(1 + \frac{1}{\varepsilon^{3/2}})\,(\int_Y |\mathfrak{q}|)^2 \geq \frac{1}{\varepsilon}\int_Y |\mathfrak{q}|^2$ .

The proof of Lemma 2.3 invokes an analysis lemma which is stated first. To set the stage for this upcoming lemma, fix for the moment a positive number to be denoted by $x$ and let $\mathcal{K}_x$: $Y \times Y \to (0, \infty)$ denote the heat kernel for the Laplacian $d^\dagger d$ on Y evaluated at 'time' $x$. This function $\mathcal{K}_x$ can be written as follows: Let $\{\phi_k\}_{k=0,1,2,\ldots}$ denote an $L^2$ orthonormal basis of eigenvectors for the Laplacian $d^\dagger d$ on Y, labeled so that the corresponding eigenvalues are non-decreasing with increasing k. The eigenvalue of $\phi_k$ is denoted in what follows by $\lambda_k$; so $\phi_0 = 1$ and $\lambda_0 = 0$; and then $\lambda_{k-1} \leq \lambda_k$ for $k > 0$. With this notation in hand, here is $\mathcal{K}_\varepsilon$;

$$\mathcal{K}_x(p,q) = \sum_{k\in\{0,1,2,\cdots\}} e^{-\lambda_k x}\phi_k(p)\phi_k(q) \tag{2.24}$$

Supposing that $f$ is a continuous function on Y, define the function $f_x$ by

$$f_x(\mathrm{p}) = \int_Y \mathcal{K}_x(\mathrm{p},\cdot) f$$

(2.25)

This is a smeared version of $f$. What is left from $f$ is denoted by $f_{x\perp}$ which is $f - f_x$.

**Lemma 2.4**: *There exists* $\kappa > 1$ *with the following significance: Fix* $x > 0$ *and, given a continuous function,* $f$*, on* Y *define* $f_x$ and $f_{x\perp}$ *as instructed in the preceding paragraph.*

- $\int_Y |\mathrm{d}f|^2 \geq \frac{1}{\kappa\, x} \int_Y f_{x\perp}^2$ .
- *If* $f$ *is non-negative, then so is* $f_x$*.*
- $|f_x| \leq \kappa\,(1 + \frac{1}{x^{3/2}}) \int_Y |f|$ *(which is* $\kappa(1 + \frac{1}{x^{3/2}}) \int_Y f$ *if* $f$ *is non-negative.)*

***Proof of Lemma 2.4***: To prove the top bullet: Write function $f$ as $\sum \alpha_k \phi_k$ with $\alpha_k$ denoting the integral of $f\phi_k$ on Y. (Here and subsequently, $\sum$ signifies that a sum over the orthonormal basis labels $k \in \{0, 1, \ldots\}$ is taken.) The integral of $|\mathrm{d}f|^2$ can be written using the numbers $\{\alpha_k\}_{k=1,2,\ldots}$ as

$$\int_Y |\mathrm{d}f|^2 = \sum \lambda_k |\alpha_k|^2.$$

(2.26)

Meanwhile: Since $f_{x\perp} = \sum(1 - e^{-\frac{1}{E}\lambda_k})\, \alpha_k\, \phi_k$, the square of its $L^2$ norm is

$$\int_Y f_{x\perp}^2 = \sum(1 - e^{-\lambda_k x})^2 |\alpha_k|^2 .$$

(2.27)

The first bullet's assertion follows from (2.26) and (2.27) because of the following fact: If $u$ is a positive number (and in particular, $\lambda_k x$), the $u \geq c_0^{-1}(1 - e^{-u})^2$.

The second bullet follows directly from (2.25) because the function $\mathcal{K}_x$ is non-negative. The third bullet follows because $|\mathcal{K}_x| \leq c_0 \frac{1}{x^{3/2}}$ .

***Proof of Lemma 2.3***: Fix $x > 0$ and then the inequality $|\nabla_A{}^{\perp}\mathrm{q}| \geq |\,\mathrm{d}|\mathrm{q}|\,|$ and the top bullet of Lemma 2.4 using $f = |\mathrm{q}|$ to see that

$$\int_Y |\nabla_A \mathrm{q}|^2 \geq c_0^{-1} x \int_Y |\mathrm{q}|_{x\perp}^2 .$$

(2.28)

Meanwhile, the last two bullets of Lemma 2.4 lead to the bound

$$\int_Y |\mathfrak{q}|_x^2 \leq c_0 (1 + \tfrac{1}{x^{3/2}}) ( \int_Y |\mathfrak{q}| )^2 .$$

(2.29)

Since $|\mathfrak{q}|_x$ is $L^2$-orthogonal to $|\mathfrak{q}|_{x\perp}$, these last to bounds imply the bound in Lemma 2.3.

## 3. Nahm pole solutions where t is near zero

The purpose of this section is to establish some preliminary bounds for Nahm pole solutions where the coordinate t is nearly zero on $(0, \infty) \times Y$. These will be used in subsequent sections. (These bounds are also used in the appendix to show that the definition here of a Nahm pole solution is equivalent to the original definition in [MW].)

Keep in mind in what follows that Definition 1.1 is equivalent to the following assertion: Given $\varepsilon \in (0, 10^{-6}]$, there is some positive time (denoted by $t_\varepsilon$) such that

- $|\langle \mathfrak{a} \otimes \mathfrak{a} \rangle - \frac{1}{4t^2} \mathfrak{g}| < \frac{\varepsilon}{t^2}$ *on the whole of* $\{t\} \times Y$ *when* $t < t_\varepsilon$.
- $\int_0^{t_\varepsilon} (|B_A|^2 + |E_A|^2 + |\nabla_A^{\perp} \mathfrak{a}|^2) < \varepsilon$ .

(3.1)

Assume in what follows that any choice of $\varepsilon$ is less than $10^{-6}$.

### a) The isomorphism $\tau$ and endomorphism $\mathfrak{c}$

Let $\mathbb{V}_0$ and $\mathbb{V}_1$ denote two copies of $\mathbb{R}^3$ which are viewed here as inner product spaces (both have the Euclidean metric). Let $\mathrm{End}(\mathbb{V}_0, \mathbb{V}_1) = \mathbb{V}_1 \otimes \mathbb{V}_0{}^*$ denote the space of endomorphism from $\mathbb{V}_0$ to $\mathbb{V}_1$. Inside this sits the subspace of isometries; it is a smooth submanifold with 2 components. The group SO(3) acts freely on this submanifold either on the right via its action on $\mathbb{V}_0$ or on the left via its action on $\mathbb{V}_1$ and each component is the orbit via this action of any one of its components. This submanifold of isomorphisms has a bi-invariant tubular neighorhood with the following property: If $\mathfrak{q}$ is in the tubular neighborhood, then there is a unique isomorphism minimizing (over all isomorphisms) the function

$$\mathfrak{s} \longrightarrow \mathrm{trace}((\mathfrak{q} - \mathfrak{s})^T (\mathfrak{q} - \mathfrak{s}))$$

(3.2)

which is the square of the distance to $\mathfrak{q}$. (The exponent $(\cdot)^T$ indicates the transpose.) Note that this minimizer, $\tau_{\mathfrak{q}}$, will vary smoothly as $\mathfrak{q}$ varies over the tubular neighborhood.

Now suppose that $\mathfrak{q}$ is in the aforementioned tubular neighborhood and that $\tau_{\mathfrak{q}}$ is the unique isomormorphism that minimizes (3.2). The assertion that $\tau_{\mathfrak{q}}$ is a critical point of (3.2) (let alone a minimizer) is the assertion that $\mathfrak{q}$ can be written as

$$\mathfrak{q} = \tau_{\mathfrak{q}} + \mathfrak{z} \ , \tag{3.3}$$

where $\mathfrak{z}$ is such that $\mathfrak{z}^T\tau_{\mathfrak{q}}$ defines (via the metric) a *symmetric* element in $V_0{}^* \otimes V_0{}^*$. Note in this regard that $\tau_{\mathfrak{q}}{}^{-1} = \tau_{\mathfrak{q}}{}^T$ because $\tau_{\mathfrak{q}}$ is an endomorphism. (It is also the case that $\mathfrak{z}\tau_{\mathfrak{q}}{}^T$ is symmetric.)

In the context at hand, it is useful to view $\mathfrak{a}$ as a homomorphism from TY to ad(P). Viewed in this light, the Nahm pole condition in the top bullet of (3.1) implies that $\mathfrak{q} = -2t\mathfrak{a}$ is nearly an isomorphism for small t. Therefore (invoking the preceding paragraphs), it follows that $\mathfrak{a}$ can be written (for small t) as

$$\mathfrak{a} = -\tfrac{1}{2t}\,\tau + \mathfrak{c} \ . \tag{3.4}$$

where $\tau$ and $\mathfrak{c}$ are t-dependent, ad(P) valued 1-form on Y with the following properties:

- $\tau$ *is an isometry when viewed as a homomorphism from* TY *to* ad(P). *This is to say that* $\langle \tau \otimes \tau \rangle = \mathfrak{g}$.
- $\langle \mathfrak{c} \otimes \tau \rangle$ *is a symmetric section of* $T^*Y \otimes T^*Y$.
- $|\mathfrak{c}| \le \frac{100\,\varepsilon}{t}$ *where* $t < t_\varepsilon$ .

$$\tag{3.5}$$

By way of a parenthetical remark, the writing of $\mathfrak{a}$ as in (3.4) is equivariant with respect to the action of the group Aut(P) (which is the group of gauge transformations). In this regard: A gauge transformation acts on any given section of ad(P) by conjugation: If u denotes the gauge transformation and $\sigma$ denotes the section, then the action of $u$ sends $\sigma$ to $u\sigma u^{-1}$. Thus, it preserves the trace on ad(P) and thus the conditions in (3.5).

Write $\langle \mathfrak{c} \otimes \tau \rangle$ (which is a symmetric section of $\otimes^2 T^*Y$) as

$$\langle \mathfrak{c} \otimes \tau \rangle = \tfrac{1}{\sqrt{3}}\, c\, \mathfrak{g} + \mathfrak{c}^+ \tag{3.6}$$

where $\mathfrak{c}^+$ is the traceless part. The section $\mathfrak{t}$ from (2.13), a section of $\otimes^2 T^*Y$, can be written using $\mathfrak{c}^+$ and $c$:

$$\mathfrak{t} = -\tfrac{1}{t}\,\mathfrak{c}^+ + \tfrac{2}{\sqrt{3}}\, c\,\mathfrak{c}^+ + \mathfrak{c}^+ \bullet \mathfrak{c}^+ - \tfrac{1}{3}|\mathfrak{c}^+|^2 \mathfrak{g} \tag{3.7}$$

where the notation uses $\mathfrak{g}$ to denote the metric on Y and $x \bullet y$ to denote the section of $\otimes^2 T^*Y$ that is obtained from sections $x$ and $y$ of $\otimes^2 T^*Y$ by contracting (using the metric) the right most index of $x$ with the left most of $y$. (To be clear, fix an orthonormal frame $\{e^i\}$ for $T^*Y$. Then write $x$ with this frame as $x_{ij}\, e^i \otimes e^j$ and similarly $y$. Then $x \bullet y = x_{ik} y_{kj}\, e^i \otimes e^j$.) The third bullet of (3.5) and (3.7) lead to the inequalities

$$(1 - c_0\varepsilon)\, \tfrac{1}{t}\, |\mathfrak{c}^+| \le |\mathfrak{t}| \le (1 + c_0\varepsilon)\, \tfrac{1}{t}\, |\mathfrak{c}^+| \tag{3.8}$$

where $t < t_\varepsilon$. Thus, $\frac{1}{t}|\mathfrak{c}^+|$ can be used as a proxy for $|\mathfrak{t}|$ and vice-versa.

**b) Covariant derivatives of $\tau$ and $\mathfrak{c}$**

The covariant derivative of $\mathfrak{a}$ can be written using (3.4) in terms of covariant derivatives of $\tau$ and $\mathfrak{c}$. To this end, it is convenient to introduce (by way of notation) $\mathscr{b}$ to denote the t-dependent, $\otimes^2 T^*Y$ valued 1-form on $(0, t_\varepsilon] \times Y$ given by the rule

$$\mathscr{b} = \langle \tau \otimes \nabla_A \tau \rangle \,. \tag{3.9}$$

To be clear about the 1-form indices: These are associated with the covariant derivative, not with the indices from $\tau$. Thus, if $e$ denotes a tangent vector to $(0, t_\varepsilon] \times Y$ and $\nabla_{Ae}$ denotes the directional covariant derivative in the direction $e$. Then pairing of $\mathscr{b}$ with $e$ (which is a section of $\otimes^2 T^*Y$) is $\langle \tau \otimes \nabla_{Ae} \tau \rangle$. An important point is this: The $\otimes^2 T^*Y$ - valued 1-form $\mathscr{b}$ is *anti-symmetric* with regards to the $\otimes^2 T^*Y$ factor. This is because $\langle \tau \otimes \tau \rangle = \mathfrak{g}$ and the metric covariant derivative of $\mathfrak{g}$ is zero. (Note that $\mathscr{b}$ is also gauge invariant.)

In any event, the section $\mathscr{b}$ determines $\nabla_A \tau$ and (of course) vice-versa. In particular, $|\mathscr{b}| = |\nabla_A \tau|$. The section $\mathscr{b}$ also determines the connection A because any change $A \to A + \hat{a}$ with $\hat{a}$ being an ad(P) valued 1-form changes $\mathscr{b}$ to $\mathscr{b} + \langle \tau \otimes [\hat{a}, \tau] \rangle$ and $\langle \tau \otimes [\hat{a}, \tau] \rangle$ is zero if and only if $\hat{a} = 0$. (To see that this is so, fix an orthonormal frame for $T^*Y$ and write $\tau$ as $\tau_i e^i$. Then, the (i, j) component 1-form of $\langle \tau \otimes [\hat{a}, \tau] \rangle$ is $\langle \tau_i [\hat{a}, \tau_j] \rangle$ which is $\langle \hat{a}, [\tau_j, \tau_i] \rangle$ which is $2\varepsilon_{ijk} \langle \hat{a}\, \tau_k \rangle$ with $\varepsilon_{ijk}$ denoting the completely antisymmetric 3-tensor for $\mathbb{R}^3$ with $\varepsilon_{123} = 1$.)

One can view the anti-symmetric tensor valued 1-form $\mathscr{b}$ as a measure of the failure (if any) of $\tau$ to intertwine the covariant derivative $\nabla$, which acts on sections of TY, with the covariant derivative $\nabla_A$, which acts on sections of ad(P).

Now consider the covariant derivative of $\mathfrak{c}$: The important point to make in this regard is that the covariant derivative of $\mathfrak{c}$ is determined by $\langle\nabla_A\mathfrak{c}\otimes\tau\rangle$. In particular they have the same norm. Meanwhile, $\langle\nabla_A\mathfrak{c}\otimes\tau\rangle$ can be written using $\mathscr{b}$ as

$$\langle\nabla_A\mathfrak{c}\otimes\tau\rangle = \nabla\langle\mathfrak{c}\otimes\tau\rangle - \langle\mathfrak{c}\otimes\tau\rangle\bullet\mathscr{b} \tag{3.10}$$

where $\nabla$ denotes the Levi-Civita covariant derivative, and where $x\bullet\mathscr{b}$ for $x$ a section of $T^*Y\otimes T^*Y$ here again denotes contraction, using the metric, of the right most $T^*Y$ part of $x$ with the left most part of $\mathscr{b}$.

The identity in (3.10) leads to the following norm identities. (These are written using an orthonormal frame for $T^*Y$ to avoid misunderstandings. Repeated frame indices are implicitly summed in these formulas.)

- $|\nabla_{A_t}\mathfrak{a} - \frac{1}{2t^2}\tau|^2 = \frac{1}{4t^2}|\mathscr{b}_t|^2 + |\nabla_{A_t}\mathfrak{c}|^2 - \frac{1}{2t}\langle\mathfrak{c}\otimes\tau\rangle_{ik}(\mathscr{b}_t)_{kj}(\mathscr{b}_t)_{ji}$ .
- $|\nabla_A{}^{\perp}\mathfrak{a}|^2 = \frac{1}{4t^2}\langle\mathscr{b}_i, \mathscr{b}_i\rangle + |\nabla_A{}^{\perp}\mathfrak{c}|^2 - \frac{1}{2t}\langle\langle\mathfrak{c}\otimes\tau\rangle_{ik}(\mathscr{b}_m)_{kj}(\mathscr{b}_m)_{ji}$ .

(3.11)

Here $\mathscr{b}_t$ is the dt component of $\mathscr{b}$; and any given $i \in \{1, 2, 3\}$ version of $\mathscr{b}_i$ is the component of $\mathscr{b}$ along the basis vector $e^i$. (So $\mathscr{b}_t$ and each $\mathscr{b}_i$ is an anti-symmetric 2-tensor.) An key point is that (3.10) and the third bullet of (3.5) lead to the inequalities

- $|\nabla_{A_t}\mathfrak{a} - \frac{1}{2t^2}\tau|^2 \geq \frac{1}{8t^2}|\mathscr{b}_t|^2 + |\nabla_{A_t}\mathfrak{c}|^2$ ,
- $|\nabla_A{}^{\perp}\mathfrak{a}|^2 \geq \frac{1}{8t^2}|\mathscr{b}^{\perp}|^2 + |\nabla_A{}^{\perp}\mathfrak{c}|^2$ .

(3.12)

if $\varepsilon < c_0^{-1}$ and $t < t_\varepsilon$. Here and in what follows, $\mathscr{b}^{\perp}$ denotes the part of $\mathscr{b}$ that annihilates the tangents to the $(0,\infty)$ factor of $(0,\infty)\times Y$. (When written using the frame $\{e^i\}$, this is the $\mathscr{b}_i e^i$ part of $\mathscr{b}$.) Note in particular that the second bullet of (3.12) and the second bullet of (3.1) imply that $\mathscr{b}^{\perp}$ has limit 0 as $t\to 0$ (in a suitable sense). This says in effect that the $t\to 0$ limit of $\tau$ intertwines (at least in an $L^2$ sense) the metric covariant derivative along Y with the covariant derivative $\nabla_A{}^{\perp}$ on sections of ad(P). (The upcoming Lemma 3.5 asserts in part that $\mathscr{b}_t$ also has limit 0 as $t\to 0$ in a suitable sense.)

When $\langle\mathfrak{c}\otimes\tau\rangle$ is written as in (3.16), then (3.12) and the third bullet of (3.5) lead to bounds for the metric derivatives of $c$ and $\mathfrak{c}^+$:

- $2|\nabla_{A_t}\mathfrak{a} - \frac{1}{2t^2}\tau|^2 \geq \frac{1}{8t^2}|\mathscr{b}_t|^2 + |\frac{\partial}{\partial t}c|^2 + |\nabla_t\mathfrak{c}^+|^2$ .
- $2|\nabla_A{}^{\perp}\mathfrak{a}|^2 \geq \frac{1}{8t^2}|\mathscr{b}^{\perp}|^2 + |dc|^2 + |\nabla^{\perp}\mathfrak{c}^+|^2$ .

(3.13)

These hold where $t < t_\varepsilon$ when $\varepsilon < c_0^{-1}$. This last inequality is particularly useful because the metric covariant derivatives are (of course) independent of A.

The following lemma asserts asserts pointwise bounds for the norms of $\mathfrak{b}$, the covariant derivatives of $\mathfrak{c}$ and the curvatures $B_A$ and $E_A$. (They are not great bounds.)

**Lemma 3.1**: *There exists* $\kappa > 10^6$ *which is independent of* $(A, \mathfrak{a})$ *and has the following significance: Take* $\varepsilon < \kappa^{-1}$. *Then the bounds listed below hold where* $t < t_\varepsilon$.

- $|\nabla_A \mathfrak{a}| + |B_A| + |E_A| \le \kappa \frac{1}{t^2}$ .
- $|\mathfrak{b}| \le \kappa \frac{1}{t}$ .

***Proof of Lemma 3.1***: Fix $t \in (0, \frac{1}{2} t_\varepsilon)$ and a point in the slice $\{t\} \times Y$. Then, the second bullet of Proposition 2.1 is in play with $r = c_0^{-1} t$ because $|\mathfrak{a}| \le c_0 \frac{1}{t}$ (which follows from the top bullet of (3.1).) To say more: Supposing that $z > 100$, set $r = z^{-1} t$ and let B denote the radius r ball centered at the givien point in $\{t\} \times Y$. The value of $\hat{\kappa}$ (r) from Proposition 2.1 will be smaller than $c_0 \frac{1}{t}$ where as $\frac{1}{r}$ is $z \frac{1}{t}$ which is much larger than $\hat{\kappa}$ (r). In particular, if $z = c_0$, then the requirement for the second bullet of Proposition 2.1 will be met. That bullet gives the bounds in the top bullet of the lemma.

As for the second bullet: The bound $|\mathfrak{b}| \le c_0 \frac{1}{t}$ follows from the $|\nabla_A \mathfrak{a}|$ bound in the first bullet because of the inequality in (3.13).

**c) The integral of $|\mathfrak{c}^+|$ on $\{t\} \times Y$**

The purpose of this section is to first state and then prove a lemma about the $t \to 0$ limit of the integral of $|\mathfrak{c}^+|$ along $\{t\} \times Y$ (which is defined for small t).

**Lemma 3.2**: *There exists* $\kappa > 10^6$ *which is independent of* $(A, \mathfrak{a})$ *and has the following significance: Take* $\varepsilon < \kappa^{-1}$. *Suppose that* $\mathfrak{c}^+$ *is defined via (3.6) from a Nahm pole solution* $(A, \mathfrak{a})$. *If* $t > 0$ *but sufficiently small, then*

- $\int_{\{t\} \times Y} |\mathfrak{c}^+| \le \kappa t |\ln t|$ ,
- $\int_{\{t\} \times Y} |\mathfrak{c}^+|^2 \le \kappa t$ .

***Proof of Lemma 3.2***: The proof has five parts. The proof of the top bullet constitutes Parts 1-4 and the second bullet is proved in Part 5.

*Part 1*: The top bullet follows directly from (3.8) if

$$\int_{\{t\}\times Y} |\mathfrak{t}| \le c_0 |\ln t| \tag{3.14}$$

for $t < t_\varepsilon$ given that $\varepsilon < c_0^{-1}$. To prove (3.14), invoke (2.14) to see that

$$\int_{\{t\}\times Y} |\mathfrak{t}| \le \int_{\{t_\varepsilon\}\times Y} |\mathfrak{t}| + c_0 \int_{[t,t_\varepsilon]\times Y} |\mathfrak{a}||B_A| \;. \tag{3.15}$$

if $t < t_\varepsilon$ when $\varepsilon < c_0^{-1}$. Since $\varepsilon$ is fixed and independent of t, the integral of $|\mathfrak{t}|$ on $\{t_\varepsilon\} \times Y$ is independent of t; it is denoted henceforth as $z_\varepsilon$. (This could be huge, but no worries: The important point is that it is independent of t.) Keeping in mind that $|\mathfrak{a}|$ at any time s between t and $t_\varepsilon$ is at most $c_0 \frac{1}{s}$, the bound in (3.13) has the following implication: Fix a time $t_{\varepsilon 1} < t_\varepsilon$ and suppose in what follows that $t < t_{\varepsilon 1}$. Then

$$\int_{\{t\}\times Y} |\mathfrak{t}| \le z_{\varepsilon 1} + |\ln t| \cdot \sup_{s\in[t,\, t_{\varepsilon 1}]} \int_{\{s\}\times Y} |B_A| \;. \tag{3.16}$$

where $z_{\varepsilon 1}$, although possibly huge, is independent of t. Therefore, the bound in (3.14) follows from an s-independent bound by $c_0$ for the $L^1$ norm of $B_A$ on $\{s\} \times Y$ when s is less than some fixed $t_{\varepsilon 1}$.

*Part* 2: The identity in (2.3) for $F_A$ when written using $B_A$ and $E_A$ leads to an equation for $B_A$. Taking the inner product of the latter with $B_A$ leads in turn (where $t < t_\varepsilon$ when $\varepsilon < c_0^{-1}$) to a differential inequality for $|B_A|$ (where $|B_A| > 0$) that has form

$$-\frac{\partial^2}{\partial t^2} |B_A| + d^\dagger d|B_A| + (1 - c_0\varepsilon)\frac{2}{t^2} |B_A| \le c_0 \frac{1}{t^2} + c_0(|B_A|^2 + |E_A|^2 + |\nabla_A{}^{\perp}\mathfrak{a}|^2) \;. \tag{3.17}$$

(A key point to note is that the $B_A$ component of (2.3) lacks terms with $\nabla_{At}\mathfrak{a}$ on its right hand side. This is not the case for the $E_A$ part.)

Integrating (3.17) over the slice $\{t\} \times Y$ leads to the differential inequality

$$-\frac{d^2}{dt^2}\Big(\int_{\{t\}\times Y} |B_A|\Big) + \frac{7}{4t^2} \int_{\{t\}\times Y} |B_A| \le c_0 \frac{1}{t^2} + c_0 \int_{\{t\}\times Y} (|B_A|^2 + |E_A|^2 + |\nabla_A^{\perp}\mathfrak{a}|^2) \tag{3.18}$$

where $t < t_\varepsilon$ (supposing that $\varepsilon < c_0^{-1}$). If $\lambda$ is fixed and greater than $c_0$, then the function

$$f(t) = \int_{\{t\}\times Y} |B_A| - \lambda \tag{3.19}$$

obeys the equation

$$(-\frac{d^2}{dt^2} + \frac{7}{4t^2})f(t) \le c_0 \int_{\{t\}\times Y} (|B_A|^2 + |E_A|^2 + |\nabla_A^{\perp}\mathfrak{a}|^2)\,. \tag{3.20}$$

The plan for what follows is to use a Green's function for the differential operator $(-\frac{d^2}{dt^2} + \frac{7}{4t^2})$ to obtain a bound for $f$ from (3.20).

*Part* 3: This part of the proof constitutes a digression to introduce a Green's function for the operator $-\frac{d^2}{dt^2} + \frac{7}{4t^2}$ on $(0, t_\varepsilon)$. A suitable version with the 'pole' at a given $s \in (0, t_\varepsilon)$ is the function defined by the following rule:

- $G_s(t) = s^{-\sqrt{2}+1/2}\, t^{\sqrt{2}+1/2}$ *where* $s > t$ .
- $G_s(t) = s^{\sqrt{2}+1/2}\, t^{-\sqrt{2}+1/2}$ *where* $s < t$.

(3.21)

This Green's function obeys

$$(-\frac{d^2}{dt^2} + \frac{7}{4t^2})G_s(t) = \delta(t-s)\ . \tag{3.22}$$

Keep in mind when using $G_s$ that it is positive (for $t > 0$), and its maximum is s which is taken at the point $t = s$. Also,

$$G_s \sim t^{\sqrt{2}+1/2} \quad \textit{and} \quad |\tfrac{d}{dt} G_s| \sim t^{\sqrt{2}-1/2} \quad \textit{as} \ \ t \to 0. \tag{3.23}$$

The plan now is to multiply both sides of (3.20) by $G_s(t)$ and a suitable cut-off function and then integrate by parts to obtain a bound for the integral of $|B_A|$ over $\{s\}\times Y$. To define this cut-off function, reintroduce the bump function $\chi$ from Section 1f. (Remember that $\chi$ is a favorite non-increasing, smooth function on $\mathbb{R}$ that is equal to 1 on $(-\infty, \frac{1}{4}]$ and equal to 0 on $[\frac{3}{4}, 1)$.) Fix $s \in (0, \frac{1}{8} t_\varepsilon)$ and then fix $\delta > 0$ but much less than s. (A $\delta \to 0$ limit will be considered in any event.) With $\delta$ in hand, define $\chi_\delta$ to be the function on $(0,t_\varepsilon)$ that is given by the rule $t \to \chi(1 - t/\delta)$. This function is equal to 1 where $t > \delta$ and it is equal to 0 where $t < \frac{1}{4}\delta$. Meanwhile, let $\beta$ denote $\chi(2t/t_\varepsilon)$. The latter function equals 1 where $t < \frac{1}{8} t_\varepsilon$ and it is zero where $t > \frac{3}{8} t_\varepsilon$.

*Part 4*: Multiply both sides of (3.20) by $\chi_\delta \beta G_s$ and then integrate over $(0, t_\varepsilon)$. Integrate by parts and use (3.21), (3.22) and the bounds in (3.23) to see that

$$f(s) \le s\, z_{2\varepsilon} + c_0\, s^{-\sqrt{2}+1/2} \delta^{\sqrt{2}-3/2} \int_{\delta/4}^{\delta} |B_A| \ ,$$

(3.24)

where $z_{2\varepsilon}$ is s and $\delta$ independent. (The $s\, z_{2\varepsilon}$ contribution to (3.24) comes from the integration by parts terms with support near $t_\varepsilon$ and from the integral of $G_s$ times the term on far right on the right hand side of (3.18).)

The important point is that the integral of $|B_A|$ that appears in (3.24) is no greater than $c_0 \delta^{1/2}$ (because the integral of $|B_A|^2$ is bounded). Therefore, the term with $\delta$ on the right hand side of (3.24) is $\mathcal{O}(\delta^{\sqrt{2}-1})$ as $\delta \to 0$ (for fixed s); and this has limit zero as $\delta > 0$ because $\sqrt{2} \ge 1.4$ which is greater than 1. Therefore, after taking this limit, (3.24) says that $f(s) \le s\, z_{2\varepsilon}$ which is less than 1 when s is sufficiently small. As a consequence,

$$\int_{\{s\}\times Y} |B_A| \le c_0 \ \textit{ when s is sufficiently small.}$$

(3.25)

In particular, there is some postive $t_{\varepsilon 1}$ such that $\int_{\{s\}\times Y} |B_A| \le c_0$ holds when $s < t_{\varepsilon 1}$.

*Part 5*: This last part of the proof gives the argument for the lower bullet in Lemma 3.2. To start: It is sufficent by virtue of (3.8) to prove that

$$\Big( \int_{\{t\}\times Y} |\mathfrak{t}|^2 \Big)^{1/2} \le z_{\varepsilon 3} + \tfrac{1}{\sqrt{t}} \ .$$

(3.26)

when t sufficiently small, with $z_{\varepsilon 3}$ being independent of t. To do this, note first that

$$\Big| \tfrac{\partial}{\partial t} \Big( \int_{\{t\}\times Y} |\mathfrak{t}|^2 \Big)^{1/2} \Big| \le c_0 \Big( \int_{\{t\}\times Y} |\alpha|^2 |B_A|^2 \Big)^{1/2}$$

(3.27)

which is less than $c_0 \frac{1}{t} \Big( \int_{\{t\}\times Y} |B_A|^2 \Big)^{1/2}$. Granted the latter bound, fix $t_{\varepsilon 1} < t_\varepsilon$ for the moment and, supposing that $t < t_{\varepsilon 1}$, integrate (3.27) to see that

$$\Big( \int_{\{t\}\times Y} |\mathfrak{t}|^2 \Big)^{1/2} \le \Big( \int_{\{t_{\varepsilon 1}\}\times Y} |\mathfrak{t}|^2 \Big)^{1/2} + \tfrac{1}{\sqrt{t}} \Big( \int_{[0,t_{\varepsilon 1}]\times Y} |B_A|^2 \Big)^{1/2} .$$

(3.28)

Now, if $t_{\varepsilon 1}$ is sufficiently small (and so t is even smaller), then the integral of $|B_A|^2$ that appears in (3.28) will be less than 1. Meanwhile, the integral of $|\mathfrak{t}|^2$ on $\{t_{\varepsilon 1}\} \times Y$ is independent of t so it is some fixed (maybe huge) number which is $z_{\varepsilon 3}{}^2$ in (3.26).

**d) The integral of $|c|$ on $\{t\} \times Y$**

The purpose of this section is to first state and then prove a lemma about the integral of the function $c = \frac{1}{\sqrt{3}}\,\text{trace}(\langle \mathfrak{c} \otimes \tau \rangle)$ along $\{t\} \times Y$ (which is defined for small t).

**Lemma 3.3**: *There exists* $\kappa > 10^6$ *which is independent of* $(A, \mathfrak{a})$ *and has the following significance: Take* $\varepsilon < \kappa^{-1}$. *Suppose that* $c$ *is defined via (3.6) from a Nahm pole solution* $(A, \mathfrak{a})$. *Then, for small positive* t,

- $\int_{\{t\}\times Y} |c| \le \kappa t$ ,
- $\int_{\{t\}\times Y} c^2 \le \kappa t$ ,

***Proof of Lemma 3.3***: Taking the inner product of the middle equation in (2.7) with $\tau$ leads to an equation for $\frac{\partial}{\partial t} c$:

$$\frac{\partial}{\partial t} c + \frac{2}{t} c = \frac{1}{\sqrt{3}} \langle \tau, B_A \rangle + \frac{1}{\sqrt{3}} (2c^2 - |\mathfrak{c}^+|^2) \,. \tag{3.29}$$

(Keep in mind that $\mathfrak{b}_t = \langle \tau \otimes \nabla_{A_t} \tau \rangle$ is anti-symmetric so its trace is zero.) The equation in (3.29) leads in turn to a differential inequality for $|c|$ which can be integrated over any small t slice $\{t\} \times Y$ to obtain the inequality

$$\frac{d}{dt} \int_{\{t\}\times Y} |c| + \frac{2-\varepsilon}{t} \int_{\{t\}\times Y} |c| \le c_0 \Big( \int_{\{t\}\times Y} |B_A| + \int_{\{t\}\times Y} |\mathfrak{c}^+|^2 \Big) \tag{3.30}$$

where $t < t_\varepsilon$ (supposing that $\varepsilon < c_0{}^{-1}$). Multiplying both sides by $t^{2-\varepsilon}$ leads to the inequality

$$\frac{d}{dt} \Big( t^{2-\varepsilon} \int_{\{t\}\times Y} |c| \Big) \le c_0\, t^{2-\varepsilon} \tag{3.31}$$

for small t because of (3.25) and the second bullet of Lemma 3.2. Integrating this from 0 to any given small t leads to the bound in the top bullet of the lemma. (Keep in mind that $|c| \le 100\,\varepsilon\,\frac{1}{t}$ so there is nothing to worry about by way of a boundary contribution to the integration by parts at t = 0.)

To obtain the second bullet, multiply both sides of (3.29) by $c$ to obtain:

$$\frac{\partial}{\partial t} c^2 + \frac{4}{t} c^2 = \frac{1}{\sqrt{3}}\, c\langle \tau, B_A\rangle + \frac{1}{\sqrt{3}}\, c\,(2c^2 - |\mathfrak{c}^+|^2)\,. \tag{3.32}$$

The middle term on the right is at most $\frac{1}{2t} c^2 + c_0\, t\,|B_A|^2$ and the far right term on the right is at most $\frac{2\varepsilon}{t} c^2 + \frac{2\varepsilon}{t} |\mathfrak{c}^+|^2$. Thus, (3.32) leads to

$$\frac{\partial}{\partial t} c^2 + \frac{3}{t} c^2 \le c_0\, t\,|B_A|^2 + \frac{2\varepsilon}{t} |\mathfrak{c}^+|^2\,. \tag{3.33}$$

And, multipling both sides of this by $t^3$ and integrating over $\{t\} \times Y$ (and invoking the second bullet of Lemma 3.2) leads to the following when t is very small

$$\frac{\partial}{\partial t}\Big(t^3 \int_{\{t\}\times Y} c^2\Big) \le c_0\, t^4 \int_{\{t\}\times Y} |B_A|^2 + c_0\,\varepsilon\, t^3\,. \tag{3.34}$$

Integrate (3.34) form 0 to t and integrate by parts. The resulting inequality for small t implies the inequality in the second bullet of Lemma 3.3. Remember in this regard that the integral of $|B_A|^2$ from 0 to t will be less than 1 if t is small enough.

**e) The integrals of $|\nabla_{A_t}\mathfrak{c}|^2$ and $\frac{1}{t^2}|\mathfrak{c}|^2$ on $(0, t_\varepsilon] \times Y$.**

The lemma that follows asserts these $(0, t_\varepsilon] \times Y$ integrals are finite.

**Lemma 3.4**: *There exists $\kappa > 10^6$ which is independent of $(A, \mathfrak{a})$ and has the following significance: Fix $\varepsilon < \kappa^{-1}$. Define $\mathfrak{c}$ via (3.4) from a given Nahm pole solution $(A, \mathfrak{a})$. The integrals $\int_{[0,t_\varepsilon]\times Y} \frac{1}{s^2}|\mathfrak{c}|^2$ and $\int_{[0,t_\varepsilon]\times Y} |\nabla_{A_t}\mathfrak{c}|^2$, and $\int_{[0,t_\varepsilon]\times Y} |\frac{\partial}{\partial t} c|^2$ and $\int_{[0,t_\varepsilon]\times Y} |\nabla_t \mathfrak{c}^+|^2$ are finite*

***Proof of Lemma 3.4***: The proof has three parts. The first considers the integral of $\frac{1}{t^2}|\mathfrak{c}^+|^2$, the second considers the integral of $\frac{1}{t^2} c^2$. The third considers the integrals of $|\nabla_{A_t}\mathfrak{c}|^2$ and $|\nabla_t \mathfrak{c}^+|^2$ and $|\frac{\partial}{\partial t} c|^2$.

*Part 1*: To prove that the $(0, t_\varepsilon] \times Y$ of $\frac{1}{t^2}|\mathfrak{c}^+|^2$ is finite, it is sufficient to consider the integral on $(0, t] \times Y$ for any positive, small t. Start the proof by using the middle

bullet in (2.7) and (3.10) to obtain the following differential equation for the component $\mathfrak{c}^+$ from (3.6):

$$\nabla_{A_t}\mathfrak{c}^+ - \tfrac{1}{t}\,\mathfrak{c}^+ = \langle \tau \otimes B_A \rangle^+ - \langle \tau \otimes *(\mathfrak{c} \wedge \mathfrak{c}) \rangle^+ \,, \tag{3.35}$$

where $\langle \cdot \rangle^+$ indicates the symmetric, traceless part of the relevant section of $\otimes^2 T^*Y$. Taking the inner product of both sides of the preceding equation with $\frac{1}{t}\,\mathfrak{c}^+$ leads to the differential inequality

$$\tfrac{1}{2t}\,\tfrac{\partial}{\partial t}|\mathfrak{c}^+|^2 - \tfrac{1}{t^2}|\mathfrak{c}^+|^2 \geq -c_0\big(\tfrac{1}{t}|\mathfrak{c}_+||B_A| + \tfrac{1}{t}|\mathfrak{c}||\mathfrak{c}^+|^2\big) \tag{3.36}$$

Keeping in mind that $|\mathfrak{c}| \leq 100\,\varepsilon\,\frac{1}{t}$, this inequality implies the following one for small t:

$$\tfrac{1}{2t}\,\tfrac{\partial}{\partial t}|\mathfrak{c}^+|^2 - \tfrac{3}{4t^2}|\mathfrak{c}^+|^2 \geq -c_0|B_A|^2 \ . \tag{3.37}$$

Fix $t_{1\varepsilon} < t_\varepsilon$ so that (3.37) holds for $t < t_{1\varepsilon}$. Then, integrating (3.37) from t to $t_{1\varepsilon}$ and integrating by parts gives the bound

$$\int_{[t,t_{1\varepsilon}]\times Y} \tfrac{1}{4s^2}|\mathfrak{c}^+|^2 + \tfrac{1}{2t}\int_{\{t\}\times Y} |\mathfrak{c}^+|^2 \leq c_0 \int_{(0,t_{1\varepsilon}]\times Y} |B_A|^2 + \tfrac{1}{2t_{1\varepsilon}} \int_{\{t_{1\varepsilon}\}\times Y} |\mathfrak{c}^+|^2 \tag{3.38}$$

Taking $t \to 0$ with $t_{1\varepsilon}$ fixed proves that $\frac{1}{t^2}|\mathfrak{c}^+|^2$ is integrable.

*Part 2*: This part proves that the $(0,t_\varepsilon]\times Y$ integral of $\frac{1}{t^2}\,c^2$ is finite. To this end: Multiply both sides of (3.33) by $\frac{1}{t}$ to obtain an inequality of the form

$$\tfrac{1}{t}\,\tfrac{\partial}{\partial t}c^2 + \tfrac{3}{t^2}c^2 \leq c_0|B_A|^2 + c_0\tfrac{2\varepsilon}{t^2}|\mathfrak{c}^+|^2 \tag{3.39}$$

which holds when t is small. Fix $\delta > 0$ but much less than t (the $\delta \to 0$ limit will be taken) and integrate (3.39) from $\delta$ to t (with it understood that $t < t_\varepsilon$). Then integrate by parts to obtain the inequality

$$\int_{[\delta,t]\times Y} \tfrac{4}{s^2}c^2 \leq c_0\,\tfrac{1}{\delta}\int_{\{\delta\}\times Y} c^2 + c_0\int_{[0,t]\times Y} |B_A|^2 + \int_{[0,t]\times Y} \tfrac{1}{s^2}|\mathfrak{c}|^2 \tag{3.40}$$

The left hand side of (3.40) is bounded by $c_0$ no matter how small $\delta$ is because of the second bullet in Lemma 3.3 and because of what was just said in Part 1 of this proof.

Therefore, the right hand side of (3.40) is bounded independently of $\delta$; and so the $\delta \to 0$ limit on the left hand side of (3.40) exists (invoke the dominated convergence theorem).

*Part* 3: This last part considers the integral of $|\nabla_{A_t}\mathfrak{c}|^2$. The middle equation in (2.7) when written using $\mathfrak{c}$ leads to the bound

$$|\nabla_{A_t}\mathfrak{a} - \tfrac{1}{2t^2}\,\tau|^2 \le c_0\,|B_A|^2 + c_0(\tfrac{1}{t^2}\,|\mathfrak{c}|^2 + |\mathfrak{c}|^4)\ . \tag{3.41}$$

Since $|\mathfrak{c}| \le 100\varepsilon\,\frac{1}{t}$ when $t < t_\varepsilon$, the right hand side of this is at most $c_0|B_A|^2 + c_0\,\frac{1}{t^2}\,|\mathfrak{c}|^2$. Since the latter function is integrable on $(0, t_\varepsilon] \times Y$, so is the left hand side (invoke (3.1) and invoke Parts 1 and 2 of this proof). This implies (by virtue of (3.12)) that $|\nabla_{A_t}\mathfrak{c}|^2$ is integrable on $(0, t_\varepsilon] \times Y$, and that both $|\frac{\partial}{\partial t}\,c|^2$ and $|\nabla_t\mathfrak{c}^+|^2$ are also (by virtue of (3.13).

**f) Remarks about the curvature**

As mentioned previously, the anti-symmetric tensor valued 1-form $\mathscr{b} = \langle\tau \otimes \nabla_A\tau\rangle$ can be said to quantify the failure of $\tau$ to intertwine the Levi-Civita covariant derivative on sections of TY with the covariant derivative $\nabla_A$ on sections of ad(P). The formula in the upcoming (3.42) for the (4-dimensional) covariant exterior derivative (using the metric connection) of $\mathscr{b}$ suggest this. To set the notation, the formula uses an oriented, orthonormal frame $\{e^i\}_{i=1,2,3}$ for $T^*Y$ to give a basis $\{e^i \otimes e^j\}_{i,j=1,2,3}$ for $\otimes^2 T^*Y$. Thus, each (i, j) version of the formula that follows is an equality between 2-forms. What is denoted by $\mathcal{R}_{ijkm}$ in this formula is the Riemann curvature tensor for the metric on Y.

$$(D_\Gamma\mathscr{b} + \mathscr{b} \wedge \mathscr{b})_{ij} = 2\,\varepsilon_{ijk}\langle F_A\tau_k\rangle - \tfrac{1}{2}\,\mathcal{R}_{ijmn}\,e^m \wedge e^n\ . \tag{3.42}$$

Thus, $\langle F_A\tau_k\rangle = \frac{1}{4}\,\varepsilon_{ijk}(D_\Gamma\mathscr{b} + \mathscr{b} \wedge \mathscr{b})_{ij} + \frac{1}{8}\,\varepsilon_{ijk}\,\mathcal{R}_{ijmn}\,e^m \wedge e^n$. (What is denoted by $D_\Gamma$ is the $(0,\infty)\times Y$ covariant exterior derivative on $\otimes^2 T^*Y$ valued 1-forms that is defined using the metric connection.) The ad(P) valued 2-form $\frac{1}{4}\,\varepsilon_{ijk}(D_\Gamma\mathscr{b} + \mathscr{b} \wedge \mathscr{b})_{ij}\,\tau_k$ is, in effect, $\mathscr{b}$'s contribution to the curvature $F_A$.

Let $\mathcal{B}$ denote $*(\frac{1}{4}\,\varepsilon_{ijk}(d_\Gamma\mathscr{b}^\perp + \mathscr{b}\wedge \mathscr{b}^\perp)_{ij}\tau_k$ which is viewed in what follows as a t-dependent, ad(P) valued 1-form on Y. (The symbol $d_\Gamma$ denotes the Levi-Civita connection's covariant exterior derivative along the Y factor of $(0, \infty) \times Y$.) This $\mathcal{B}$ is, in effect, $\mathscr{b}$'s contribution to $B_A$ (which is the Hodge star of the curvature of A along the constant t slices in $(0,\infty)\times Y$.)

**g) The $(0, t_\varepsilon] \times Y$ integrals of $\frac{1}{t^2}|\mathfrak{b}|^2$ and $|\nabla \mathfrak{b}|^2$**

The following lemma makes a formal statement to the effect that functions $\frac{1}{t^2}|\mathfrak{b}|^2$ and $|\nabla\mathfrak{b}|^2$ are integrable on $[0, t_\varepsilon] \times Y$.

**Lemma 3.5**: *There exists $\kappa > 10^6$ which is independent of $(A, \mathfrak{a})$ and has the following significance: Take $\varepsilon < \kappa^{-1}$. If $\mathfrak{b}$ is defined as in (3.9), then $\int_{(0,t_\varepsilon]\times Y} (|\nabla\mathfrak{b}|^2 + \frac{1}{t^2}|\mathfrak{b}|^2)$ is finite.*

***Proof of Lemma 3.5***: The proof has three parts.

*Part 1*: The assertion that $\frac{1}{t^2}|\mathfrak{b}^\perp|^2$ is integrable on $(0, t_\varepsilon] \times Y$ follows from (3.1) and (3.12). The assertion that $\frac{1}{t^2}|\mathfrak{b}_t|^2$ is integrable on $(0, t_\varepsilon] \times Y$ follows from the top bullet of (3.12) given the fact noted in Part 3 of the proof of Lemma 3.4 that $|\nabla_{A_t}\mathfrak{a} - \frac{1}{2t^2}\tau|^2$ is integrable on the same domain.

*Part 2*: By virtue of (3.42) and the second bullet of Lemma 3.1, the Levi-Civita's covariant exterior derivative of $\mathfrak{b}$ enjoys the norm bound

$$|D_\Gamma\mathfrak{b}| \le c_0(1 + |F_A| + \tfrac{1}{t}|\mathfrak{b}|) . \tag{3.43}$$

As explained next, the exterior covariant derivative of $*\mathfrak{b}$ obeys a similar bound:

$$|D_\Gamma *_X \mathfrak{b}| \le c_0(1 + \tfrac{1}{t}|\mathfrak{b}| + \tfrac{1}{t}|c|) . \tag{3.44}$$

(In this context $*_X$ is the Hodge star operator for $(0, \infty) \times Y$.) To see where (3.44) comes from, rewrite (2.1) using $\tau$ where $t < t_\varepsilon$ as

$$\langle \tau \otimes \nabla_A{}^\dagger\nabla_A\mathfrak{a}\rangle + \langle \tau \otimes [\mathfrak{a}_*[\mathfrak{a}, \mathfrak{a}_*]]\rangle + \mathrm{Ric}(\langle \tau \otimes \mathfrak{a}\rangle) = 0 \tag{3.45}$$

so that each term on the left hand side is a t-dependent 1-form on Y with values in the vector bundle $\otimes_2 T^*Y$. Each term on the left hand side can be written as a sum of a symmetric section of $\otimes_2 T^*Y$ and an anti-symmetric section. The symmetric part and the anti-symmetric part must both vanish separately. Meanwhile, the anti-symmetric part can be written schematically

$$\tfrac{1}{2t}(1 - t\mathbb{M}) *_X D_\Gamma(*_X\mathfrak{b}) + \mathbb{L} = 0 \tag{3.46}$$

where $\mathbb{M}$ is an endomorphism that obeys $|\mathbb{M}| \le c_0|\mathfrak{c}|$ and where the norm of $\mathbb{L}$ is bounded by $c_0(|\mathcal{b}|\,|\nabla_A\mathfrak{c}| + |\mathcal{b}|^2|\mathfrak{c}| + \frac{1}{t}|\mathfrak{c}|^2 + |\mathfrak{c}|^3)$. The key point is that are no terms with two Levi-Civita covariant derivatives acting on $\langle\tau\otimes\mathfrak{c}\rangle$ appear in the antisymmetric part of $\langle\tau\otimes\nabla_A{}^{\dagger}\nabla_A\mathfrak{a}\rangle$ because $\langle\tau\otimes\mathfrak{c}\rangle$ is a symmetric section of $\otimes_2 T^*Y$. The bound in (3.44) follows from (3.46) and the preceding bounds for $\mathbb{M}$ and $\mathbb{L}$ it $\varepsilon < c_0^{-1}$ (and where $t < \frac{1}{2}t_\varepsilon$) because $|\mathfrak{c}| \le \frac{100\varepsilon}{t}$ , and because of Lemma 3.1.

*Part 3*: Fix $\delta \in (0, \frac{1}{100}t_\varepsilon]$ and reintroduce the function $\chi_\delta$ from Part 3 of the proof of Lemma 3.2. Remember that this function is equal to 1 where $t > \delta$ and it is equal to 0 where $t < \frac{1}{4}\delta$. And reintroduce the function $\beta$ from this same Part 3 of Lemma 3.2's proof. The latter function equals 1 where $t < \frac{1}{8}t_\varepsilon$ and it is zero where $t > \frac{3}{8}t_\varepsilon$. Multiply both sides of (3.43) and both sides of (3.44) by $\chi_\delta\beta$; then square both sides and integrate over $(0, t_\varepsilon]\times Y$. Having done this, integration by parts to rearrange the derivatives on $\mathcal{b}$ leads to the following inequality:

$$\int_{[\delta,\frac{1}{4}t_\varepsilon]\times Y} |\nabla\mathcal{b}|^2 \le c_0 \int_{(0,t_\varepsilon)\times Y} (|F_A|^2 + |\nabla_A\mathfrak{c}|^2 + \tfrac{1}{t^2}|\mathfrak{c}|^2 + \tfrac{1}{t^2}|\mathcal{b}|^2)\ . \tag{3.47}$$

The assertion of Lemma 3.5 follows from this last inequality because the right hand side of (3.47) is finite and independent of $\delta$.

## h) The $\{t\}\times Y$ integral of $\langle\mathfrak{a},\mathcal{B}\rangle$

Remember that $\mathcal{B}$ denotes $*(\frac{1}{4}\varepsilon_{ijk}(d_\Gamma\mathcal{b}^\perp + \mathcal{b}^\perp\wedge\mathcal{b}^\perp)_{ij}\tau_k$ which is viewed as a t-dependent, ad(P) valued 1-form on Y. Note that the respective $\{t\}\times Y$ integrals of $\langle\mathfrak{a},\mathcal{B}\rangle$ and $\langle\mathfrak{a}\wedge\frac{1}{4}\varepsilon_{ijk}(D\mathcal{b} + \mathcal{b}\wedge\mathcal{b})_{ij}\tau_k\rangle$ are the same because $*\mathcal{B}$ on Y is the restriction to TY of $\frac{1}{4}\varepsilon_{ijk}(D\mathcal{b} + \mathcal{b}\wedge\mathcal{b})_{ij}\tau_k\rangle$.

**Lemma 3.6**: *There exists* $\kappa > 10^6$ *which is independent of* $(A, \mathfrak{a})$ *and has the following significance: Take* $\varepsilon < \kappa^{-1}$. *If* $t \in (0, t_\varepsilon]$ *is sufficiently small, then*

- $\displaystyle\int_{\{t\}\times Y} \langle\mathfrak{a},\mathcal{B}\rangle \le \int_{(0,t]\times Y} (|B_A|^2 + |E_A|^2 + |\nabla_A^\perp\mathfrak{a}|^2) + \kappa t$
- $\displaystyle\int_{\{t\}\times Y} \langle\mathfrak{a},\mathcal{B}\rangle \ge \int_{(0,t]\times Y} (|B_A|^2 + |E_A|^2 + |\nabla_A^\perp\mathfrak{a}|^2) - \kappa t$

*In particular,* $\displaystyle\lim_{t\to 0} \int_{\{t\}\times Y} \langle\mathfrak{a},\mathcal{B}\rangle = 0.$

***Proof of Lemma 3.6***: The proof has four parts.

*Part 1*: Write $\mathfrak{a}$ as (3.4) to decompose the integral of $\langle \mathfrak{a}, \mathcal{B}\rangle$ into two parts:

$$-\tfrac{1}{2t} \int_{\{t\}\times Y} \langle \tau \wedge \tfrac{1}{4}\varepsilon_{ijk}(d_\Gamma \mathscr{b}^\perp + \mathscr{b}^\perp \wedge \mathscr{b}^\perp)_{ij}\tau_k\rangle \ + \int_{\{t\}\times Y} \langle \mathfrak{c} \wedge \tfrac{1}{4}\varepsilon_{ijk}(d_\Gamma \mathscr{b}^\perp + \mathscr{b}^\perp \wedge \mathscr{b}^\perp)_{ij}\tau_k\rangle \ .$$

(3.48)

With regards to the left hand integral: It follows from the definitions of $\tau$ (and an integration by parts) that the $d_\Gamma \mathscr{b}^\perp$ part of the left hand integral in (3.48) is zero. The left hand integral in (3.48) is therefore bounded by

$$c_0 \tfrac{1}{t} \int_{\{t\}\times Y} |\mathscr{b}^\perp|^2$$

(3.49)

with $\mathscr{b}^\perp$ again denoting the part of $\mathscr{b}$ that annilates the tangent vectors to the $(0, \infty)$ factor of $(0,\infty) \times Y$. Meanwhile, the function on $(0,t_\varepsilon)\times Y$ given by $\frac{1}{t^2} |\mathscr{b}^\perp|^2$ is integrable (because of Lemma 3.5), and this implies the following: Given $\varepsilon > 0$ and a positive integer n, let $\mu_{\varepsilon,n}$ denote the fraction of the interval $[2^{-n-1}t_\varepsilon, 2^{-n}t_\varepsilon]$ where the number depicted in (3.49) is greater than $\varepsilon$. Then $\lim_{n\to\infty} \mu_{\varepsilon,n} = 0$.

*Part 2*: With regards to the right hand integral in (3.48): It is bounded by

$$c_0\Big(\int_{\{t\}\times Y} |\mathfrak{c}|^2\Big)^{1/2}\Big(1 + \int_{\{t\}\times Y} |B_A|^2\Big)^{1/2} \ .$$

(3.50)

(This bound follows from (3.42).) And, because of Lemmas 3.2 and 3.3, what is written (3.50) is bounded in turn by

$$c_0 \sqrt{t}\,\Big(1 + \int_{\{t\}\times Y} |B_A|^2\Big)^{1/2}.$$

(3.51)

This last expression defines a function on $(0, t_\varepsilon]$ which is denoted $\mathscr{h}$. Then, by virtue of the second bullet in (3.1)

$$\lim_{t\to 0} \tfrac{1}{t} \int_0^t \mathscr{h} \ = 0 \ .$$

(3.52)

This implies in particular the following: Given $\varepsilon > 0$ and a positive integer n, let $\nu_{\varepsilon,n}$ denote the fraction of of the interval $[2^{-n-1}t_\varepsilon, 2^{-n}t_\varepsilon]$ where $\mathscr{h} > \varepsilon$. Then $\lim_{n\to\infty} \nu_{\varepsilon,n} = 0$.

*Part 3*: Given $\varepsilon > 0$ and a positive integer n, let $\mathrm{M}_{\varepsilon,n}$ denote the fraction of the interval $[2^{-n-1}t_\varepsilon, 2^{-n}t_\varepsilon]$ where the norm of the $\{t\}\times Y$ integral of $\langle \mathfrak{a}, \mathcal{B}\rangle$ is greater than $\varepsilon$. Then $\lim_{n\to\infty} \mathrm{M}_{\varepsilon,n} = 0$ because of what is said in Parts 1 and 2.

The preceding assertion is almost what is said by Lemma 3.6, but not quite. As is explained next, Lemma 3.6's assertion does follows from the fact that $\lim_{n\to\infty} \mathrm{M}_{\varepsilon,n} = 0$ using extra input from the identity in (2.19). To see why, fix times $t_0 < t_1$, both from $(0,t_\varepsilon]$ and integrate (2.19) on $[t_0,t_1]\times Y$ to obtain the following identity:

$$\int_{\{t_1\}\times Y} \langle \mathfrak{a}, B_A\rangle - \int_{\{t_0\}\times Y} \langle \mathfrak{a}, B_A\rangle = \int_{[t_0,t_1]\times Y} (|B_A|^2 + \tfrac{1}{2}|E_A|^2 + \tfrac{1}{2}|\nabla_A^\perp \mathfrak{a}|^2) + \tfrac{1}{2}\int_{[t_0,t_1]\times Y} \langle \mathrm{Ric}, \langle \mathfrak{a}\otimes\mathfrak{a}\rangle\rangle \ . \tag{3.53}$$

To exploit this, write $\langle \mathfrak{a}, B_A\rangle$ using (3.42) as

$$\langle \mathfrak{a}, B_A\rangle = \langle \mathfrak{a}, \mathcal{B}\rangle + \langle \mathfrak{a}, B_\Gamma\rangle, \tag{3.54}$$

where the notation has $B_\Gamma$ denoting the contribution to $B_A$ from the Riemann curvature term in (3.42).

As explained in Part 4, using (3.54) for $\langle \mathfrak{a}, B_A\rangle$ on the left hand side of (3.53) leads to the following observation:

*If* $t_0 < t_1$ *are from* $(0, t_\varepsilon]$, *then the difference between the respective* $\{t_1\} \times Y$ *and* $\{t_0\} \times Y$ *integrals of* $\langle \mathfrak{a}, \mathcal{B}\rangle$ *differs by at most* $c_0 t_1$ *from*

$$\int_{[t_0,t_1]\times Y} (|B_A|^2 + |E_A|^2 + |\nabla_A^\perp \mathfrak{a}|^2) \ . \tag{3.55}$$

Lemma 3.6 follows from the preceding because: Given $\varepsilon > 0$, the time $t_0$ in (3.55) can be chosen less than any give positive time so that the $\{t_0\}\times Y$ integral of $\langle \mathfrak{a}, \mathcal{B}\rangle$ is less than $\varepsilon$.

*Part 4*: To prove the claim in (3.55), consider first the difference between the respect $\{t_0\}\times Y$ and $\{t_1\}\times Y$ integrals of $\langle \mathfrak{a}, B_\Gamma\rangle$ that appear on the left hand side of (3.53) when $\langle \mathfrak{a}, B_A\rangle$ is written as in (3.54). The computation for this difference uses the following formula for $B_\Gamma$:

$$B_\Gamma = \tfrac{1}{4} R\tau_k e^k - \tfrac{1}{2} \mathrm{Ric}_{ik}\tau_i e^k \ , \tag{3.56}$$

with R denoting the metric's scalar curvature. This formula can now be used to write the difference between the respect $\{t_0\}\times Y$ and $\{t_1\}\times Y$ integrals of $\langle \mathfrak{a}, B_\Gamma\rangle$ as a sum of the three terms

- $\frac{1}{8}\left(\frac{1}{t_0} - \frac{1}{t_1}\right)\int_Y R$
- $-\frac{1}{2}\int_{\{t_1\}\times Y}\langle \mathrm{Ric}, \mathfrak{c}^+\rangle + \frac{1}{2}\int_{\{t_0\}\times Y}\langle \mathrm{Ric}, \mathfrak{c}^+\rangle$ .
- $\frac{1}{2\sqrt{3}}\int_{\{t_1\}\times Y} c\mathrm{R} - \frac{1}{2\sqrt{3}}\int_{\{t_0\}\times Y} c\mathrm{R}$ .

(3.57)

(These correspond to the respective $-\frac{1}{2t}\tau$ and $c$ and $\mathfrak{c}^+$ contributions to $\mathfrak{a}$.)

Now, the $c$ integrals in the third bullet of (3.57) are bounded by $c_0t_1$ and $c_0t_0$ respectively (see Lemma 3.3), so they can be ignored for the purposes of (3.55) and Lemma 3.6. Meanwhile, the $\mathfrak{c}^+$ integrals in the middle bullet of (3.57) are bounded apriori only by $c_0t_1|\ln t_1|$ and $c_0t_0|\ln t_0|$ respectively (see Lemma 3.2), so they could be ignored if one can live with the error in (3.55) being $c_0t_1|\ln t_1|$ and a corresponding $\kappa\, t|\ln t|$ error term in Lemma 3.6. In fact, for applications to come, the extra $|\ln t|$ factor is permissable, but even so, the arguments that follow will explain why there is no such $|\ln t|$ factor in (3.55) and in Lemma 3.6. Of course, the norm of the expression in the top bullet of (3.57) is huge when $t_0$ is small (unless the integral of R is zero). As explained next, up to a term with norm bounded by $c_0t_1$, the expressions in the top bullet and second bullet in (3.57) are precisely accounted for by the $\langle \mathrm{Ric}, \langle\mathfrak{a}\otimes\mathfrak{a}\rangle\rangle$ integral on the right hand side of (3.53). The rest of the right hand side of (3.53) (the $|B_A|^2 + \frac{1}{2}|E_A|^2 + \frac{1}{2}|\nabla_A{}^{\perp}\mathfrak{a}|^2$ integral) accounts for the $|B_A|^2 + |E_A|^2 + |\nabla_A{}^{\perp}\mathfrak{a}|^2$ integral in (3.55).

To see that the $[t_0,t_1]\times Y$ integral of $\langle \mathrm{Ric}, \langle\mathfrak{a}\otimes\mathfrak{a}\rangle\rangle$ account for the terms in the first and second bullets of (3.57) (up to a $c_0t_1$ error), first use (3.6) to write

$$\langle\mathfrak{a}\otimes\mathfrak{a}\rangle = \frac{1}{4t^2}\left(1 - \frac{2}{\sqrt{3}}\,t\,c\right)^2\mathfrak{g} - \frac{1}{t}\,\mathfrak{c}^+ + \mathfrak{c}^+\otimes\mathfrak{c}^+ \ . \tag{3.58}$$

Now, Lemmas 3.2 and 3.3 imply that the contributions from term linear in $c$ and from the quadratic terms in $c$ and $\mathfrak{c}^+$ are bounded by $c_0t_1$ when $t_1$ is small. With this understood, then up to an $\mathcal{O}(t_1)$ correction,

$$\frac{1}{2}\int_{[t_0,t_1]\times Y}\langle \mathrm{Ric}, \langle\mathfrak{a}\otimes\mathfrak{a}\rangle\rangle = \frac{1}{8}\left(\frac{1}{t_0} - \frac{1}{t_1}\right)\int_Y \mathrm{R} - \frac{1}{2}\int_{[t_0,t_1]\times Y}\frac{1}{t}\langle \mathrm{Ric}, \mathfrak{c}^+\rangle \ . \tag{3.59}$$

The left most term on the right hand side of (3.59) is exactly what appears in the top bullet of (3.57). As for the right most term in (3.59): Use (3.35) to replace $\frac{1}{t}\langle \mathrm{Ric}, \mathfrak{c}^+\rangle$ in the integrand by $\frac{\partial}{\partial t}\langle \mathrm{Ric}, \mathfrak{c}^+\rangle$ plus expressions that are bounded by $c_0(|\mathfrak{c}|^2 + |\mathscr{b}_t|^2 + |B_A|)$. Having done this, then integrate by parts to see that contribution from $\frac{\partial}{\partial t}\langle \mathrm{Ric}, \mathfrak{c}^+\rangle$ precisely accounts for the expression in the second bullet of (3.57). Meanwhile, the norm of the contribution from the expressions that are bounded by $c_0(|\mathfrak{c}|^2 + |\mathfrak{b}_t|^2 + |B_A|)$ is at most $c_0 t_1$ (invoke Lemmas 3.4 and 3.5 to bound the respective $[t_0, t_1] \times Y$ integrals of $|\mathfrak{c}|^2$ and $|\mathfrak{b}_t|^2$; use (3.25) to bound the $[t_0, t_1] \times Y$ integral of $|B_A|$.)

### i) A priori $C^k$ bounds

Fix $t > 0$. The equations in (2.7) when written in terms of $\hat{\mathfrak{c}} = \langle \mathfrak{c} \otimes \tau \rangle$ and $\mathfrak{b}$ with the addition of (3.44) are uniformly elliptic equations for the pair $(\hat{\mathfrak{c}}, \mathfrak{b})$ on an open set with compact closure in $(0, t_\varepsilon) \times Y$ when $\varepsilon < c_0^{-1}$. Because of this, they lead (using standard arguments) to a priori $C^k$ bounds (for any given non-negative integer k) on the components of the tensors $(\hat{\mathfrak{c}}, \mathfrak{b})$ on any domain of the form $(t, 2t) \times Y$ for $t \in (0, \frac{1}{2} t_\varepsilon)$. These bounds will depend on k and t, but not on the chosen Nahm pole solution. The following lemma makes a formal statement to this effect.

**Lemma 3.7**: *There exists* $\kappa > 10^6$ *with the following significance*: *Supposing that* $(A, \mathfrak{a})$ *is a Nahm pole solution, fix* $\varepsilon < \frac{1}{\kappa}$ *so as to specify the time* $t_\varepsilon$ *and to define the pair* $(\hat{\mathfrak{c}}, \mathfrak{b})$ *on* $(0, t_\varepsilon) \times Y$. *Fix a non-negative integer* k *and a time* $t \in (0, \frac{1}{2} t_\varepsilon)$. *The* $C^k$ *norm of* $(\hat{\mathfrak{c}}, \mathfrak{b})$ *on* $[t, 2t] \times Y$ *is bounded by number that depends on* t *and* k *but not on* $(A, \mathfrak{a})$.

As noted, the proof is a standard application of elliptic regularity arguments using the equations in (2.7) and the equation in (3.44).

## 4. The frame change for the (relatively) large t analysis

The isomorphism $\tau$ from TY to ad(P) that is defined in Section 3a and appears in (3.4) can only be defined at values of t where $t^2 \langle \mathfrak{a} \otimes \mathfrak{a} \rangle$ is close to the Riemannian metric on Y. This constraint is unfortunate. A second isomorphism from TY to ad(P) is defined momentarily for use where t is not so very small. To this end, fix $t^- > 0$ so that $\tau$ is defined for $t < t^-$. Choose it so that it is less than $10^{-6} t_\varepsilon^{\,4}$ for $\varepsilon = 10^{-6}$. It therefore obeys (3.5) with this version of $\varepsilon$. Make $t^-$ smaller if necessary so that the conclusions of the Lemmas 3.1–3.7 hold for $t < t^-$. And, make it even smaller (if needed) so that the absolute values of the integrals considered by Lemma 3.6 are less than 1 when $t < t^-$. One additional constraint on $t^-$ is needed to invoke the upcoming Lemma 4.1.

### a) A second isomorphism from TY to ad(P).

For times $t \geq t^-$, define an isomorphism from TY to ad(P) (to be denoted by $\sigma$) is defined by the following two rules:

- $\sigma = \tau$ at $t = t^-$ .
- $\nabla_{At}\sigma = 0$ for $t > t^-$ .

$$(4.1)$$

Thus, $\sigma$ at any $(t, x) \in [t^-, \infty) \times Y$ is obtained from $\tau$ at $(t^-,x)$ by parallel transport along the path $[t^-, t] \times \{x\}$.

This isomorphism $\sigma$ can also be viewed as a t-dependent, ad(P) valued 1-form on Y. View it in this light to write $\mathfrak{a}$ for $t \geq t^-$ as:

$$\mathfrak{a} = -\tfrac{1}{2t}\sigma + \mathfrak{c} \tag{4.2}$$

with $\mathfrak{c}$ being an ad(P) valued 1-form. (This $\mathfrak{c}$ and that in (3.4) are identical at $t = t^-$. The version in (3.4) will be considered for only where $t < t^-$ so there should be no confusion as to which is which.) Note that the section $\langle \mathfrak{c} \otimes \sigma\rangle$ of $\otimes^2 T^*Y$ need not be symmetric. Even so, there is still an analog of (3.6):

$$\langle \mathfrak{c} \otimes \sigma\rangle = \tfrac{1}{\sqrt{3}}\, c\, \mathfrak{g} + \mathfrak{c}^+ + \mathfrak{c}^- , \tag{4.3}$$

where $\mathfrak{c}^+$ is traceless and symmetric, and where $\mathfrak{c}^-$ is antisymmetric. The component $\mathfrak{c}^-$ vanishes at $t = t^-$, where as $c$ and $\mathfrak{c}^+$ are equal to their namesakes in (3.6) at $t = t^-$. The norm squared of $\mathfrak{a}$ and the tensor $\mathfrak{t}$ when written using $c$, $\mathfrak{c}^+$ and $\mathfrak{c}^-$ are

- $|\mathfrak{a}|^2 = \frac{3}{4t^2} - \frac{\sqrt{3}}{t} c + c^2 + |\mathfrak{c}^+|^2 + |\mathfrak{c}^-|^2$ .
- $\mathfrak{t} = -\frac{1}{t}\mathfrak{c}^+ + \frac{2}{\sqrt{3}} c\mathfrak{c}^+ + \mathfrak{c}^+\cdot\mathfrak{c}^+ + \mathfrak{c}^-\cdot\mathfrak{c}^+ - \mathfrak{c}^+\cdot\mathfrak{c}^- - \mathfrak{c}^-\cdot\mathfrak{c}^- - \frac{1}{3}(|\mathfrak{c}^+|^2+|\mathfrak{c}^-|^2)\mathfrak{g}$ .

$$(4.4)$$

Remember that $\mathfrak{t}$ is the traceless part of $\langle \mathfrak{a} \otimes \mathfrak{a}\rangle$.

**b) The derivative of $\sigma$**

What is denoted for $t \geq t^-$ as $\mathcal{b}$ is (henceforth) the anti-symmetric tensor valued 1-form given by

$$\mathcal{b} = \langle \sigma \otimes \nabla_A \sigma\rangle . \tag{4.5}$$

Note that the component $\mathcal{b}_t$ proportional to dt is zero. The other components are equal at $t = t^-$ to the $\mathcal{b}^\perp$ part of what is denoted by $\mathcal{b}$ in (3.9) at $t = t^-$. The analog of (3.10) is

$$\langle \nabla_A \mathfrak{c} \otimes \sigma\rangle = \nabla\langle \mathfrak{c} \otimes \sigma\rangle - \langle \mathfrak{c} \otimes \sigma\rangle\cdot \mathcal{b} . \tag{4.6}$$

It proves useful in what follows to denote by $\mathcal{B}$ the t-dependent, ad(P) valued 1-form on Y that is defined by the rule

$$\mathcal{B} = \tfrac{1}{4} *(d_\Gamma \mathscr{b} + \mathscr{b} \wedge \mathscr{b})_{ij}\, \varepsilon_{ijk} \sigma_k \tag{4.7}$$

where the notation is as follows: What is $d_\Gamma$ denotes the metric connection's exterior, covariant derivative along Y. The subsript $(\cdot)_{ij}$ indicates the $\otimes^2$ component along the basis 2-tensor $e^i \otimes e^j$, so each $(\cdot)_{ij}$ is a 2-form on Y (the Hodge star makes this 2-form into a 1-form.) Meanwhile, $\sigma_k$ is the component of $\sigma$ along $e^k$ (so it is a section of ad(P).) The $\sigma$ analog of (3.42) can be written using the 1-forms $\{\mathcal{B}_k\}_{k=1,2,3}$ and the curvature of Y as

- $B_A = \mathcal{B} - \tfrac{1}{2} \mathrm{Ric}_{km} e^m \sigma^k + \tfrac{1}{4} R e^k \sigma^k$ ,
- $E_A = \tfrac{1}{4} \varepsilon_{ijk} \tfrac{\partial}{\partial t} \mathscr{b}_{ij}\, \sigma^k$ ,

$$\tag{4.8}$$

with R denoting the scalar curvature of Y. (The top bullet in (4.8) can also be written as $\langle B_A \otimes \sigma \rangle = \langle \mathcal{B}, \sigma \rangle - \tfrac{1}{2} \mathrm{Ric} + \tfrac{1}{4} R\mathfrak{g}$.)

**c) Equations for $\langle \mathfrak{c} \otimes \sigma \rangle$**

Equations for the t-derivatives of $c$, $\mathfrak{c}^+$ and $\mathfrak{c}^-$ are implied by the middle equation in (2.7) and (4.5). For future reference, these are:

- $\tfrac{\partial}{\partial t} c + \tfrac{2}{t} c = \tfrac{1}{\sqrt{3}} \langle \sigma, B_A \rangle + \tfrac{2}{\sqrt{3}} c^2 + \tfrac{1}{\sqrt{3}} |\mathfrak{c}^-|^2 - \tfrac{1}{\sqrt{3}} |\mathfrak{c}^+|^2$ .
- $\tfrac{\partial}{\partial t} \mathfrak{c}^+ - \tfrac{1}{t} \mathfrak{c}^+ = \langle \sigma \otimes B_A \rangle^+ - \langle \sigma \otimes *(\mathfrak{c} \wedge \mathfrak{c}) \rangle^+$ .
- $\tfrac{\partial}{\partial t} \mathfrak{c}^- + \tfrac{1}{t} \mathfrak{c}^- = \langle \sigma \otimes B_A \rangle^- - \langle \sigma \otimes *(\mathfrak{c} \wedge \mathfrak{c}) \rangle^-$.

$$\tag{4.9}$$

In the second and third bullets (and in what follows), the notation has $\langle\ \rangle^\pm$ denoting the symmetric traceless (with +) and anti-symmetric parts of the $\otimes^2 T^*Y$ valued 1-form.

Some parenthetical remarks about the second and third bullet equations: The equation in the second bullet in (4.9) can be written as

$$\tfrac{\partial}{\partial t} \mathfrak{c}^+ - \tfrac{1}{t} \mathfrak{c}^+ + \tfrac{2}{\sqrt{3}} c\, \mathfrak{c}^+ = \langle \sigma \otimes B_A \rangle^+ - \langle \sigma \otimes *(\tilde{\mathfrak{c}} \wedge \tilde{\mathfrak{c}}) \rangle^+ \tag{4.10}$$

where $\tilde{\mathfrak{c}} = \mathfrak{c} - \tfrac{1}{\sqrt{3}} c \sigma$. (Thus, $\langle \tilde{\mathfrak{c}} \otimes \sigma \rangle = \mathfrak{c}^+ + \mathfrak{c}^-$ which is the traceless part of $\langle \mathfrak{c} \otimes \sigma \rangle$.) Using the second bullet of (4.4), the identity in (4.9) can be also be written as

$$\tfrac{\partial}{\partial t} \mathfrak{c}^+ + \mathfrak{t} = \langle \sigma \otimes B_A \rangle^+ + \langle \tilde{\mathfrak{c}} \otimes \tilde{\mathfrak{c}} \rangle^+ + \langle \sigma \otimes *(\tilde{\mathfrak{c}} \wedge \tilde{\mathfrak{c}}) \rangle^+ . \tag{4.11}$$

Meanwhile, the equation in the third bullet of (4.9) can be written as

$$\frac{\partial}{\partial t}\mathfrak{c}^- + \frac{1}{t}\mathfrak{c}^- - \frac{2}{\sqrt{3}}\, c\,\mathfrak{c}^- = \langle \sigma\otimes B_A\rangle^- - \mathcal{T}'(\mathfrak{c}^-, \mathfrak{c}^+) \tag{4.12}$$

where $\mathcal{T}'$ is a certain canonical homomorphism from $(\otimes^2 T^*Y) \otimes (\otimes^2 T^*Y)$ to $\otimes^2 T^*Y$. The particular form of $\mathcal{T}'$ is not relevant except for the fact that its norm is bounded by $c_0$.

**d) The integral of $\langle \mathfrak{a}, \mathcal{B}\rangle$**

The identity in (2.19) is less useful than a variant for the t-derivative of the function of t that is given by the $\{t\} \times Y$ integral of $\langle \mathfrak{a}, \mathcal{B}\rangle$. (See also [MW], [LT].) For the present purposes, the important point is that the t-derivative of the integral of $\langle \mathfrak{a}, \mathcal{B}\rangle$ on $\{t\} \times Y$ integral can be written schematically as

$$\frac{d}{dt}\int_{\{t\}\times Y} \langle \mathfrak{a}, \mathcal{B}\rangle = \int_{\{t\}\times Y} |B_A|^2 + \frac{1}{2}\int_{\{t\}\times Y} |E_A|^2 + \frac{1}{2}\int_{\{t\}\times Y} |\nabla_A^{\perp}\mathfrak{a}|^2 + \mathcal{E}\ , \tag{4.13}$$

where the term denoted by $\mathcal{E}$ is

$$\mathcal{E} = -\frac{1}{12}\int_{\{t\}\times Y} R\langle \sigma, B_A\rangle + \frac{1}{2}\int_{\{t\}\times Y} \langle \mathrm{Ric}, \langle \sigma\otimes B_A\rangle^+\rangle + \frac{1}{12}\int_{\{t\}\times Y} R(3|\mathfrak{c}^+|^2 + |\mathfrak{c}^-|^2)$$
$$+ \frac{1}{2}\int_{\{t\}\times Y} \langle \mathrm{Ric}, \langle \tilde{\mathfrak{c}}\otimes\tilde{\mathfrak{c}}\rangle^+ + \langle \sigma\otimes *(\tilde{\mathfrak{c}}\wedge\tilde{\mathfrak{c}})\rangle\rangle\ . \tag{4.14}$$

The important point is that $\mathcal{E}$ has no terms with powers of $\frac{1}{t}$ and that it obeys

$$|\mathcal{E}| \le c_0\Big(\int_{\{t\}\times Y} |B_A|^2\Big)^{1/2} + c_0\int_{\{t\}\times Y} (|\mathfrak{c}^+|^2 + |\mathfrak{c}^-|^2)\ . \tag{4.15}$$

(Remember that $\tilde{\mathfrak{c}}$ is identified by $\sigma$ with $\mathfrak{c}^+ + \mathfrak{c}^-$.)

As for the derivation: The identity in (4.13) is derived by writing $B_A$ as $\mathcal{B} + B_\Gamma$ using the top bullet of (4.8). Then subtract the $\{t\} \times Y$ integral of $\langle \frac{\partial}{\partial t}\mathfrak{a}, B_\Gamma\rangle$ from the right hand side of (2.19). But, use the formula in the second bullet of (2.7) to write this subtracted term as the $\{t\} \times Y$ integral of $\langle B_A, B_\Gamma\rangle - \langle \mathfrak{a}\wedge\mathfrak{a}\wedge B_\Gamma\rangle$. The integral of the $\langle B_A, B_\Gamma\rangle$ term is accounted for by the two left most two integrals in (4.14). As for the rest, take note of all the cancellations that occur between the $-\langle \mathfrak{a}\wedge\mathfrak{a}\wedge B_\Gamma\rangle$ integral and the $\langle \mathrm{Ric}, \langle \mathfrak{a}\otimes\mathfrak{a}\rangle\rangle$ integral when $\mathfrak{a}$ is written using (4.2) and (4.3).

**e) Integrals of $|\mathfrak{b}|^2$ and $|\nabla_A{}^{\perp}\mathfrak{c}|$**

The lemma that follows momentarily asserts a useful inequality concerning integrals of $|E_A|^2$ and $|\nabla_A{}^{\perp}\mathfrak{a}|^2$ on domains $[t^-, t] \times Y$. This lemma introduces the following notation: A set $I \subset (0, 1]$ is said to have asymptotic full measure if the fraction of $(0, t]$ containing I limits to 1 as $t \to 0$. (Formally: Given $\varepsilon > 0$, there exists $t_\varepsilon > 0$ such the measure of $I \cap (0, t_\varepsilon]$ is greater $(1 - \varepsilon) t_\varepsilon$.)

**Lemma 4.1**: *There exist* $\kappa > 0$ *with the following significance: Let* $(A, \mathfrak{a})$ *denote a Nahm pole solution. If* $t^-$ *is sufficiently small and chosen from a certain asymptotically full measure set (which is determined by* $(A, \mathfrak{a})$*), then the corresponding versions of* $\mathfrak{b}$ *and* $\mathfrak{c}$ *on* $[t^-, t]$ *for* $t > t^-$ *obey*

$$\frac{1}{t}\int_{\{t\}\times Y} |\mathfrak{b}|^2 + \int_{[t^-,t]\times Y} \left(\frac{1}{s^2}|\mathfrak{b}|^2 + |\nabla_A^{\perp}\mathfrak{c}|^2\right) \leq \kappa \int_{[0,t]\times Y} \left(|E_A|^2 + |\nabla_A^{\perp}\mathfrak{a}|^2\right)$$

To be sure: It is not a mistake that the integration domain on the right is $[0, t] \times Y$ whereas that for the right most integral on the left is $[t^-, t] \times Y$.

***Proof of Lemma 4.1***: Given $\varepsilon > 0$ but less than $10^{-6}$, fix $t_\varepsilon > 0$ so that (3.1) holds on the interval $(0, t_\varepsilon]$. Then, by virtue of (3.1) and the second bullet in (3.13):

$$\int_{(0,t_\varepsilon]\times Y} \frac{1}{s^2}|\mathfrak{b}^{\perp}|^2 \leq 16 \int_{[0,t_\varepsilon]\times Y} |\nabla_A^{\perp}\mathfrak{a}|^2 \ . \tag{4.16}$$

Now, given $n \in \{1, 2, \ldots\}$, let $J_n \subset [2^{-n}t_\varepsilon, 2^{-n+1}t_\varepsilon]$ denote the set of times t where

$$\frac{1}{t}\int_{\{t\}\times Y} |\mathfrak{b}^{\perp}|^2 \geq 128 \int_{[2^{-n}t_\varepsilon, 2^{-n+1}t_\varepsilon]\times Y} |\nabla_A^{\perp}\mathfrak{a}|^2 \ . \tag{4.17}$$

For each n, let $|J_n|$ denote the measure of $J_n$ and let $\mu_n = 2^n t_\varepsilon^{-1} |J_n|$ which is the fraction of the interval between $2^{-n}t_\varepsilon$ and $2^{-n+1}t_\varepsilon$ that is accounted for by the points of $J_n$. Then,

$$\sum\nolimits_{n=2,3,\ldots} \mu_n \leq \frac{1}{4} \ , \tag{4.18}$$

so as not to run afoul of (4.16). Therefore, the complement (call it I) in $(0, t_\varepsilon]$ of the union of the $J_n$'s has asymptotically full measure.

With the preceding understood, agree to choose $t^-$ from the complement of the union of the $J_n$'s so that the $t = t^-$ version of (4.17) holds with the inequality going the other way. Meanwhile for $t > t^-$, one has

$$\frac{1}{t} \int_{\{t\}\times Y} |\mathfrak{b}^\perp|^2 + \frac{1}{2} \int_{[t^-,t]\times Y} \frac{1}{s^2} |\mathfrak{b}|^2 \leq \frac{1}{t^-} \int_{\{t^-\}\times Y} |\mathfrak{b}^\perp|^2 + 2 \int_{[t^-,t]\times Y} |\tfrac{\partial}{\partial t} \mathfrak{b}|^2 \ , \tag{4.19}$$

which is called Hardy's inequality. But it is probably Hardy's least significant achievement given that its proof amounts to writing the integral over $[t^-, t]$ of a 1-form $\frac{1}{s^2} f(s)\,ds$ as that of $-f(s)d(\frac{1}{s})$ and then integrating by parts and then using the algebraic inequality $2|ab| \leq \frac{1}{2} a^2 + 2|b|^2$. (The function $f$ in this case is the integral of $|\mathfrak{b}|^2$ over the slice $\{s\} \times Y$; and in this case, $a = \frac{1}{s}|\mathfrak{b}|$ and $b = |\frac{\partial}{\partial t} \mathfrak{b}|$.) By virtue of the choice for $t^-$, and because of the second bullet in (4.8), the inequality in (4.19) implies in turn:

$$\frac{1}{t} \int_{\{t\}\times Y} |\mathfrak{b}^\perp|^2 + \frac{1}{2} \int_{[t^-,t]\times Y} \frac{1}{s^2} |\mathfrak{b}|^2 \leq c_0 \Big( \int_{[0,t_\varepsilon]\times Y} |\nabla_A^\perp \mathfrak{a}|^2 + \int_{[t^-,t]\times Y} |E_A|^2 \Big) . \tag{4.20}$$

Meanwhile, supposing that $s \geq t^-$, one has, for any $\delta \in (0, \frac{1}{10}\,]$, the inequality

$$|\nabla_A{}^\perp \mathfrak{a}|^2 \geq \delta |\nabla_A{}^\perp \mathfrak{c}|^2 - c_0 \delta \frac{1}{s^2} |\mathfrak{b}|^2 . \tag{4.21}$$

This last inequality (for $\delta = c_0{}^{-2}$) with (4.20) imply what is asserted by Lemma 4.1.

## 5. The behavior of K, N and the integrals of $\langle \mathfrak{a}, \mathcal{B} \rangle$ and $|\mathfrak{t}|$ on $\{t\} \times Y$

The first subsection in this section considers the behavior of the functions K and N at times t where t can be as large as $\mathcal{O}(1)$. The other subsections derive relations between K and the integrals of $\langle \mathfrak{a}, \mathcal{B} \rangle$ and $|\mathfrak{t}|$ on the slices $\{t\} \times Y$ for different values of $t \geq t^-$. With regards to this time $t^-$; assume henceforth that it is chosen so that the conclusions of Lemma 4.1 hold for all $t < 10^8 t^-$.

### a) The behavior of the functions K and N

The function K is defined by (2.8) and the function N is defined implicitly by writing the derivative of K as in (2.11).

**Lemma 5.1**: *There exists* $\kappa > 4$ *with the following significance: Let* $(A, \mathfrak{a})$ *denote a given Nahm pole solution. The associated versions of* K *and* N *have the properties listed below.*

- *If* t *is sufficiently small, then* $|N(t) - 1| \le \kappa t^2$ .
- *There exists* $t_\Delta \in (0, \frac{1}{\kappa^4}] \cup \{1\}$ *such that*
  a) *If* $t_\Delta < \frac{1}{\kappa^4}$ *then* $N(t) \le -\sqrt{t}$ *where* $t \in [t_\Delta, \frac{1}{\kappa^4}]$, *and* $K(s) \ge K(t)$ *if both* $t \in [t_\Delta, \frac{1}{\kappa^4})$ *and* $s \in [t, \frac{1}{\kappa^4}]$ ;
  b) $N(t) > -\sqrt{t}$ *where* $t \le t_\Delta$; *and* $K(t) \ge (1 - \frac{1}{\kappa}\sqrt{s})K(s)$ *when* $t \in (0, t_\Delta)$ *and* $s \in (t, t_\Delta)$.

Note that the first bullet implies that $\lim_{t\to 0} N(t) = 1$. This is why $t_\Delta$ is positive.

Lemma 5.1 implies this: If $t < \frac{1}{\kappa^4}$ and $s \in (t, \frac{1}{\kappa^4}]$, then K on $[t, s]$ is at most $c_0$ times the maximum of K(t) and K(s). This is used in what follows (often implicitly).

The proof of Lemma 5.1 uses a formula for the derivative of the function N:

$$\frac{d}{dt} N = \frac{N+2N^2}{t} - \frac{t}{K^2} \int_{\{t\}\times Y} (|\nabla_{A_t}\mathfrak{a}|^2 + 2|\mathfrak{a}\wedge\mathfrak{a}|^2 + |\nabla_A^\perp \mathfrak{a}|^2 + \langle Ric, \langle \mathfrak{a} \otimes \mathfrak{a}\rangle\rangle) \ . \tag{5.1}$$

This identity is obtained by directly differentiating the expression in (2.11) and invoking (2.10). This formula is also invoked in subsequent proofs.

***Proof of Lemma 5.1***: The first bullet of the lemma follows from Lemmas 3.2-3.4 and (3.29). In both a) and b) of the second bullet, the assertions about K follow using (2.11) from the assertions about N. To prove the assertions about N, fix $k \ge 1$ and let $t \in (0, 1]$ denote a time where $N(t) = -\sqrt{t}$ assuming such time exists. If not, set $t_\Delta$ equal to 1. By virtue of (5.1),

$$\frac{d}{dt}(N + \sqrt{t}) \le -\frac{1}{2\sqrt{t}} + 2 + c_0 t \ . \tag{5.2}$$

Hence, if it is the case that $\sqrt{t} < \frac{1}{8} c_0^{-1}$, then the function $N + \sqrt{t}$ is decreasing at t. Therefore, the function N never becomes greater than $-\sqrt{t}$ once it equals $-\sqrt{t}$ (except perhaps when $t \ge \frac{1}{8} c_0^{-1}$).

The next lemma relates N to integrals of $|\nabla_A \mathfrak{a}|^2$ and $|\mathfrak{a} \wedge \mathfrak{a}|^2$. This lemma uses $\kappa_*$ to denote Lemma 5.1's version of $\kappa$.

**Lemma 5.2**: *There exists* $\kappa > 1$ *with the following significance*: *Let* $(A, \mathfrak{a})$ *denote a given Nahm pole solution. Use this Nahm pole solution to define* $t_\Delta$ *as in Lemma 5.1.*

- *If* $t \in (0, \frac{1}{4\kappa_*^4}]$ *and* $t \le t_\Delta$, *then* $|N(t)|$ *differs by at most* $\kappa t t_\Delta$ *from*

$$\frac{t}{\kappa^2(t)} \int_{[t,t_\Delta]\times Y} (|\nabla_A \mathfrak{a}|^2 + 2|\mathfrak{a}\wedge\mathfrak{a}|^2)$$

- *If* $t \in (0, \frac{1}{4\kappa_*^4}]$ *and* $t \geq t_\Delta$, *then* $N(t)$ *differs by at most* $\kappa t^2$ *from*

$$\frac{t}{\kappa^2(t)} \int_{[t_0,t]\times Y} (|\nabla_A \mathfrak{a}|^2 + 2|\mathfrak{a}\wedge\mathfrak{a}|^2)$$

*where* $t_0 \in (0, t_\Delta)$ *is the largest time where* $N$ *is zero.*

***Proof of Lemma 5.2***: The bullets are proved in reverse order. To prove the second bullet, use (2.12) with t in (2.12) replaced by $t_\Delta$ and with s in (2.12) replaced by Lemma 5.2's version of t. Lemma 5.1 is used to bound the norm of the integral on the right hand side with the Ricci curvature by $c_0 t \kappa^2(t)$.

To prove the second bullet: If $t_\Delta \leq \frac{1}{\kappa_*^4}$, then the top bullet follows from the version of (2.12) with s taken to be $t_\Delta$. Lemma 5.1 is invoked in this case to bound the norm of the integral of the Ricci curvature term that appears on the right hand side of (2.12) by $c_0 t_0 \kappa^2(t)$.

Suppose now that $N > 0$ on the whole of $(0, \frac{1}{\kappa_*^4}]$ (so that $t_\Delta = 1$). Consider now (2.12) with s from $[\frac{1}{16\kappa_*^4}, \frac{1}{4\kappa_*^4}]$. The Ricci curvature integral on the right hand side is treated just as before. But, there is now a negative term on the left hand side of (2.12), which is the $-N(s)\frac{\kappa(s)^2}{s}$ term. To bound the norm of this term: Use $\chi$ from Section 1f to construct a smooth, non-negative function with compact support on $[\frac{1}{16\kappa_*^4}, \frac{1}{4\kappa_*^4}]$ whose derivative has norm bounded by $c_0^{-1}$. (Keep in mind that $t_* \geq c_0^{-1}$.) Denote this function by $\chi_*$ and let $\mu_*$ denote the integral of $\chi_*$ on $[\frac{1}{16\kappa_*^4}, \frac{1}{4\kappa_*^4}]$. Now take (2.12) with s in the support of $\chi_*$ and multiply both sides of (2.12) for this choice of s by $\mu_*^{-1}\chi_*(s)$. Then, integrate the result with respect to the variable s. Write the integral of $-\mu_*^{-1}\chi_*(s)N(s)\frac{\kappa(s)^2}{s}$ as the integral of $\mu_*\chi_*(s)\frac{1}{2}(\frac{d}{ds}\kappa^2)$ and integrate by parts. The resulting of doing that is an integral that is bounded by $c_0\kappa^2(t)$ (by virtue of Lemma 5.1).

The last lemma in this subsection gives an a priori lower bound for the time it takes $N$ to increase from a given value to twice that value.

**Lemma 5.3**: *There exists* $\kappa > 8$ *with the following significance: Let* $(A, \mathfrak{a})$ *denote a Nahm pole solution; and let* $N$ *denote the associated function from (2.11). Fix* $t > 0$ *where* $N(t) \geq \kappa t^2$. *If* s *is between* t *and* $(\frac{2+2N(t)}{1+2N(t)})^{1/2}t$, *then* $N(s) \leq 2N(t)$.

Here is an example that is used later: If $N(t) \le 2$, then $N(s) \le 4$ if $s \in [t, \frac{101}{100} t]$. Another example: If $N(t) \le 1$, then $N(s) \le 2$ if $s \in [t, \frac{11}{10} t]$.

***Proof of Lemma 5.3***: The identity in (5.1) leads to the differential inequality

$$\frac{d}{dt} N \le \frac{N + 2N^2}{t} + c_0 t \,. \tag{5.7}$$

And, if $N(t) \ge c_0 t^2$, then this in turn implies that

$$\frac{d}{dt} N \le 2 \frac{N + N^2}{t} \,. \tag{5.8}$$

The inequality in (5.8) can be directly integrated to see that $N(s) < 2N(t)$ when s is between t and $(\frac{2 + 2N(t)}{1 + 2N(t)})^{1/2} t$.

A second formula for the derivative of $N$ is obtained from (5.1) by writing $\nabla_{At}\mathfrak{a}$ in the integral on the right hand side as $(\nabla_{At}\mathfrak{a} + \frac{N}{t}\mathfrak{a}) - \frac{N}{t}\mathfrak{a}$, and also writing $\mathfrak{a} \wedge \mathfrak{a}$ in this integral as $(\mathfrak{a} \wedge \mathfrak{a} - * \frac{N}{t}\mathfrak{a}) + * \frac{N}{t}\mathfrak{a}$. The result of doing these substitutions (and appealing to (2.9)) is the promised second formula:

$$\frac{d}{dt} N = \frac{N(1 - N)}{t} - \frac{t}{K^2} \int_{\{t\}\times Y} (|\nabla_t \mathfrak{a} + \tfrac{N}{t}\mathfrak{a}|^2 + 2|\mathfrak{a}\wedge\mathfrak{a} - \tfrac{N}{t} * \mathfrak{a}|^2 + |\nabla_A^{\perp}\mathfrak{a}|^2 + \langle \mathrm{Ric}, \langle \mathfrak{a} \otimes \mathfrak{a}\rangle\rangle) - \frac{4N}{K^2} \int_{\{t\}\times Y} \langle \mathfrak{a}, B_A \rangle \,. \tag{5.9}$$

This second formula is more or less useful than that in (5.1) depending on what can be said about the integral on $\{t\} \times Y$ of $\langle \mathfrak{a}, B_A \rangle$.

**b) Initial bounds for the integral of $\langle \mathfrak{a}, \mathcal{B} \rangle$ on $\{t\} \times Y$**

The upcoming lemmas refer to an $\mathcal{O}(1)$ time $t_*$ that is defined as follows: Let $\kappa_*$ denote the version of the number $\kappa$ that appears in Lemma 5.1 and set $t_* = \frac{1}{8\kappa_*^4}$. By definition, the conclusions of Lemmas 5.1 and 5.2 hold on $(0, t_*]$.

The first lemma also refer to $\mathfrak{cs}_\infty$ which is the $t \to \infty$ limit of (1.3)'s function $\mathfrak{cs}$.

**Lemma 5.4**: *Let* $(A, \mathfrak{a})$ *denote Nahm pole solution. If* $t \in [t^-, t_*]$ *and* $N(t) \le 0$, *then*

- $\int_{\{t\}\times Y} \langle \mathfrak{a}, B_A \rangle \le \frac{3}{2} \mathfrak{cs}_\infty$.
- $\int_{\{t\}\times Y} \langle \mathfrak{a}, \mathcal{B} \rangle \le \kappa K(t) + \frac{3}{2} \mathfrak{cs}_\infty$ *with* $\kappa$ *being independent of* $(A, \mathfrak{a})$.

***Proof of Lemma 5.4***: The second bullet follows from the top because $|B_A - \mathcal{B}| \le c_0$. To prove the top bullet: If $\mathrm{N}(t) \le 0$, then the second bullet of (2.9) and (2.11) are compatible only in the event that the $\{t\} \times Y$ integral of $\langle \mathfrak{a}, B_A\rangle$ is greater than that of $\langle \mathfrak{a} \wedge \mathfrak{a} \wedge \mathfrak{a}\rangle$. This with (1.3) implies that the $\{t\} \times Y$ integral of $\langle \mathfrak{a}, B_A\rangle$ is less than $\frac{3}{2}\,\mathfrak{cs}(t)$. Meanwhile, $\frac{3}{2}\,\mathfrak{cs}(t)$ is less than $\frac{3}{2}\,\mathfrak{cs}_\infty$ because the function $\mathfrak{cs}$ is increasing.

The next lemma concerns the $\{t\} \times Y$ integrals of $\langle \mathfrak{a}, \mathcal{B}\rangle$ with no a priori assumption about the sign of $\mathrm{N}$.

**Lemma 5.5**: *There exists* $\kappa > 8$ *which is independent of* $(A, \mathfrak{a})$ *and such that if* $t \le s$ *and both are in* $[t^-, t_*]$*, then*

- $\displaystyle\int_{\{t\}\times Y} \langle \mathfrak{a}, \mathcal{B}\rangle + \frac{1}{4}\int_{[t,s]\times Y} (|B_A|^2 + |E_A|^2 + |\nabla_A^{\perp}\mathfrak{a}|^2) \le \int_{\{s\}\times Y} \langle \mathfrak{a}, \mathcal{B}\rangle + \kappa\Big((s-t) + \int_{[t,s]\times Y} (|\mathfrak{c}^+|^2 + |\mathfrak{c}^-|^2)\Big).$
- $\displaystyle\int_{\{t\}\times Y} \langle \mathfrak{a}, \mathcal{B}\rangle + \frac{1}{4}\int_{[t,s]\times Y} (|B_A|^2 + |E_A|^2 + |\nabla_A^{\perp}\mathfrak{a}|^2) \le \int_{\{s\}\times Y} \langle \mathfrak{a}, \mathcal{B}\rangle + \kappa\Big(s + \int_{[t,s]} \mathrm{K}^2\Big).$
- $\displaystyle\int_{\{t\}\times Y} \langle \mathfrak{a}, \mathcal{B}\rangle + \frac{1}{4}\int_{[t,s]\times Y} (|B_A|^2 + |E_A|^2 + |\nabla_A^{\perp}\mathfrak{a}|^2) \le \int_{\{s\}\times Y} \langle \mathfrak{a}, \mathcal{B}\rangle + \kappa s(1 + \mathrm{K}^2(t) + \mathrm{K}^2(s)).$

***Proof of Lemma 5.5***: To prove the top bullet of the lemma, fix $t < s$ with both in $[t^-, t_*]$. Integrate (4.13) between t and s and use the bound in (4.15). The second bullet follows from the top bullet because the integral of $|\mathfrak{c}^+|^2 + |\mathfrak{c}^-|^2$ that appears on the far right in the top bullet is not greater than that of $|\mathfrak{a}|^2$ on the same domain, which is the integral of $\mathrm{K}^2$ from t to s. The third bullet follows from the second bullet because the integral of $\mathrm{K}^2$ on $[t, s]$ is no greater than $c_0\, s\,(\mathrm{K}^2(t) + \mathrm{K}^2(s))$ due to Lemma 5.1.

The final lemma in this subsection is for the most part a corrollary of Lemma 5.5. This lemma uses $\mathrm{K}_*$ to denote the number $\frac{1}{2}(\mathrm{K}(\frac{1}{8}t_*) + \mathrm{K}(t_*))$.

**Lemma 5.6**: *There exists* $\kappa > 8$ *with the following significance: Suppose that* $(A, \mathfrak{a})$ *is a Nahm pole solution. If* $t \le s$ *and both are in* $[t^-, \frac{1}{8}t_*]$ *then*

- $\displaystyle\int_{\{t\}\times Y} \langle \mathfrak{a}, \mathcal{B}\rangle \le \kappa(1 + s\mathrm{K}^2(t) + \mathrm{K}^2(s) + \mathrm{K}_*^{\,2}).$
- $\displaystyle\int_{\{t\}\times Y} \langle \mathfrak{a}, \mathcal{B}\rangle \le \kappa(1 + s\,\mathrm{K}^2(t) + s\,\mathrm{K}^2(s)) + \mathrm{K}(s) + \mathfrak{cs}_\infty)$ *if* $\mathrm{N}(s) \le 0$.

*In general, if* $\int_{\{t\}\times Y} \langle \mathfrak{a}, \mathcal{B}\rangle \geq 0$, *then*

- $\int_{[t,s]\times Y} (|B_A|^2 + |E_A|^2 + |\nabla_A^{\perp}\mathfrak{a}|^2) \leq \kappa(1 + s\kappa^2(t) + \kappa^2(s) + \kappa_*^2)$ *when* $\int_{\{s\}\times Y} \langle \mathfrak{a}, \mathcal{B}\rangle \geq 0$ .
- $\int_{[t,s]\times Y} (|B_A|^2 + |E_A|^2 + |\nabla_A^{\perp}\mathfrak{a}|^2) \leq \kappa s(1 + \kappa^2(t) + \kappa^2(s))$ *when* $\int_{\{s\}\times Y} \langle \mathfrak{a}, \mathcal{B}\rangle \leq 0$.

***Proof of Lemma 5.6***: The bullets are proved in reverse order. The fourth bullet of this lemma follows directly from the third bullet of Lemma 5.5. To prove the third bullet: Invoke the top bullet of this lemma (supposing that it is true) with t replaced by the third bullet's version of s and with s replaced by $\frac{1}{8}t_*$. This bounds the $\{s\}\times Y$ integral of $\langle \mathfrak{a}, \mathcal{B}\rangle$ by $c_0(1 + \kappa^2(s) + \kappa_*^2)$. Now use this bound for the $\{s\}\times Y$ integral of $\langle \mathfrak{a}, \mathcal{B}\rangle$ in the third bullet of Lemma 5.5. To prove the second bullet: Use Lemma 5.4's second bullet to bound the $\{s\}\times Y$ integral of $\langle \mathfrak{a}, \mathcal{B}\rangle$ in that appears on the right hand side of the third bullet of Lemma 5.5. The proof of the top bullet of Lemma 5.6 has two parts.

*Part 1*: This first part gives a bound for the $\{t\}\times Y$ integral of $\langle \mathfrak{a}, \mathcal{B}\rangle$ for certain values of t from the interval $[\frac{1}{4}t_*, \frac{1}{2}t_*]$. To start, use the function $\chi$ from Section 1f to construct a function that is equal to 1 where $t \in [\frac{1}{4}t_*, \frac{1}{2}t_*]$, and equal to zero where $t \geq t_*$ and where $t \leq \frac{1}{8}t_*$. This function can and should be constructed so that the norm of its first and second derivatives are bounded respectively by $c_0$. (Remember that $t_* = \frac{1}{8\kappa_*^4}$ .) Denote this function by $\chi_*$ (it is not the same as the previous incarnation of $\chi_*$). Now take the inner product of both sides of (2.4) with $\chi_*\mathfrak{a}$ and integrate the resulting identity on $[\frac{1}{8}t_*, t_*]\times Y$. Integration by parts and an appeal to Lemma 5.1 leads to the bound

$$\int_{[\frac{1}{4}t_*,\frac{1}{2}t_*]\times Y} (|\nabla_A\mathfrak{a}|^2 + |[\mathfrak{a}_*, \mathfrak{a}]|^2) \leq c_0\kappa_*^2. \tag{5.10}$$

What with the second bullet in (2.7), this in turn implies that

$$\int_{[\frac{1}{4}t_*,\frac{1}{2}t_*]\times Y} |B_A|^2 \leq c_0\kappa_*^2 , \tag{5.11}$$

and therefore (by invoking Lemma 5.1), that

$$\int_{[\frac{1}{4}t_*,\frac{1}{2}t_*]\times Y} |\langle \mathfrak{a}, \mathcal{B}\rangle| \leq c_0\kappa_*^2 . \tag{5.12}$$

This implies that $\int_{\{t\}\times Y} \langle \mathfrak{a}, \mathcal{B}\rangle \le c_0 K_*^{\,2}$ for at least half of the times t in $[\frac{1}{4} t_*, \frac{1}{2} t_*]$.

*Part 2*: Fix $t \in [t^-, \frac{1}{4} t_*]$ and suppose that

$$\int_{\{t\}\times Y} \langle \mathfrak{a}, \mathcal{B}\rangle \ge 0 \tag{5.13}$$

because there is nothing to prove for the top bullet of Lemma 5.6 otherwise. Using the third bullet of Lemma 5.5 with s chosen appropriately in $[\frac{1}{4} t_*, \frac{1}{2} t_*]$ gives

$$\int_{\{t\}\times Y} \langle \mathfrak{a}, \mathcal{B}\rangle \le c_0 K_*^{\,2} + c_0 (1 + K^2(t)) \,. \tag{5.14}$$

Change notation now and use s instead of t in (5.14) to conclude from (5.14) that if $s \le \frac{1}{8} t_*$, then the integral of $\langle \mathfrak{a}, \mathcal{B}\rangle$ on $\{s\} \times Y$ is no greater than $c_0 (1 + K^2(s) + K_*^{\,2})$. Use the latter bound in the third bullet of Lemma 5.5 to obtain the first bullet of Lemma 5.6.

**c) Initial bounds for the integral of $|\mathfrak{t}|$**

The next lemma concerns the integral of $|\mathfrak{t}|$ over $\{t\}\times Y$. (Remember that the tensor $\mathfrak{t}$ is $\langle \mathfrak{a} \otimes \mathfrak{a}\rangle - \frac{1}{3} |\mathfrak{a}|^2 \mathfrak{g}$.)

**Lemma 5.7**: *There exists* $\kappa > 8$ *which is independent of* $(A, \mathfrak{a})$ *and has the following significance: If* t and s are in $(t^-, \frac{1}{8} t_*]$ *with* $s > t$, *then*

- $\left| \int_{\{t\}\times Y} |\mathfrak{t}| - \int_{\{s\}\times Y} |\mathfrak{t}| \right| \le \kappa \left( \int_{[t,s]} K^2 \right)^{1/2} \left( \int_{\{s\}\times Y} \langle \mathfrak{a}, \mathcal{B}\rangle - \int_{\{t\}\times Y} \langle \mathfrak{a}, \mathcal{B}\rangle + \kappa \left( (s-t) + \int_{[t,s]\times Y} (|\mathfrak{c}^+|^2 + |\mathfrak{c}^-|^2) \right) \right)^{1/2}$ .
- $\left| \int_{\{t\}\times Y} |\mathfrak{t}| - \int_{\{s\}\times Y} |\mathfrak{t}| \right| \le \kappa \left( \int_{[t,s]} K^2 \right)^{1/2} \left( \int_{\{s\}\times Y} \langle \mathfrak{a}, \mathcal{B}\rangle - \int_{\{t\}\times Y} \langle \mathfrak{a}, \mathcal{B}\rangle + \kappa \left( s + \int_{[t,s]} K^2 \right) \right)^{1/2}$ .

*In particular, if the* $\{t\}\times Y$ *integral of* $\langle \mathfrak{a}, \mathcal{B}\rangle$ *is non-negative, then*

- $\int_{\{t\}\times Y} |\mathfrak{t}| \le \kappa (1 + K^2(s) + s K^2(t) + K_*^{\,2})$

*In particular, if* $\alpha \ge 0$ *and the integral of* $\langle \mathfrak{a}, \mathcal{B}\rangle$ *on* $\{t\}\times Y$ *is greater than* $-\alpha K^2(t)$, *then*

- $\int_{\{t\}\times Y} |\mathfrak{t}| \le \kappa (1 + K^2(s) + \sqrt{s}\, \alpha^{1/2} K^2(t) + K_*^{\,2})$

***Proof of Lemma 5.7***: Use the identity in (2.18) to see that

$$\Big| \int_{\{t\}\times Y} |\mathfrak{t}| - \int_{\{s\}\times Y} |\mathfrak{t}| \ \Big| \le c_0 \int_{[t,s]\times Y} |\mathfrak{a}||B_A| \ . \tag{5.15}$$

The right most integral on the right hand side of (5.15) is at most the product

$$\Big( \int_{[t,s]\times Y} |\mathfrak{a}|^2 \Big)^{1/2} \Big( \int_{[t,s]\times Y} |B_A|^2 \Big)^{1/2} . \tag{5.16}$$

The first and second bullets of Lemma 5.7 follows directly from this and the respective first and second bullets of Lemma 5.5. The third and fourth bullets of the lemma follows directly from the second bullet using Lemma 5.1 and the first bullet of Lemma 5.5. With regards to these last two bullets: The $\{s\} \times Y$ integral of $|\mathfrak{t}|$ is not larger than $c_0 K^2(s)$). Also, Lemma 5.1 is used to bound the integral of $K^2$ that appears on the right hand side of the third bullet of this lemma by $c_0\, s\, (K^2(t) + K^2(s))$. Meanwhile, the top bullet of Lemma 5.6 with s replaced by t and with t replaced by $t_*$ is used to bound the $\{s\}\times Y$ integral of $\langle \mathfrak{a}, \mathcal{B}\rangle$ that appears on the right hand side of the third bullet of this lemma.

## d) Initial bounds for $K$

The lemma that follows supplies an a priori initial bound for the function $K$. The lemma uses $R_*$ to denote the supremum over Y of the norm of the metric's Riemann curvature tensor.

**Lemma 5.8**: *There exists* $\kappa > 16(1 + R_*)$ *which is independent of* $(A, \mathfrak{a})$ *and has the following significance: Fix* $t \in (t^-, \frac{1}{\kappa} t_*]$ *and suppose that*

a) $K(t) \ge \kappa(1 + K(\frac{1}{\kappa} t_*) + K_*)$

b) $N(t) \le 2$ .

c) $\int_{\{t\}\times Y} \langle \mathfrak{a}, \mathcal{B}\rangle \ \ge -K^2(t)$.

*Then* $K(t) \le \kappa \frac{1}{t}$

This lemma is proved momentarily. Some remarks follow directly.

Some remarks on Requirement a): Supposing for the moment that $\kappa$ is a given number greater than 16, then the requirement $K(t) \ge \kappa(1 + K(\frac{1}{\kappa} t_*) + K_*)$ can not hold if $t > t_\Delta$ with $t_\Delta$ defined by Lemma 5.1. This is because $K$ is increasing between $t_\Delta$ and $t_*$. (But note that since $K$ is increasing between $t_\Delta$ and $t_*$, if $t > t_\Delta$ (supposing that $t_\Delta < t_*$), then $K(t) < K(t_*)$; and this is less than $\frac{1}{t}$ when $t < \frac{1}{K(t_*)}$ .)

An observation concerning Requirement a) of the lemma and subsequent requirements of the form $\text{K}(\cdot) \geq \kappa(1+\text{K}(\frac{1}{\kappa}t_*)+\text{K}_*)$ is used repeatedly in this paper:

*If* t *is such that* $\text{K}(t) \geq 4\kappa((1+\text{K}(\frac{1}{\kappa}t_*)+\text{K}_*)$ *then* $\text{K}(s) \geq \kappa(1+\text{K}(\frac{1}{\kappa}t_*)+\text{K}_*)$ *for all* $s \leq t$.

(5.17)

This also follows from Lemma 5.1. The contexts for the (mostly implicit) appeals to (5.17) is this: If a number $\kappa \geq 16$ is given and if it has been established that a certain property holds at any time s provided that $\text{K}(s) \geq \kappa(1+\text{K}(\frac{1}{\kappa}t_*)+\text{K}_*)$ (perhaps subject to other constraints), then that property will automatically hold at all $s \leq t$ (subject to those other constraints) if $\text{K}(t)$ is not less than $4\kappa((1+\text{K}(\frac{1}{\kappa}t_*)+\text{K}_*)$.

A remarks about Requirements b) and c): The $\text{N}(t) < 2$ requirement can be replaced by any a priori upper bound on N at the expense of changing $\kappa$. Note also that if $\text{N}(t) > 2$, then $\text{K}(s) \leq (\frac{t}{s})^2\text{K}(t)$ for $s > t$ until s is such that $\text{N}(s) = 2$. Likewise, the required lower bound for the $\{t\} \times Y$ integral of $\langle \mathfrak{a}, \mathcal{B}\rangle$ can be changed to be any fixed, non-positive multiple of $\text{K}^2$ at the expense of changing $\kappa$.

With regards to the first bullet of the lemma: This bullet and Requirement a) are consistent only in the event that

$$t \leq \frac{1}{\kappa^2(1+\text{K}(\frac{1}{\kappa}t_*)+\text{K}_*)} . \tag{5.18}$$

With regards to the second bullet of the lemma: It is not useful if $\text{N} \leq 0$.

***Proof of Lemma 5.8***: The proof of the lemma has two parts.

*Part 1*: Fix $\delta \in (0, \frac{1}{16}]$ and suppose for the moment that $t \in (0, \delta t_*]$. Assume in what follows that $\text{K}(t) \geq \frac{1}{\sqrt{\delta}}(1+\text{K}(\delta t_*) + \text{K}_*)$ (this requires $t < t_\Delta$); and assume that the $\{t\}\times Y$ integral of $\langle \mathfrak{a}, \mathcal{B}\rangle$ is greater than $-\text{K}^2(t)$. Given these assumptions, then the $s = \delta t_*$ version of Lemma 5.7's second bullet leads to this:

$$\int_{\{t\}\times Y} |\mathfrak{t}| \leq c_0\sqrt{\delta}\,\text{K}^2(t) . \tag{5.19}$$

Granted (5.19), and granted that $\delta \leq c_0^{-1}$, then the first bullet of (2.16) implies that

$$\text{K}^2(t) \leq 4 \int_{\{t\}\times Y} |\mathcal{Q}|^{2/3} . \tag{5.20}$$

As a consequence, a bound on $\text{K}(t)$ follows given a bound for the integral of $|\mathcal{Q}|^{2/3}$.

With regards to (5.20): The bound holds if $\delta \le \delta_0$ with $\delta_0 = \frac{1}{16} c_0^{-1}$. This said, take $\delta$ to be $\frac{1}{100}\delta_0$ so that (5.20) holds if $\mathrm{K}(t) \ge \frac{1}{2\sqrt{\delta}}(1 + \mathrm{K}(\delta t_*) + \mathrm{K}_*)$. With $\delta$ so defined, then constrain t so that $t \in (0, \delta t_*]$ and $\mathrm{K}(t) \ge \frac{1}{\sqrt{\delta}}(1 + \mathrm{K}(\delta t_*) + \mathrm{K}_*)$.

*Part* 2: The condition from Part 1 to the effect that $\mathrm{K}(t) \ge \frac{1}{\sqrt{\delta}}(1 + \mathrm{K}(\delta t_*) + \mathrm{K}_*)$ can not hold if $t > t_\Delta$ so t must be less than $t_\Delta$. More to the point, t must be less than $\frac{100}{101} t_\Delta$. To see why, suppose that this were not the case. Since $\mathrm{N}(t) < 2$ (by assumption), it follows from Lemma 5.3 that $\mathrm{N} \le 4$ between t and $t_\Delta$. As a consequence $\mathrm{K}(t_\Delta) \ge 4\mathrm{K}(t)$. On the otherhand, $\mathrm{K}(t_\Delta) \le \mathrm{K}_*$ so $\mathrm{K}_* \ge 4\mathrm{K}(t)$. This is impossible if $\delta > 16$.

Because $|\mathrm{N}(t)| \le 2$, the bound on $|\mathrm{N}|$ from the first bullet of Lemma 5.2 requires a time $s \in [t, \frac{101}{100} t]$ such that

$$\int_{\{s\}\times Y} |[\mathfrak{a}_*, \mathfrak{a}]|^2 \le c_0 \tfrac{1}{t^2} \mathrm{K}^2(t) \ . \tag{5.21}$$

Meanwhile, it follows from this first bullet of Lemma 5.2 that $\mathrm{N}(s) \le 4$ on $[t, s]$. Thus, because $s \le \frac{101}{100} t$, and because the inequality $\mathrm{K}(t) \ge \frac{1}{\sqrt{\delta}}(1 + \mathrm{K}(\delta t_*) + \mathrm{K}_*)$ holds, it follows that $\mathrm{K}(s) \ge \frac{1}{2\sqrt{\delta}}(1 + \mathrm{K}(\delta t_*) + \mathrm{K}_*)$ (integrate (2.11) from t to s). As a consequence, (5.19) can be invoked with t replaced by s on both sides. Likewise for (5.20). Meanwhile, (5.21) holds with t replaced by s on its right hand side.

With the preceding in mind, note that

$$\int_{\{s\}\times Y} |\mathcal{Q}|^{2/3} \le c_0 \mathrm{K}(s)^{2/3} \Big( \int_{\{s\}\times Y} |[\mathfrak{a}_*, \mathfrak{a}]|^2 \Big)^{1/3} , \tag{5.22}$$

and therefore (because of (5.21) with t replaced by s on the right hand side),

$$\int_{\{s\}\times Y} |\mathcal{Q}|^{2/3} \le c_0 \tfrac{1}{s} \mathrm{K}(s) \tag{5.23}$$

Granted this bound, then the t = s version of (5.20) leads to the bound $\mathrm{K}(s) \le c_0 \frac{1}{s}$. And, because $\mathrm{N} \le 4$ on $[t, s]$ and $s \le \frac{101}{100} t$, this in turn implies that $\mathrm{K}(t) \le c_0 \frac{1}{t}$.

**e) Bounds for times downstream from t**

This section first states and then proves a lemma that gives bounds for $\mathrm{K}$ at times from t to $\mathcal{O}(\sqrt{t})$ given only that certain bounds hold at time t. By way of notation, $c_*$ is used in the upcoming lemma and subsequently to denote Lemma 5.8's version of $\kappa$.

**Lemma 5.9**: *There exists* $\kappa > 1$ *with the following significance: Fix a Nahm pole solution* $(A, \mathfrak{a})$*; and fix* $t \in (t^-, \frac{1}{\kappa} t_*]$ *where Requirements* a), b) *and* c) *below are met.*

a) $K(t) > 4c_*(1 + K(\frac{1}{c_*} t_*) + K_*)$ ,

b) $N(t) \le 2$ ,

c) $\int_{\{t\}\times Y} \langle \mathfrak{a}, B_A \rangle > -K^2(t)$ .

*With* t *as above, fix a non-negative* $\Theta$ *so that* $\int_{\{t\}\times Y} \langle \mathfrak{a}, \mathcal{B} \rangle \ge -\frac{\Theta}{t}$. *Then*

- $K(s) \le \frac{2c_*}{s}$
- $\int_{\{s\}\times Y} \langle \mathfrak{a}, \mathcal{B} \rangle \ge -(\Theta + \kappa\, c_*^{\,2}) \frac{1}{t}$ .

*for all* $s \in [t, t_\Theta]$ *with* $t_\Theta$ *being no less than the smaller of* $\frac{1}{\kappa} t_*$, $\frac{c_*}{\sqrt{(\Theta + \kappa c_*^2)}} \sqrt{t}$ *and the largest time where* $K(\cdot) = 4c_*(1 + K(\frac{1}{c_*} t_*) + K_*)$.

To be sure, there is a factor of $\frac{1}{s}$ on the right hand side of the first bullet, but not on the right hand side of the second; the factor there is $\frac{1}{t}$ which is not $\frac{1}{s}$.

***Proof of Lemma 5.9***: The proof has three parts

*Part 1*: Fix $\Theta$ to be greater than $c_*$. By virtue of Lemma 5.8 (see (5.17)), one has $K(t) < \frac{c_*}{t}$ ; so there is a maximal $t' \in (t, \frac{1}{8} t_*]$ such that $K(t') \ge 4c_*(1 + K(\frac{1}{c_*} t_*) + K_*)$ and such that $K(s) \le \frac{2c_*}{s}$ for all $s \in [t, t']$. Let $t_{\ddagger}$ denote this maximal time. Assume henceforth that $K(t_{\ddagger})$ is strictly greater than $4c_*(1 + K(\frac{1}{c_*} t_*) + K_*)$ because there is nothing to prove otherwise.

*Part 2*: If $N(t_{\ddagger}) \le 2$ and if the $\{t_{\ddagger}\}\times Y$ integral of $\langle \mathfrak{a}, \mathcal{B} \rangle$ is greater than $-K^2(t_{\ddagger})$, then Lemma 5.8 can be invoked to see that $K(t_{\ddagger}) \le \frac{c_*}{t_{\ddagger}}$ which is less than $\frac{2c_*}{t_{\ddagger}}$. Therefore, one or the other of these two conditions can't hold at $t_{\ddagger}$. As explained directly, it can't be that $N(t_{\ddagger}) > 2$. In fact, $N(t_{\ddagger})$ can't be greater than 1; because if $N(t_{\ddagger})$ were greater than 1, then $N$ would be greater than 1 on some interval of the form $[(1 - \varepsilon) t_{\ddagger}, t_{\ddagger}]$ for some positive $\varepsilon$. And, if this were true, and if s were in this interval, then $K(s)$ would obey $K(s) > \frac{2c_*}{s}$ .

*Part 3*: Granted the preceding, the question now is this: Can the $\{t_{\ddagger}\}\times Y$ integral of $\langle \mathfrak{a}, \mathcal{B} \rangle$ be as large as $-K^2(t_{\ddagger})$ if $t_{\ddagger}$ is less than $\mathcal{O}(\sqrt{t})$? To answer this, consider the second bullet of Lemma 5.5: If the $\{t\} \times Y$ integral of $\langle \mathfrak{a}, \mathcal{B} \rangle$ is greater than $-\frac{\Theta}{t}$, and if $s \in [t, t_\Theta]$, then the second bullet of Lemma 5.5 (with s there the same as s here) and the fact that $K(s)$ for $s \in [t, t_\Theta]$ is at most $\frac{2c_*}{s}$ implies that

$$\int_{\{s\}\times Y} \langle \mathfrak{a}, \mathcal{B}\rangle \geq -(\Theta + c_0 c_*^2)\tfrac{1}{t} \ .$$

(5.26)

This holds for $s = t_‡$ in particular, and since $\kappa^2(t_‡)$ is $(\frac{2c_*}{t_‡})^2$, and since the $\{t_‡\}\times Y$ integral of $\langle \mathfrak{a}, \mathcal{B}\rangle$ is $-\kappa^2(t_‡)$, it follows that $t_‡$ can't be less than $\frac{2c_*}{\sqrt{(\Theta + c_0 c_*^2)}}\sqrt{t}$.

## 6. Where the $\{t\} \times Y$ integral of $\langle \mathfrak{a}, \mathcal{B}\rangle$ is non-negative (Part 1)

Let $\Omega^+$ denote the set of times $t \in [t^-, \frac{1}{8} t_*]$ where the integral of $\langle \mathfrak{a}, \mathcal{B}\rangle$ on $\{t\} \times Y$ is non-negative; and let $\Omega^-$ denote the set of $t \in [t^-, \frac{1}{8} t_*]$ where the integral of $\langle \mathfrak{a}, \mathcal{B}\rangle$ on $\{t\} \times Y$ is negative. (Remember that $t_*$ is an $\mathcal{O}(1)$ time that is defined without regard for any given Nahm pole solution. It is $\frac{1}{8\kappa_*^4}$ with $\kappa_*$ denoting Lemma 5.1's version of $\kappa$.) This section and the next (which is Section 7) derive bound for the integrals of $|\mathfrak{c}|$ and $|\mathfrak{c}|^2$ on $\{t\} \times Y$ for $t > t^-$ and $t \in \Omega^+$. These sections also derive bounds for integrals of $|B_A|^2$, $|E_A|^2$, $|\nabla_A{}^{\perp}\mathfrak{a}|^2$ and $|\mathfrak{c}|^2$ on domains of the form $I\times Y$ with $I$ being an interval in $\Omega^+$. The proposition below summarizes the salient observations from this section and Section 7.

With regards to notation: Given an interval $[t_0, t_1]$, the proposition uses $E_0$ to denote the $(0, t_0] \times Y$ integral of $|B_A|^2 + |E_A|^2 + |\nabla_A{}^{\perp}\mathfrak{a}|^2$. Also, to say that $t_1$ is a boundary point of $\Omega^-$ is to say that the $\{t_1\}\times Y$ integral of $\langle \mathfrak{a}, \mathcal{B}\rangle$ is zero.

A final note about notation: Lemma 5.8's version of $\kappa$ is denoted again by $c_*$.

**Proposition 6.1**: *There exists* $\kappa > 8$ *with the following significance: Fix a Nahm pole solution* $(A, \mathfrak{a})$*, and let* $[t_0, t_1]$ *denote an interval in the corresponding* $\Omega^+$*. Assume that* $t_0$ *is such that* $N(t_0) < 2$ *and* $K(t_0) \geq \frac{1}{100\, t_0}$*; and that* $t_1$ *is less than* $\frac{1}{\kappa} t_*$ *and is such that* $K(\sqrt{t_1}) \geq 16c_*(1 + K(\frac{1}{c_*} t_*) + \kappa_*)$. *Distinguish the three cases listed below.*

- CASE 1: *The time* $t_1$ *is a boundary point of* $\Omega^-$.
- CASE 2: *There is a boundary point of* $\Omega^-$ *between* $t_1$ *and* $\frac{1}{\kappa}\sqrt{t_1}$.
- CASE 3: *The time* $t_1$ *is not a boundary point of* $\Omega^-$ *and* $K(\sqrt{t_1}) = 16c_*(1 + K(\frac{1}{c_*} t_*) + \kappa_*)$.

*Fix* $z_* > \kappa$ *so that* $N(t_0) < 1 + z_* t_0^2$. *There exists* $\lambda > \kappa$ *which can depend on* $z_*$ *but is otherwise independent of* $t_0$, $t_1$ *and* $(A, \mathfrak{a})$*; this number* $\lambda$ *appears below.*

A) *In all cases, the inequalities listed below hold at each* $t \in [t_0, t_1]$.

- ○ $\displaystyle\int_{\{t\}\times Y} |\mathfrak{c}| < \lambda\Big(1 + \int_{\{t_0\}\times Y} (|c| + |\mathfrak{c}^-|)\Big)$ .
- ○ $\displaystyle\int_{\{t\}\times Y} |\mathfrak{c}|^2 < \tfrac{\lambda}{t}\Big(1 + \Big|\int_{\{t_0\}\times Y} c\,\Big|\Big)$ .

B) *The inequalities listed below hold in* CASE 1:

- $\displaystyle\int_{\{t_1\}\times Y} (|c| + |\mathfrak{c}^-|) < \lambda(t_1 - t_0)\Big(1 + \big(\int_{\{t_0\}\times Y} |c|\,\big)^2\Big) + \lambda E_0 + \int_{\{t_0\}\times Y} (|c| + |\mathfrak{c}^-|)\,.$
- $\displaystyle\int_{[t_0,t_1]\times Y} (|B_A|^2 + |E_A|^2 + |\nabla_A^{\perp}\mathfrak{a}|^2 + |\mathfrak{c}|^2) < \lambda(t_1 - t_0)\Big(1 + \big(\int_{\{t_0\}\times Y} (|c| + |\mathfrak{c}^-|)\,\big)^2\Big) + \lambda E_0\;.$

C) *The inequality below hold in* CASE 2 *and* CASE 3:

$$\int_{[t_0,t_1]\times Y} (|B_A|^2 + |E_A|^2 + |\nabla_A^{\perp}\mathfrak{a}|^2 + |\mathfrak{c}|^2) < \lambda\Big(1 + \big(\int_{\{t_0\}\times Y} (|c| + |\mathfrak{c}^-|)\,\big)^2 + E_0\Big)$$

D) *In all cases, the function* $N$ *at each* $t \in [t_0, t_1]$ *obeys* $N(t) \le 1 + z_* t^2$.

As noted, Proposition 6.1 is a summary of subsequent lemmas: The assertions in Part A follow from Lemma 7.4. The first assertion from Part B is from Lemma 7.6 and the second is from Lemma 7.5. The assertion of Part C is from Lemma 7.4 and the top bullet of Lemma 7.5 for $t = t_0$ and for s chosen as follows: In CASE 2, take s to be the point in $[t_1, \frac{1}{\kappa}\sqrt{t_1}]$ where the $\{\cdot\}\times Y$ integral of $\langle \mathfrak{a}, B_A\rangle$ is zero; and in CASE 3, take s to be the point $\sqrt{t_1}$. The assertion in Part D about $N(t_0)$ follows from Lemma 6.2.

## a) Bounds for N and K on $\Omega^+$

The first lemma in this subsection concerns the function $N$; it says in effect that $N$ can hardly increase across an interval in $\Omega^+$ where it 1 or more unless $K$ is very small. By way of notation, the lemma uses $(N - 1)_+$ to denote the maximum of 0 and $N - 1$. The lemma also uses $R_*$ to denote the norm of the Ricci curvature tensor.

**Lemma 6.2**: *There exists* $\kappa > 8$ *with the following significance: Suppose that* $(A, \mathfrak{a})$ *is a Nahm pole solution; and then suppose that* $[t_0, t_1] \subset \Omega^+$. *Also assume that* $N(t_0) \le 2$ *and that* $K(t) \ge c_* t$ *for all* $t \in [t_0, t_1]$. *Then* $N(t) < 2$ *on the whole of* $[t_0, t_1]$. *Moreover,*

$$(N(t) - 1)_+ \le (N(t_0) - 1)_+ + \tfrac{1}{2} R_* (t^2 - t_0^2) + 8R_* \int_{t_0}^{t} \frac{1}{K(\cdot)}\;.$$

*for all* $t \in [t_0, t_1]$.

The proof comes momentarily, first a remark: If $K(t') \ge 4c_*(1 + K(\frac{1}{c_*} t_*) + K_*)$ for any $t' \in [t_1, \frac{1}{c_*} t_*]$, then by virtue of (5.17), the bound $K(\cdot) \ge c_*(1 + K(\frac{1}{c_*} t_*) + K_*)$ will hold on $[t_0, t_1]$ which implies in particular that $K(t) \ge c_* t$ on $[t_0, t_1]$.

***Proof of Lemma 6.2***: The assertions will be seen to follow from the identity in (5.9) for the derivative of $N$. To start, take (5.9) and keep only the term $\frac{N(1-N)}{t}$ on the right hand side and the two terms that are not manifestly negative the term with the Ricci curvature and the term with the integral of $\langle \mathfrak{a}, B_A\rangle$. For the first: Since $|\langle \mathrm{Ric}, \langle \mathfrak{a} \otimes \mathfrak{a}\rangle|$ is no greater than $R_*|\mathfrak{a}|^2$, the term with the Ricci curvature in (5.9) is no larger than $R_* t$. For the second: The $\{t\} \times Y$ integral of $\langle \mathfrak{a}, B_A\rangle$ for $t \in \Omega^+$ is no smaller than $-R_* K$ because the $\{t\} \times Y$ integral of $\langle \mathfrak{a}, \mathcal{B}\rangle$ is non-negative. As a consequence, the term in (5.9) with $\langle \mathfrak{a}, B_A\rangle$ is no larger than $4R_* \frac{N}{K}$. Thus, the identity in (5.9) leads to this inequality:

$$\frac{d}{dt} N \leq \frac{N(1-N)}{t} + R_* t + 4R_* \frac{N}{K} \ . \tag{6.1}$$

The inequality in (6.1) says (in part) that $\frac{d}{dt} N$ is negative when $N \geq 2$ unless $K \leq 4 R_* t$, so if $K$ is greater than $c_* t$ at each $t \in [t_0, t_1]$, then $N$ is less than 2 on the whole of $[t_0, t_1]$. And, supposing that $N$ is less than 2 on $[t_0, t_1]$, then integrating (6.1) where $1 \leq N$ leads directly to the inequality that is asserted by the lemma.

Looking ahead, there are two cases of interest with regards to this lemma: The first is when $N(t_0) < 1 + z t_0$ and $z$ is $\mathcal{O}(1)$ and the second when $N(t_0) \leq 1 + z_* t_0^2$ with $z_*$ being $\mathcal{O}(1)$. (This is the case for Proposition 6.1.) The next lemma talks about these two cases and the implications for the behavior of $K$.

**Lemma 6.3**: *Suppose that* $(A, \mathfrak{a})$ *is a Nahm pole solution. Let* $[t_0, t_1] \subset \Omega^+$ *denote an interval with* $t_0 < t_1 \leq \frac{1}{c_*} t_*$. *Suppose that* $K(t) \geq c_*(1 + K(\frac{1}{c_*} t_*) + K_*)$ *when* $t \in [t_0, t_1]$*; and suppose that* $N(t_0) < 2$. *Then* $K(t) \leq c_* \frac{1}{t}$ *when* $t \in [t_0, t_1]$. *Moreover:*

- *If* $z > 100 c_*$ *and* $N(t_0) \leq 1 + z t_0$, *then*
  - $N(t) \leq 1 + z t$ *for all* $t \in [t_0, t_1]$.
  - $K(t) \geq \frac{t_0}{t} e^{-z(t - t_0)} K(t_0)$ .
- *If* $\delta > 0$ *and* $K(t_0) \geq \delta \frac{1}{t}$ *; and if* $z_* \geq 100 c_*(1 + \frac{1}{\delta})$ *is such that both* $z_* t_1 \leq 1$ *and also* $N(t_0) \leq 1 + z_* t_0^2$, *then*
  - $N(t) \leq 1 + z_* t^2$ *for all* $t \in [t_0, t_1]$.
  - $K(t) \geq \frac{t_0}{t} e^{-\frac{1}{2} z_* (t^2 - t_0^2)} K(t_0)$ .

With regards to the assumptions: Keep in mind the note subsequent to Lemma 5.8 to the effect that the condition $K(t) \geq c_*(1 + K(\frac{1}{c_*} t_*) + K_*)$ requires that $t$ be less than Lemma 5.1's time $t_\Delta$. Also keep in mind that this condition holds on $[t_0, t_1]$ when there is a time $t' \in [t_1, t_\Delta]$ with $K(t')$ obeying the inequality $K(t') \geq 4c_*(1 + K(\frac{1}{c_*} t_*) + K_*)$.

***Proof of Lemma 6.3***: The bound $\mathrm{K}(t) \le c_* \frac{1}{t}$ follows directly from Lemma 5.8 given the lower bound assumption about $\mathrm{K}(t)$ and given that $\mathrm{N}(t_0) \le 2$. This is because Lemma 6.2 guarantees that $\mathrm{N}(t) \le 2$ on the whole of $[t_0, t_1]$. (Also: Remember that the $\{t\}\times Y$ integral of $\langle \mathfrak{a}, B_A\rangle$ is no smaller than $-\mathrm{R}_*\mathrm{K}$. To prove the lemma's first bullet: The bound on $\mathrm{N}$ follows from Lemma 6.2's inequality for $(\mathrm{N}(t) - 1)_+$ because $\frac{1}{\mathrm{K}(t)} \le \frac{1}{c_*}$. The lower bound on $\mathrm{K}$ follows from the upper bound on $\mathrm{N}$ using (2.11). To prove the lemma's second bullet: Use the Lemma's first bullet and the upper bound assumption about $\mathrm{K}(t_0)$ to see that $\frac{1}{\mathrm{K}(t)} \le \frac{4}{\delta} t$ for all $t \in [t_0, t_1]$. Now use this bound in Lemma 6.2's inequality for $(\mathrm{N} - 1)_+$ to get the asserted bound on $\mathrm{N}$. The bound for $\mathrm{K}$ follows from the bound for $\mathrm{N}$ using (2.11).

**b) An $\mathcal{O}(\frac{1}{t})$ upper bound for the integral of $\langle \mathfrak{a}, \mathcal{B}\rangle$ on $\{t\} \times Y$ when $t \in \Omega^+$**

The upcoming lemma about the $\{t\}\times Y$ integral of $\langle \mathfrak{a}, B_A\rangle$ uses $c_*$ to denote the version of $\kappa$ from Lemma 5.8.

**Lemma 6.4**: *There exists* $\kappa > c_*$ *which is independent of* $(A, \mathfrak{a})$ *and has the following significance: Let* $[t_0, t_1] \subset \Omega^+$ *denote an interval* $t_1 \le \frac{1}{\kappa} t_*$. *Suppose in addition that* $\mathrm{N}(t_0) < 2$ *and that* $\mathrm{K}(\sqrt{t_1}) \ge 4c_*(1 + \mathrm{K}(\frac{1}{c_*} t_*) + \mathrm{K}_*)$. *Then* $\int_{\{t\}\times Y} \langle \mathfrak{a}, \mathcal{B}\rangle \le \kappa \frac{1}{t}$ *on* $[t_0, t_1]$.

***Proof of Lemma 6.4***: A preliminary remark: If $t_1 < c_0^{-1} t_*$ with $c_0 > c_*$, then $\mathrm{N}$ will be less than 2 on the whole of $[t_0, t_1]$. This is guaranteed by Lemma 6.3. With this understood, Lemma 5.9 can be employed here to see that $\mathrm{K}(t) \le c_0 c_* \frac{1}{t}$ on the interval $[t_0, \frac{1}{c_*}\sqrt{t_1}]$. With this bound in hand, invoke the $s = \frac{1}{4} t_*$ version of the second bullet of Lemma 5.5 with the integral of $\mathrm{K}^2$ that appears on the right hand side of that bullet split into the contributions from the respective intervals $[t, \frac{1}{c_*}\sqrt{t_1}]$ and $[\frac{1}{c_*}\sqrt{t_1}, \frac{1}{4} t_*]$. The former contribution is no larger than $c_0 c_* \frac{1}{t}$ and the latter contribution is bounded by $c_0 \mathrm{K}^2(\frac{1}{c_*}\sqrt{t})$ using Lemma 5.1. With the preceding understood, then the second bullet of Lemma 5.5 leads to a $c_0 \frac{1}{t} + c_0 (\mathrm{K}^2(\frac{1}{c_*}\sqrt{t}) + \mathrm{K}_*^2)$ bound for the $\{t\}\times Y$ integral $\langle \mathfrak{a}, \mathcal{B}\rangle$. This bound implies the bound in the lemma because $\mathrm{K}_*$ is no greater than $c_0\, \mathrm{K}(\frac{1}{c_*}\sqrt{t})$ which in turn is at most $c_0 c_* \frac{1}{\sqrt{t}}$.

**c) An $\mathcal{O}(\frac{1}{t})$ upper bound for the integral of $|\mathfrak{t}|$ on $\{t\} \times Y$**

Remember that $\mathfrak{t}$ is the tensor $\langle \mathfrak{a} \otimes \mathfrak{a}\rangle - \frac{1}{3} \mathfrak{g} |\mathfrak{a}|^2$, the traceless part of $\langle \mathfrak{a} \otimes \mathfrak{a}\rangle$.

**Lemma 6.5**: *There exists* $\kappa > 16c_*$ *that is independent of the given Nahm pole solution* $(A, \mathfrak{a})$ *and has the following property: Let* $[t_0, t_1]$ *denote an interval in* $\Omega^+$ *with endpoints obeying* $N(t_0) < 2$ *and* $t_1 \leq \frac{1}{\kappa} t_*$ *and* $K(\sqrt{t_1}) \geq 16c_*(1 + K(\frac{1}{c_*} t_*) + K_*)$. *If* $t \in [t_0, t_1]$ *then*

$$\int_{\{t\}\times Y} |\mathfrak{t}| < \kappa \tfrac{1}{t} .$$

Here are two other points to note for the future: Supposing that $\kappa > c_*$ and that $t_1$ is such that $K(\sqrt{t_1}) \geq 16c_*(1 + K(\frac{1}{c_*} t_*) + K_*)$, then $K(t) \geq 4c_* (1 + K(\frac{1}{c_*} t_*) + K_*)$ when $t \leq \sqrt{t_1}$. This is by virtue of (5.17). Therefore, previous lemmas (Lemmas 5.8, 5.9 and 6.3 and 6.4) that require $K(t)$ to be greater than $4c_* (1 + K(\frac{1}{c_*} t_*) + K_*)$ or $c_*(1 + K(\frac{1}{c_*} t_*) + K_*)$ can be invoked when $t \leq \sqrt{t_1}$ if their other constraints are met.

***Proof of Lemma 6.5***: Suppose first that $t \in [t_0, t_1]$. Let $t_2$ denote the largest time greater than t such that $K(s) \leq \frac{2c_*}{s}$ for every $s \in [t, t_2]$ and such that $K(t_2) \geq 4c_*(1 + K(\frac{1}{c_*} t_*) + K_*)$. The $\Theta = 0$ version of Lemma 5.9 can be invoked to see that $t_2$ is no less than that lemma's $t_\Theta$. According to Lemma 5.9, this $t_\Theta$ is greater than $c_0^{-1}\sqrt{t}$ provided that $t \leq c_0^{-1} t_*$ and provided that $K(\sqrt{t}) \geq 4c_*(1 + K(\frac{1}{c_*} t_*) + K_*)$; and this is guaranteed to be the case if $t_1 \leq c_0^{-1} t_*$ and $K(\sqrt{t_1}) \geq 16c_*(1 + K(\frac{1}{c_*} t_*) + K_*)$. Assuming that these $t_1$ constraints are met, use this $t_2$ for the time s in Lemma 5.7's second bullet to see that

$$\int_{\{t\}\times Y} |\mathfrak{t}| \leq c_0 \left(\tfrac{1}{t} + K^2(t_2)\right) \quad . \tag{6.2}$$

To explain: The factor of $K^2(t_2)$ accounts for the $\{s\}\times Y$ integral of $|\mathfrak{t}|$ and also the $\{s\}\times Y$ integral of $\langle \mathfrak{a}, \mathcal{B} \rangle$ if the latter is positive. (It is bounded by $c_* K^2(t_2)$ courtesy of the version of the top bullet in Lemma 5.6 that has $t = t_2$ and $s = \frac{1}{8} t_*$.) The term on the right hand side of (6.5) proportional $\frac{1}{t}$ comes from the integrals of $K^2(\cdot)$ between t and s.

## d) Initial bounds for $\kappa^2 - \frac{3}{4t^2}$

Suppose in what follows that $[t_0, t_1] \subset \Omega^+$ and that $K(\sqrt{t_1}) \geq 16c_*(1 + K(\frac{1}{c_*} t_*) + K_*)$ and that $t_1$ is small enough (less than $\frac{1}{c_0 c_*} t_*$) so that Lemmas 5.8, 5.9 and 6.3, 6.4 can be invoked for $t \leq t_1$, and that Lemma 6.5 can be invoked for $t \leq t_1$ (if their other requirements are met). One further assumption is made here which is this: If $N(t_0) > 1$, then $N(t_0)$ can be written as $1 + z t_0$ with $z t_0 < 1$ (which implies that $N(t_0) \leq 2$). If $N(t_0) \leq 1$, set $z = 0$. (Think of $z$ as being at most $\mathcal{O}(1)$, which is what it will be for the cases relevant to the theorems in Section 1.)

The upcoming lemma introduces notation that will be used subsequently: What is denoted by $c_‡$ signifies the larger of the versions of the number $\kappa$ that appear in Lemmas 6.2–6.5.

**Lemma 6.6**: *Fix a Nahm pole solution* $(A, \mathfrak{a})$. *Suppose that* $[t_0, t_1] \in \Omega^+$ *with* $t_1 < \frac{1}{c_‡} t_*$ *and with* $\kappa(\sqrt{t_1})$ *being larger than* $16c_*(1 + \kappa(\frac{1}{c_*} t_*) + \kappa_*)$. *If these conditions are met, then there exists* $\kappa > 1$ *which, except for possible dependence on z (this is defined above), is independent of* $(A, \mathfrak{a})$ *and* $[t_0, t_1]$*; and it is such that* $\kappa^2 - \frac{3}{4t^2}$ *is at most* $\kappa \frac{1}{t}$ *on* $[t_0, t_1]$.

With regards to a lower bound: Lemma 6.3 lead directly to a lower bound for the function $\kappa^2 - \frac{3}{4t^2}$ on $[t_0, t_1]$ given a lower bound for $\kappa(t_0)$ and an upper bound for $N(t_0)$ of the form $N(t_0) \le 1 + z t_0$.

***Proof of Lemma 6.6***: The proof uses a convention whereby $c_1$ is a number greater than 1 that is allowed to depend implicitly on given upper bounds for number $z$. Otherwise, it is independent of the given Nahm pole solution $(A, \mathfrak{a})$ and $[t_0, t_1]$. The precise value of $c_1$ can be assumed to increase between successive appearances.

The desired upper bound $\kappa^2 - \frac{3}{4t^2}$ is obtained by exploiting the top bullet of (2.16) using what can be inferred about the size of the terms in this equation. The four parts that follow give the details.

*Part 1*: Suppose for the moment that $t \in [t_0, \frac{101}{100} t_1]$. The top bullet in (2.16) leads (using Lemma 6.5) to the inequality

$$\kappa^2 \le \frac{3}{4^{1/3}} \int_{\{t\}\times Y} |\mathcal{Q}|^{2/3} + c_0 \frac{1}{t} \,, \tag{6.3}$$

which in turn is less than

$$\kappa^2 \le \frac{3}{4^{1/3}} \Big( \int_{\{t\}\times Y} |\mathcal{Q}| \Big)^{2/3} + c_0 \frac{1}{t} \,. \tag{6.4}$$

This last inequality is then written as

$$\kappa^2 \le \frac{3}{4^{1/3}} \Big( \int_{\{t\}\times Y} \mathcal{Q} + \Delta \Big)^{2/3} + c_0 \frac{1}{t} \tag{6.5}$$

with $\Delta$ denoting the number

$$\Delta = \int_{\{t\}\times Y} |\mathcal{Q}| - \int_{\{t\}\times Y} \mathcal{Q} \,.$$

(6.6)

Meanwhile, by virtue of the second bullet in (2.9) and Lemma 6.4:

$$\int_{\{t\}\times Y} \mathcal{Q} \le \tfrac{NK^2}{3t} + c_0 \tfrac{1}{t} \;.$$

(6.7)

Use (6.7) in (6.5) to see that

$$K^2 \le \tfrac{3}{4^{1/3}} \left(\tfrac{NK^2}{3t} + \Delta\right)^{2/3} + c_0 \tfrac{1}{t} \;.$$

(6.8)

This inequality still leaves $N$ and $\Delta$ unspoken for.

*Part 2*: The bound $N(t) \le 1 + z t$ can be used in (6.11) if $t \in [t_0, t_1]$. Indeed, Lemma 6.3 leads to this bound for $t \in [t_0, t_1]$ because $N(t_0)$ is at most $1 + z t_0$ (the latter bound is assumed in the statement of the lemma).

To obtain a bound for $\Delta$, fix $s > t^-$ for the moment and let $W(s)$ denote the subset of points in $\{s\} \times Y$ where $\mathcal{Q} < 0$. The number $\Delta(s)$ is the difference between the integrals of $|\mathcal{Q}|$ and $\mathcal{Q}$ on $\{s\} \times Y$ which is twice the integral of $-\mathcal{Q}$ on $W(s)$. With this understood, then $\Delta$ is the function on $[t^-, \frac{1}{8} t_*] \times Y$ that is given by rule

$$\Delta(s) = -2 \int_{\{s\}\times W(s)} \mathcal{Q} \;.$$

(6.9)

(Thus, $\Delta$ is positive because $\mathcal{Q} < 0$ on $W(s)$.) The function $\Delta$ is piecewise differentiable with derivative given by the rule

$$\tfrac{d}{ds} \Delta = -2 \int_{\{s\}\times W(s)} \tfrac{\partial}{\partial s} \mathcal{Q} \;.$$

(6.10)

(There is no term from the change of the domain W because $\mathcal{Q} = 0$ on W's boundary.)

Meanwhile, the derivative of $\mathcal{Q}$ along the $(0, \infty)$ factor of $(0, \infty) \times Y$ can be computed using the middle bullet in (2.7):

$$\tfrac{\partial}{\partial s} \mathcal{Q} = -\tfrac{1}{2} |\mathfrak{a} \wedge \mathfrak{a}|^2 + *\langle B \wedge \mathfrak{a} \wedge \mathfrak{a} \rangle \;.$$

(6.11)

The preceding identity leads to the inequality

$$\frac{\partial}{\partial s}\mathcal{Q} \le \frac{1}{2}|B_A|^2 ;$$

(6.12)

and this with (6.10) implies that $\Delta$ obeys

$$\frac{d}{ds}\Delta \ge -\int_{\{s\}\times Y} |B_A|^2 \quad .$$

(6.13)

As explained in Part 3, this is the key inequality.

*Part* 3: Take $t \in [t_0, t_1]$ and integrate (6.13) on the domain $[t, \frac{1}{4}t_*]$ to see that

$$\Delta(t) \le \Delta(\tfrac{1}{4}t_*) + \tfrac{1}{2}\int_t^{\frac{1}{4}t_*} |B_A|^2 \quad .$$

(6.14)

The bound in the second bullet of Lemma 5.5 can be invoked with the help of Lemmas 5.1 and 5.9 to see that

$$\Delta(t) \le c_0 \kappa_*^2 + c_0 c_* \tfrac{1}{t} \ .$$

(6.15)

To elaborate: Lemma 5.9 says in part that $\kappa(s) \le 2c_* \frac{1}{s}$ when s is from the $[t, \frac{1}{c_\ddagger}\sqrt{t}]$ interval, and granted this, Lemma 5.1 says in part that $\kappa(\cdot) \le 2c_\ddagger^{\,2}\frac{1}{\sqrt{t}}$ on the $[\frac{1}{c_\ddagger}\sqrt{t}, \frac{1}{4}t_*]$ interval. With this understood, break the integral of $\kappa^2$ that appears on the right of the second bullet of Lemma 5.5 into the contribution from $[t, \frac{1}{c_\ddagger}\sqrt{t}]$ and from $[\frac{1}{c_\ddagger}\sqrt{t}, \frac{1}{4}t_*]$. Having done that, use the preceding bounds for $\kappa^2$ to obtain the bound in (6.15).

One last point with regards to $\Delta$: The bound in (6.15) implies in turn that

$$\Delta(t) \le c_0 \tfrac{1}{t}$$

(6.16)

because $\kappa(\frac{1}{c_\ddagger}\sqrt{t}) \ge \kappa_*$ by assumption, and because $\kappa(\frac{1}{c_\ddagger}\sqrt{t})$ is at most $c_0 c_\ddagger^{\,2}\frac{1}{\sqrt{t}}$.

*Part 4*: To finish the proof: If $t \in [t_0, t_1]$ and if $\kappa^2(t) \le \frac{3}{4t^2}$, then there is nothing to prove. Supposing that $\kappa^2 > \frac{3}{4t^2}$, then (6.8) and (6.16) and $\kappa$'s lower bound implies that

$$\kappa^2 \le (\tfrac{3}{4t^2})^{1/3}(1 + c_1 t)^{2/3}\kappa^{4/3} + c_0\tfrac{1}{t} \quad .$$

(6.17)

Dividing through by $\kappa^{4/3}$ and invoking the $\frac{\sqrt{3}}{2t}$ lower bound for $\kappa$ leads from (6.20) to this:

$$K^{2/3} \le \left(\tfrac{3}{4t^2}\right)^{1/3}(1 + c_1 t)^{2/3} \ . \tag{6.18}$$

This bound implies directly the bound in Lemma 6.6.

### e) Initial bounds for the integrals of $|\mathfrak{c}|^2$ and $|\mathfrak{c}|$ on $\{t\} \times Y$

Write $\mathfrak{a}$ as in (4.2) so as to define $\mathfrak{c}$ and its components as done in (4.3). The central lemma in this section states a preliminary bound for the function of t given by the integral of $|\mathfrak{c}|^2$ on $\{t\} \times Y$. The assumption is that t is in an interval $[t_0, t_1]$ from $\Omega^+$. The upcoming lemma defines $z$ as before, by writing $N(t_0)$ as $1 + z t_0$ if $N(t_0) > 1$. As before, the condition $z t_0 < 1$ is assumed (and so $N(t_0) \le 2$). If $N(t_0) \le 1$, then $z = 0$. Remember that $c_*$ is the value of $\kappa$ that appears Lemma 5.8. The number $c_*$ does not depend on any given Nahm pole solution. Meanwhile, $c_\ddagger$ in what follows is an $\mathcal{O}(c_0)$ number that is larger than the values of $\kappa$ in Lemmas 6.2–6.5.

**Lemma 6.7**: *Let* $(A, \mathfrak{a})$ *denote a given Nahm pole solution and let* $[t_0, t_1]$ *denote an interval in* $\Omega^+$ *such that* $t_1 \le \frac{1}{c_\ddagger} t_*$ *and such that* $K(\sqrt{t_1}) \ge 16 c_*(1 + K(\frac{1}{c_*} t_*) + K_*)$. *There exists* $\kappa > 1$ *which is independent of the given Nahm pole solution and* $[t_0, t_1]$ *except for possible dependence on z; it has the following significance: Suppose that* $t \in [t_0, t_1)$ *is given as is a later time* $s \in (t, t_1]$. *Then*

- $\displaystyle\int_{[t,s] \times Y} (|\mathfrak{c}^+|^2 + |\mathfrak{c}^-|^2) \le \kappa\Big(\ln(\tfrac{s}{t}) + 1 + \Big|\int_{\{t\}\times Y} c\,\Big|\Big)$ .
- $\displaystyle\int_{\{s\}\times Y} |\mathfrak{c}|^2 \le \kappa\, \tfrac{1}{s}\Big(\ln(\tfrac{s}{t}) + 1 + \Big|\int_{\{t\}\times Y} c\,\Big|\Big)$.
- $\displaystyle\int_{\{s\}\times Y} |\mathfrak{c}| \le \kappa\Big(\ln(\tfrac{s}{t}) + 1 + \int_{\{t\}\times Y} (|c| + |\mathfrak{c}^-|)\Big)$ .

Notice that the $\{t\} \times Y$ integrals that appears on the right hand side of these inequalities do not depend on $\mathfrak{c}^+$.

***Proof of Lemma 6.7***: The convention in this proof allows the values of the subsequent incarnations of $c_0$ to depend implicitly on upper bounds for $z$. If this is $\mathcal{O}(1)$, then so is $c_0$. The incarnations of $c_0$ are otherwise independent of the given Nahm pole solution $[A, \mathfrak{a}]$ and the interval $[t_0, t_1]$.

The proof of the lemma has three parts.

*Part 1*: This part proves the following bound:

$$\int_{\{s\}\times Y} c \; + \int_{[t,s]\times Y} (|\mathfrak{c}^+|^2 + |\mathfrak{c}^-|^2) \le c_0\big(\ln(\tfrac{s}{t}) + 1 + \big|\int_{\{t\}\times Y} c\,\big|\big) . \tag{6.19}$$

The subsequent parts explain how this leads to the assertions made by the lemma. To prove (6.19), use the top bullet in (4.4) to write (4.9)'s top equation as

$$\sqrt{3}\,\tfrac{\partial}{\partial t}\, c + 3|\mathfrak{c}^+|^2 + |\mathfrak{c}^-|^2 = \langle \sigma B_A\rangle + 2\,(|\alpha|^2 - \tfrac{3}{4t^2}) . \tag{6.20}$$

Integrate this equation on any given $\{t\} \times Y$ to obtain the equation

$$\sqrt{3}\,\tfrac{d}{dt} \int_{\{t\}\times Y} c \; + \int_{\{t\}\times Y} (3\,|\mathfrak{c}^+|^2 + |\mathfrak{c}^-|^2) \; = \int_{\{t\}\times Y} \langle \sigma, B_A\rangle \; + 2\,(\kappa^2 - \tfrac{3}{4t^2}) . \tag{6.21}$$

Then, integrate the preceding identity on the given interval [t, s]. The resulting integral identity with Lemma 6.6's bound $\kappa^2 - \frac{3}{4t^2} \le c_0 \frac{1}{t}$ leads to the following inequality:

$$\sqrt{3} \int_{\{s\}\times Y} c \; + \int_{[t,s]\times Y} (3\,|\mathfrak{c}^+|^2 + |\mathfrak{c}^-|^2) \; \le c_0 \sqrt{s}\,\Big( \int_{[t,s]\times Y} |B_A|^2 \Big)^{1/2} + \sqrt{3} \int_{\{t\}\times Y} c \; + c_{0*}\ln(\tfrac{s}{t}) . \tag{6.22}$$

Consider the integral of $|B_A|^2$ that appears in (6.22): Invoke the top bullet of Lemma 5.5 and Lemma 6.4 (with t replaced by s) to see that

$$\int_{[t,s]\times Y} |B_A|^2 \; \le c_0 \tfrac{1}{s} \; + c_0 \int_{[t,s]\times Y} (|\mathfrak{c}^+|^2 + |\mathfrak{c}^-|^2) \; . \tag{6.23}$$

Using this bound leads to the bound in (6.22) because the *square root* of the integral of $|B_A|^2$ appears in (6.23).

*Part 2*: The first bullet of Lemma 6.7 follows directly from (6.19) if the $\{s\} \times Y$ integral of $c$ is non-negative. The first bullet also follows from (6.19) if the $\{s\} \times Y$ integral of $c$ is negative because (4.4) implies the identity

$$\kappa^2(s) - \tfrac{3}{4s^2} \; = - \tfrac{\sqrt{3}}{s} \int_{\{s\}\times Y} c \; + \int_{\{s\}\times Y} |\mathfrak{c}|^2 \quad , \tag{6.24}$$

and Lemma 6.6's bound $\kappa^2(s) - \frac{3}{4s} \le c_0 \frac{1}{s}$ implies that the $\{s\}\times Y$ integral of $c$ is greater than $-c_0$. The second bullet of the lemma follows from (6.19) and (6.24) and Lemma 6.6 when the $\{s\} \times Y$ integral of $c$ is non-negative; and it follows from just (6.24) and Lemma 6.6 when said integral is negative.

*Part 3*: The third bullet of the lemma holds if it holds separately for $|c|$, $|\mathfrak{c}^-|$ and $|\mathfrak{c}^+|$. The proofs of the latter bounds follow in order. To bound the $\{s\} \times Y$ integral of $|c|$, suppose for the moment that t is any given point in $[t^-, \frac{1}{8} t_*]$. Multiply both sides of the equation in the top bullet of (4.9) by the sign of $c$ (where $c \neq 0$) and then integrate over $\{t\} \times Y$. Having done so, multiply both sides of the resulting identity by $\sqrt{3}\, t^2$ to obtain

$$\frac{d}{dt}\Big(t^2 \int_{\{t\}\times Y} |c|\Big) \le c_0\, t^2 \int_{\{t\}\times Y} |B_A| \;+\; 2t^2 \int_{\{t\}\times Y} |\mathfrak{c}|^2 \;. \tag{6.25}$$

Integrate (6.25) from Lemma 6.7's version of t to its version of s to see that

$$\sqrt{3} \int_{\{s\}\times Y} |c| \le \sqrt{3}\,\frac{t^2}{s^2} \int_{\{t\}\times Y} |c| + c_0 (s-t)^{1/2} \Big( \int_{[t,s]\times Y} |B_A|^2 \Big)^{1/2} + \frac{2}{t^2} \int_{[t,s]\times Y} x^2\, |\mathfrak{c}|^2\, dx \;. \tag{6.26}$$

Use the bound from the second bullet of Lemma 6.7 for the $\{x\} \times Y$ integral of $|\mathfrak{c}|^2$ at each $x \in [t, s]$ to bound the right most integral in (6.26). And, bound the $|B_A|^2$ integral on the right hand side of (6.26) using (6.23) and the top bullet's bound for the $[t, s] \times Y$ integral of $|\mathfrak{c}^+|^2 + |\mathfrak{c}^-|^2$. This leads to the desired bound for the $\{s\} \times Y$ integral of $|c|$.

The desired bound for the $\{t\} \times Y$ integral of $|\mathfrak{c}^-|$, can be obtained using much the same argument as for $|c|$ by exploiting the equation in the third bullet of (4.9). In particular, this equation implies the following one:

$$\frac{d}{dt}\Big(t \int_{\{t\}\times Y} |\mathfrak{c}^-|\Big) \le c_0\, t \int_{\{t\}\times Y} |B_A| \;+ 2t \int_{\{t\}\times Y} |\mathfrak{c}|^2 \;; \tag{6.27}$$

and integrating the latter from t to s leads to the desired bound.

Now consider the $\{s\} \times Y$ integral of $|\mathfrak{c}^+|$. The identity in (4.4) leads to this bound:

$$\frac{1}{s}\,|\mathfrak{c}^+| \le |\mathfrak{t}| + c_0\, |\mathfrak{c}|^2 \,; \tag{6.28}$$

and this bound, when integrated leads to the lemma's bound because Lemma 6.5 bounds the $\{s\} \times Y$ integrals of $|\mathfrak{t}|$ by $c_0 \frac{1}{s}$ and the second bullet of this lemma bounds the $\{s\} \times Y$ integrals of $|\mathfrak{c}|^2$ by $c_0(\frac{1}{s}(\ln(\frac{s}{t}) + 1 + |\int_{\{t\}\times Y} c\,|)$.

## 7. Where the $\{t\} \times Y$ integral of $\langle \mathfrak{a}, \mathcal{B}\rangle$ is non-negative (Part 2)

This section revises the bounds that are supplied by the lemmas in Section 6. The notation used in this section has $[t_0, t_1]$ being an interval in $\Omega^+$ which is subject to certain contraints with regards to its endpoints. As in the previous section, $z$ is defined by writing $N(t_0)$ as $1 + z t_0$ if $N(t_0) > 1$; and it is assumed that $z t_0 < 1$ so that $N(t_0) \leq 2$. If $N(t_0) \leq 1$, then $z$ is set equal to zero. The number $c_*$ is the version of $\kappa$ from Lemma 5.8; it is independent of any given Nahm pole solution. The number $c_\ddagger$ is an $\mathcal{O}(1)$ number that is larger than the values of $\kappa$ in Lemmas 6.2–6.5.

### a) Revised bounds for the integrals of $\langle \mathfrak{a}, \mathcal{B}\rangle$ and $|\mathfrak{t}|$ on $\{t\} \times Y$

The lemma in this subsection supplies new upper bounds for the $\{t\}\times Y$ integrals of $\langle \mathfrak{a}, \mathcal{B}\rangle$ and $|\mathfrak{t}|$ when t is from $[t_0, t_1]$.

**Lemma 7.1**: *Fix a Nahm pole solution* $(A, \mathfrak{a})$ *and let* $[t_0, t_1]$ *denote an interval in the corresponding version of* $\Omega^+$ *with* $N(t_0) < 2$*; and with* $t_1 \leq \frac{1}{c_\ddagger} t_*$ *and such that* $K(\sqrt{t_1}) \geq 16 c_*(1 + K(\frac{1}{c_*} t_*) + K_*)$. *There exists* $\kappa > 1$ *which is independent of the given Nahm pole solution and* $[t_0, t_1]$ *except possible dependence on z, and which has the following significance: Fix* t *and* s *from* $[t_0, t_1]$ *with* $s > t$.

- $\int_{\{t\}\times Y} \langle \mathfrak{a}, \mathcal{B}\rangle + \frac{1}{4} \int_{[t,s]\times Y} (|B_A|^2 + |E_A|^2 + |\nabla_A^\perp \mathfrak{a}|^2) \leq \int_{\{s\}\times Y} \langle \mathfrak{a}, \mathcal{B}\rangle + \kappa(\ln(\frac{s}{t}) + 1 + |\int_{\{t\}\times Y} c|)$.
- $\int_{\{t\}\times Y} |\mathfrak{t}| \leq \int_{\{s\}\times Y} |\mathfrak{t}| + \kappa \frac{1}{\sqrt{t}} (\int_{\{s\}\times Y} \langle \mathfrak{a}, \mathcal{B}\rangle + \ln(\frac{s}{t}) + 1 + |\int_{\{t\}\times Y} c|)^{1/2}$ .

***Proof of Lemma 7.1***: To prove the top bullet: Invoke the top bullet of Lemma 5.5 and then use the bounds from the top bullet of Lemma 6.7 for the integral $|\mathfrak{c}^+|^2 + |\mathfrak{c}^-|^2$ that appears on the right hand side of the inequality in the top bullet of Lemma 5.5.

To prove the second bullet: Invoke the top bullet of Lemma 5.7. When doing this, use Lemma 6.3 to bound the integral of $K^2$ that appears on the right hand side of the inequality in the top bullet of Lemma 5.7 by $c_0 \frac{1}{t}$; and, use the top bullet of Lemma 6.7 to bound the integral of $|\mathfrak{c}^+|^2 + |\mathfrak{c}^-|^2$ that appears on the right hand side of that same inequality from Lemma 5.7.

### b) A new upper bound for $K^2 - \frac{3}{4t^2}$

The lemma in this section asserts an upper bound for the $K^2 - \frac{3}{4t^2}$ when $t \in [t_0, t_1]$ that can be significantly smaller than the bound in Lemma 6.6. By way of new notation: The lemma introduces a non-negative number it denotes by $z_*$ which is defined to be zero if $N(t_0) \le 1$, and it is defined otherwise by writing $N(t_0)$ as $1 + z_* t_0^2$. It is assumed in what follows that $z_* t_0^2 < 1$. (But, in applications, $z_*$ will be $\mathcal{O}(1)$.) The numbers $c_*$ and $c_\ddagger$ are defined as before.

**Lemma 7.2**: *Fix a Nahm pole solution* $(A, \mathfrak{a})$ *and let* $[t_0, t_1]$ *denote an interval in the corresponding version of* $\Omega^+$ *with* $t_0$ *such that both* $N(t_0) < 2$ *and* $K(t_0) \ge \frac{1}{100\, t_0}$ *; and with* $t_1 < \frac{1}{c_\ddagger} t_*$ *and such that* $K(\surd t_1) \ge 16 c_*(1 + K(\frac{1}{c_*} t_*) + K_*)$. *There exists* $\kappa > 1$ *which is independent of the given Nahm pole solution and* $[t_0, t_1]$ *except possible dependence on* $z_*$ *(as defined above), and which has the following significance: If* $t \in [t_\ddagger, t_1]$, *then*

$$K^2 - \frac{3}{4t^2} \le \kappa \frac{1}{\sqrt{t}} \frac{1}{\sqrt{t_1}} + \kappa \frac{1}{\sqrt{t}} \Big(1 + |\ln(t)| + \Big| \int_{\{t\}\times Y} c \Big| \Big)^{1/2} .$$

***Proof of Lemma 7.2***: The convention in this proof allows the values of the incarnations of $c_0$ to depend implicitly on given upper bounds for $z_*$. If this is $\mathcal{O}(1)$, then so is $c_0$. The incarnations of $c_0$ are otherwise independent of the given Nahm pole solution $[A, \mathfrak{a}]$ and the interval $[t_0, t_1]$.

The proof starts with a definition: Given $t \in [t_0, t_1]$, define $s_t \in [t, t_1]$ to be largest time s from $[t, t_1]$ such that

$$s \Big(\ln\big(\tfrac{s}{t}\big) + 1 + \Big| \int_{\{t\}\times Y} c \Big| \Big) \le 100\, c_* . \tag{7.1}$$

Note that there is a well defined $s_t$ for every $t \in [t_0, t_1]$. This is because the $\{t\}\times Y$ integral of $|c|$ is less than $\frac{2}{t} c_*$ (which is a consequence of Lemma 6.3). With regards to the size of $s_t$ relative to t: If $s \ln(\frac{s}{t}) \ge 1$, then s must be greater than $\frac{1}{|\ln t|}$ ; so s is much bigger than any fractional power of t. In any event, either

$$s_t = t_1 \quad \textit{or} \quad s_t \ge c_0^{-1}\Big(1 + |\ln t| + \Big| \int_{\{t\}\times Y} c \Big| \Big)^{-1} . \tag{7.2}$$

By way of a look ahead, the proof repeats the arguments used in the proof of Lemma 6.6 but with the updated bounds from Lemma 7.1 for various integrals that appear in Lemma 6.6's proof. In particular, (7.1) and the bounds from the versions of Lemmas 7.1 and 6.3 and 6.4 with $s \in [t, s_t]$ replacing t lead to the following:

- $\displaystyle\int_{\{t\}\times Y} \langle \mathfrak{a}, \mathcal{B} \rangle + \frac{1}{4} \int_{[t,s]\times Y} \big(|B_A|^2 + |E_A|^2 + |\nabla_A^{\perp} \mathfrak{a}|^2\big) \le c_0 \frac{1}{s}$ .

- $\int_{\{t\}\times Y} |t| \leq c_0 \frac{1}{\sqrt{ts}}$ .

(7.3)

The proof of Lemma 7.2 has three parts.

*Part 1*: The arguments in Part 1 of the proof of Lemma 6.6 that lead from (2.16) to (6.8) can be repeated but with the second bullet in (7.3) used to bound the $\{t\} \times Y$ integral of $|t|$. Doing this leads to the inequality

$$K^2 \leq \frac{3}{4^{1/3}} \Big( \int_{\{t\}\times Y} \mathcal{Q} + \Delta\Big)^{2/3} + c_0 \frac{1}{\sqrt{ts_t}} .$$

(7.4)

Here, $\Delta$ is defined by the right hand side of (6.9).

What with the top bullet in (7.3), the inequality in (6.10) can be replaced by this:

$$\int_{\{t\}\times Y} \mathcal{Q} \leq \frac{NK^2}{3t} + c_0 \frac{1}{s_t} .$$

(7.5)

Using this bound in (7.4) leads to a revised version of (6.11) which is this:

$$K^2 \leq \frac{3}{4^{1/3}} \Big(\frac{NK^2}{3t} + \Delta\Big)^{2/3} + c_0 \frac{1}{\sqrt{ts_t}} .$$

(7.6)

*Part 2*: This part of the proof provides suitable bounds for the values of $N$ and $\Delta$ that appear in (7.6). Start with $N$: Define $z_*$ by writing $N(t_0) = 1 + z_* t_0^2$ when $N(t_0) > 1$; and set $z_* = 0$ when $N(t_0) < 1$. Since it is assumed that $K(t_0) \geq \frac{1}{100\, t_0}$ , Lemma 6.3 can be brought to bear to see that $N(t) \leq 1 + z_* t^2$ for $t \in [t_0,t_1]$ if $z_* \geq c_0$.

Now consider $\Delta$: According to (6.16), the number $\Delta$ is at most $c_0 \frac{1}{t}$ which is small enough for the purposes at hand since it implies already that

$$\Big(\frac{NK^2}{3t} + \Delta\Big)^{2/3} \leq \frac{1}{3^{2/3} t^{2/3}} K^{4/3} (1 + c_0 t^2)^{2/3}$$

(7.7)

where $N \leq 1 + z_* t^2$ and $K^2(t) \geq \frac{3}{4t^2}$ . (It follows from (6.14) and (7.3) that $\Delta$ is in fact no greater than $c_0 \frac{1}{s_t}$ .)

*Part 3*: For the purposes of the proof of the lemma, it is sufficient to assume that $\kappa^2$ at a given time t is greater than $\frac{3}{4t^2}$. At such t, the bounds in (7.6) and (7.7) imply that

$$\kappa^2 \leq \left(\tfrac{3}{4}\right)^{1/3} \kappa^{4/3} \frac{1}{t^{2/3}} (1 + c_0 t^2) + c_0 \frac{1}{\sqrt{t s_t}} \tag{7.8}$$

This bound implies can hold with $\kappa^2 \geq \frac{3}{4t^2}$ only if

$$\kappa^2 - \frac{3}{4t^2} \leq c_0 \left(1 + \frac{1}{\sqrt{t s_t}}\right); \tag{7.9}$$

and the latter bound implies directly the one in Lemma 7.3.

**c) A priori bounds for the {t} × Y integrals of $c$, $|\mathfrak{c}|$ and $|\mathfrak{c}|^2$**

The assumptions for the first lemma in this section are the same as for the previous section.

**Lemma 7.3**: *Fix a Nahm pole solution* $(A, \mathfrak{a})$ *and let* $[t_0, t_1]$ *denote an interval in the corresponding version of* $\Omega^+$ *with* $t_0$ *such that both* $N(t_0) < 2$ *and* $\kappa(t_0) \geq \frac{1}{100\, t_0}$ *; and with* $t_1 < \frac{1}{c_\ddagger} t_*$ *and such that* $\kappa(\sqrt{t_1}) \geq 16 c_*(1 + \kappa(\frac{1}{c_*} t_*) + \kappa_*)$. *There exists* $\kappa > 1$ *which is independent of the given Nahm pole solution and* $[t_0, t_1]$ *except possible dependence on* $z_*$, *and it is such that* $\left| \int_{\{t\}\times Y} c \right| \leq \kappa + \left| \int_{\{t_0\}\times Y} c \right|$ *for all* $t \in [t_0, t_1]$.

***Proof of Lemma 7.3***: As in the previous proofs, the convention here allows the incarnations of $c_0$ to depend on the number $z_{**}$. The upcoming incarnations of $c_0$ are otherwise independent of the Nahm pole solution $(A, \mathfrak{a})$ and the interval $[t_0, t_1]$.

To start the proof, note that if the {t} × Y integral of $c$ is negative, then it is no smaller than $-c_0$. This follows from (6.31) and Lemma 6.6's bound $\kappa^2 - \frac{3}{4t^2} \leq c_0 \frac{1}{t}$.

To bound the {t} × Y integral of $c$ from above, assume that the maximum of the {t}×Y integrals of $c$ for $t \in [t_0, t_m]$ is positive and let $m$ denote this maximum. According to Lemma 7.2:

$$\kappa^2 - \frac{3}{4t^2} \leq c_0 \frac{1}{\sqrt{t}} \frac{1}{\sqrt{t_1}} + c_0 \frac{1}{\sqrt{t}} (1 + |\ln(t)| + m)^{1/2} \tag{7.10}$$

Letting $t_m$ denote the smallest $t \in [t_0, t_1]$ where the $\{t\}\times Y$ integral of $c$ has this maximal value, integrate the differential identity in (6.25) from $t_0$ to $t_m$ to see (using (6.30)) that

$$m \le \int_{\{t_0\}\times Y} c \;+ c_0\sqrt{t_m}\,\Big(\tfrac{1}{\sqrt{t_1}} + |\ln t_m|^{1/2} + c_0\,|m|^{1/2}\Big) \quad . \tag{7.11}$$

Since $t_m \le t_1$, this leads directly to Lemma 7.3's bound.

The next lemma uses Lemma 7.3 to update the bounds in Lemma 6.7.

**Lemma 7.4**: *Fix a Nahm pole solution* $(A, \mathfrak{a})$ *and let* $[t_0, t_1]$ *denote an interval in the corresponding version of* $\Omega^+$ *with* $t_0$ *such that both* $N(t_0) < 2$ *and* $K(t_0) \ge \frac{1}{100\,t_0}$ *; and with* $t_1 < \frac{1}{c_\ddagger} t_*$ *and such that* $K(\sqrt{t_1}) \ge 16c_*(1 + K(\frac{1}{c_*} t_*) + K_*)$. *There exists* $\kappa > 1$ *which is independent of the given Nahm pole solution and* $[t_0, t_1]$ *except possible dependence on* $z_*$, *and which has the following significance:*

- $\displaystyle\int_{[t_0,t_1]\times Y} (|\mathfrak{c}^+|^2 + |\mathfrak{c}^-|^2) \le \kappa\Big(1 + \Big|\int_{\{t_0\}\times Y} c\,\Big|\Big)$ .
- $\displaystyle\int_{\{t\}\times Y} |\mathfrak{c}|^2 \le \kappa\,\tfrac{1}{t}\Big(1 + \Big|\int_{\{t_0\}\times Y} c\,\Big|\Big)$ *for all* $t \in [t_0, t_1]$ .
- $\displaystyle\int_{\{t\}\times Y} |\mathfrak{c}| \le \kappa\Big(1 + \int_{\{t_0\}\times Y} (|c| + |\mathfrak{c}^-|)\Big)$ *for all* $t \in [t_0, t_1]$ .

***Proof of Lemma 7.4***: To prove the top bullet, integrate (6.25) from $t_0$ to $t_1$ using (6.27), and using (7.10) with Lemma 7.3's bound for $m$ to bound $\kappa^2 - \frac{3}{4t^2}$. The second bullet follows from (6.28) using Lemma 7.3's bound. The third bullet is proved by considering separately the $\{t\}\times Y$ integrals of $|c|$, $|\mathfrak{c}^-|$ and $|\mathfrak{c}^+|$. The asserted bounds for the first two (the integrals of $|c|$ and $|\mathfrak{c}^-|$) are obtained by integrating (6.29) and (6.31) using the bound in (6.27) and the bounds in the first and second bullets of this lemma. The asserted bound for the $\{t\}\times Y$ integral of $|\mathfrak{c}^+|$ is obtained by integrating both sides of (6.32) and then using Lemma 6.5 with the second bullet of this lemma.

### d) Bounds for the $[t_0, t_1]\times Y$ integrals of $|B_A|^2$, $|E_A|^2$, $|\nabla_A{}^{\perp}\mathfrak{a}|^2$ and $|\mathfrak{c}|^2$

This subsection considers integrals on $[t_0, t_1]\times Y$ when $t_1 \le \frac{1}{c_*} t_*$ and when it is either a boundary point of $\Omega^-$ with $K(\sqrt{t_1}) \ge 16c_*(1 + K(\frac{1}{c_*} t_*) + K_*)$, or else it is not a boundary point of $\Omega^-$ in which case $K(\sqrt{t_1}) = 16c_*(1 + K(\frac{1}{c_*} t_*) + K_*)$ (which is to say that $t_1$ is as large as possible for appeals to preceding lemmas). If $t_1$ is as just described, then one of the two bullets that follow hold:

- $\int_{\{t_1\}\times Y} \langle \mathfrak{a}, \mathcal{B}\rangle = 0$ *(which is to say that* $t_1$ *is a boundary point of* $\Omega^-$ ).
- $\int_{\{t_1\}\times Y} \langle \mathfrak{a}, \mathcal{B}\rangle \le c_{\ddagger}\frac{1}{t_1}$ .

$$(7.12)$$

As before, $c_{\ddagger}$ signifies a number that is larger than the versions of $\kappa$ that appear in Lemmas 6.2–6.5. (The second bullet follows from Lemma 6.4.) See the comment after Lemma 6.4.) With the conditions in (7.12) understood, and what with the top bullet of Lemma 7.4, an appeal to the top bullet in Lemma 5.5 leads to this:

$$\int_{\{t_0\}\times Y} \langle \mathfrak{a}, \mathcal{B}\rangle + \int_{[t_0,t_1]\times Y} (|B_A|^2 + |E_A|^2 + |\nabla_A^{\perp}\mathfrak{a}|^2) \le x + c_0(1 + \Big|\int_{\{t_0\}\times Y} c\,\Big|). \tag{7.13}$$

where $x$ is either 0 or the $\{t_1\}\times Y$ integral of $\langle \mathfrak{a}, \mathcal{B}\rangle$ depending on whether the top bullet of (7.13) describes $t_1$ or the lower bullet does.

The next lemma gives a stronger bound for the right hand side of (7.13) when the number $x$ in (7.13) is zero. This lemma introduces new notation, a number $E_0$:

$$E_0 = \int_{(0,t_0]\times Y} (|B_A|^2 + |E_A|^2 + |\nabla_A^{\perp}\mathfrak{a}|^2)\ , \tag{7.14}$$

This number $E_0$ in conjunction with the $\{t_0\}\times Y$ integral of $c$ can be viewed as 'initial conditions' for the interval $[t_0, t_1]$.

**Lemma 7.5**: *Fix a Nahm pole solution* $(A, \mathfrak{a})$ *and let* $[t_0, t_1]$ *denote an interval in the corresponding version of* $\Omega^+$ *with* $t_0$ *such that both* $N(t_0) < 2$ *and* $K(t_0) \ge \frac{1}{100\,t_0}$ *; and with* $t_1$ *being less than* $\frac{1}{c_*}t_*$ *and a boundary point of* $\Omega^-$ *such that* $K(\sqrt{t_1}) \ge 16c_*(1 + K(\frac{1}{c_*}t_*) + K_*)$. *There exists* $\kappa > 1$ *which is independent of the given Nahm pole solution and* $[t_0, t_1]$ *except possible dependence on z, and which has the following significance:*

$$\int_{\{t_0\}\times Y} \langle \mathfrak{a}, \mathcal{B}\rangle + \int_{[t_0,t_1]\times Y} (|B_A|^2 + |E_A|^2 + |\nabla_A^{\perp}\mathfrak{a}|^2 + |\mathfrak{c}|^2) \le$$

$$\kappa(t_1 - t_0)\Big(1 + \big(\int_{\{t_0\}\times Y} (|c| + |\mathfrak{c}^-|)\big)^2\Big) + \kappa(E_0 + x)\,.$$

***Proof of Lemma 7.5***: The convention in this proof again has $c_0$ depending possibly on the numbers $z_*$. And, as before, the incarnations of $c_0$ are otherwise independent of the given Nahm pole solution $[A, \mathfrak{a}]$ and the interval $[t_0, t_1]$.

The proof starts with the inequality in the top bullet of Lemma 5.5 for the case where $t = t_0$ and $s = t_1$:

$$\int_{\{t_0\}\times Y}\langle\mathfrak{a},\mathcal{B}\rangle+\tfrac{1}{4}\int_{[t_0,t_1]\times Y}(|B_A|^2+|E_A|^2+|\nabla_A^{\perp}\mathfrak{a}|^2)\le c_0(t_1-t_0)+c_0\int_{[t_0,t]\times Y}(|\mathfrak{c}^+|^2+|\mathfrak{c}^-|^2)+c_0 x\,. \tag{7.15}$$

Add $\frac{1}{4}E_0$ to the right and left hand sides of this inequality; then use Lemma 4.1 to see that

$$\int_{[t_0,t_1]\times Y}|\nabla_A^{\perp}\mathfrak{c}|^2\le c_0(t_1-t_0)+c_0\int_{[t_0,t]\times Y}(|\mathfrak{c}^+|^2+|\mathfrak{c}^-|^2)+c_0(E_0+x)\,. \tag{7.16}$$

Now invoke the $\mathfrak{q} = \mathfrak{c}$ and $\varepsilon = c_0^{-1}$ version of the inequality in Lemma 2.3. Doing this leads from (7.16) to this:

$$\int_{[t_0,t_1]\times Y}|\mathfrak{c}|^2\le c_0(t_1-t_0)+c_0\int_{[t_0,t_1]}\Big(\int_{\{s\}\times Y}|\mathfrak{c}|\Big)^2 ds\ +\ c_0(E_0+x)\,. \tag{7.17}$$

And, given what is said by the third bullet of Lemma 7.4, this in turn leads to:

$$\int_{[t_0,t_1]\times Y}|\mathfrak{c}|^2\le c_0(t_1-t_0)\Big(1+\Big(\int_{\{t_0\}\times Y}(|c|+|\mathfrak{c}^-|)\Big)^2\Big)+c_0(E_0+x). \tag{7.18}$$

Using the latter bound for the integral on the right hand side of (7.15) gives the bound that is asserted by the lemma.

## e) Refined bounds for the $\{t_1\}\times Y$ integrals of $|c|$ and $|\mathfrak{c}^-|$

The next lemma refines the bounds implied by the third bullet in Lemma 7.4 for the $\{t_1\}\times Y$ integral of $|c|$ and $|\mathfrak{c}^-|$. (It is silent about the corresponding integral of $|\mathfrak{c}^+|$.)

**Lemma 7.6**: *Fix a Nahm pole solution* $(A,\mathfrak{a})$ *and let* $[t_0, t_1]$ *denote an interval in the corresponding version of* $\Omega^+$ *with* $t_0$ *such that both* $N(t_0) < 2$ *and* $K(t_0) \ge \frac{1}{100\,t_0}$ *; and with* $t_1 \le \frac{1}{c_*}t_*$ *such that* $K(\sqrt{t_1}) \ge 16c_*(1 + K(\frac{1}{c_*}t_*) + K_*)$. *There exists* $\kappa > 1$ *which is independent of the* $(A,\mathfrak{a})$ *and* $[t_0, t_1]$ *except for possible dependence on z with the following significance: If the* $\{t_1\}\times Y$ *integral of* $\langle\mathfrak{a},\mathcal{B}\rangle$ *is zero, then*

$$\int_{\{t_1\}\times Y} (|c| + |c^-|) \le \kappa(t_1 - t_0)\,(1+(\int_{\{t_0\}\times Y} |c|\,)^2) + \kappa E_0 + \tfrac{t_0}{t_1} \int_{\{t_0\}\times Y} (|c| + |c^-|)\,.$$

***Proof of Lemma 7.6***: As was the case in the last proof, the convention in what follows allows $c_0$ to depend implicitly on given upper bounds for the number $z_*$.

To prove the asserted bound for the $\{t_1\}\times Y$ integral of $|c|$, start with the $t = t_0$ and $s = t_1$ version of (6.30). Then use Lemma 7.5's bound for the integrals of $|B_A|^2$ and $|c|^2$ that appear on (6.30)'s right hand side.

To prove the asserted bound for the $\{t_1\}\times Y$ integral of $|c^-|$, integrate both sides of (6.31) on $[t_0, t_1]$ and then use Lemma 7.4's bounds as in the preceding paragraph.

## 8. Where the $\{t\} \times Y$ integral of $\langle \mathfrak{a}, \mathcal{B}\rangle$ is negative

Remember that $\Omega^-$ denotes the set of $t \in [t^-, \frac{1}{8} t_*]$ where the integral of $\langle \mathfrak{a}, \mathcal{B}\rangle$ on $\{t\} \times Y$ is negative The plan for this section is to derive a priori bounds for the $\{t\}\times Y$ integrals of $|c|, |c^-|, |c^+|$ and $|c|^2$ for $t \in \Omega^-$. The section also derives bounds for the $I \times Y$ integrals of $|B_A|^2$, $|E_A|^2$ and $|\nabla_A{}^{\perp}\mathfrak{a}|^2$ and $|c|^2$ when I is an interval in $\Omega^-$. The final bounds in this section are analogous to those in Proposition 6.1. The upcoming Proposition 8.1 makes the precise statements. As in the previous sections, $c_*$ denotes the version of $\kappa$ that appears in Lemma 5.8.

**Proposition 8.1**: *There exists* $\kappa > 32$ *with the following significance: Having fixed a Nahm pole solution* $(A, \mathfrak{a})$*, let* $(t_0, t_1)$ *denote an open set in the corresponding version of* $\Omega^-$ *with* $t_0$ *being the lower boundary point of a component of* $\Omega^-$*. Assume that* $N(t_0) < 2$*; and assume that* $t_1 \le \frac{1}{c_*} t_*$ *and* $K(\surd t_1) \ge 16 c_*(1+K(\frac{1}{c_*} t_*)+K_*)$*. Define* $E_0$ *using the rule in (7.14); and define* $m_0$ *to be the* $\{t_0\}\times Y$ *integral of* $|c|+|c^-|$*. Assume that both* $E_0$ *and* $m_0$ *are less than* $\kappa^{-2}$*. Fix* $t \in [t_0, t_1]$*. Under these circumstances,*

- $\int_{\{t\}\times Y} |c^+| \le \kappa.$
- $\int_{\{t\}\times Y} (|c| + |c^-|) \le \kappa((t - t_0) + E_0 + m_0).$
- $|\int_{\{t\}\times Y} \langle \mathfrak{a}, \mathcal{B}\rangle| + \int_{[t_0, t]\times Y} (|B_A|^2 + |E_A|^2 + |\nabla_A^{\perp}\mathfrak{a}|^2 + |c|^2) \le \kappa((t - t_0) + E_0 + m_0)\,,$
- *If* $z_* \ge \kappa$ *and* $N(t_0) \le 1 + z_* t_0^2$ *and* $K(t_0) \ge \frac{1}{100\, t_0}$*, then* $N(t) \le 1 + z_* t^2$ *for all* $t \in [t_0, t_1]$ .

The proof of this proposition is in Section 8e. Looking ahead, it differs from the proof of the analogous assertions in Proposition 6.1 for $t \in \Omega^+$ at the very start with regards to the $\Omega^-$ analogs of the lemmas in Section 6. The differences can be traced back to Lemma 5.7 and its third bullet in particular: That bullet is not in play on $\Omega^-$. However, the fourth bullet of Lemma 5.7 is in play *if* there exists a not so very negative lower bound for the $\{t\} \times Y$ integral of $\langle\mathfrak{a},\mathcal{B}\rangle$ for $t \in \Omega^-$. Thus, the task (stated simplistically) is this: Derive such a bound. This is essentially what is done below (although the arguments are framed differently).

Keep in mind for what is to come that the interval $(t_0, t_1)$ is such that $t_0$ is a lower boundary point of a component of $\Omega^-$. In particular, this means that if $t_0$ is not $t^-$, then the $\{t_0\}\times Y$ integral of $\langle\mathfrak{a},\mathcal{B}\rangle$ is zero. Also, if $t_0 = t^-$, then the $\{t_0\}\times Y$ integral of $\langle\mathfrak{a},\mathcal{B}\rangle$ is no smaller than $-c_0 t_0 + \mathrm{E}_0$ (because of Lemma 3.6). Also keep in mind that $t_1$ is no greater than $\frac{1}{c_*} t_*$ *and that* $\mathrm{K}(\sqrt{t_1}) \geq 16c_*(1+\mathrm{K}(\frac{1}{c_*}t_*)+\mathrm{K}_*)$. These constraints on $(t_0, t_1)$ are for the most part implicit in what follows.

## a) Initial $\Omega^-$ bounds

Fix $t \in (t_0, t_1)$ and then invoke the top bullet of Lemma 5.5 using $t_0$ for what that bullet calls t and using t for what that bullet calls s. The result is the following inequality:

$$\Big|\int_{\{t\}\times Y}\langle\mathfrak{a},\mathcal{B}\rangle\Big| + \tfrac{1}{4}\int_{[t_0,t]\times Y}(|B_A|^2+|E_A|^2+|\nabla_A^{\perp}\mathfrak{a}|^2) \leq c_0\Big((t-t_0)+\mathrm{E}_0+c_0\int_{[t_0,t]\times Y}(|\mathfrak{c}^+|^2+|\mathfrak{c}^-|^2)\Big) \tag{8.1}$$

(The $\mathrm{E}_0$ term on the right hand side can be dropped in the event that $t_0 > t^-$.) Note in particular that the absolute value of the $\{t\} \times Y$ integral of $\langle\mathfrak{a},\mathcal{B}\rangle$ is on the left hand side. This is because the $\{t\}\times Y$ integral of $\langle\mathfrak{a},\mathcal{B}\rangle$ is negative. With (8.1) understood, the argument used in the proof of Lemma 7.5 from (7.15) through (7.17) can be repeated with only notational changes to derive the next inequality from (8.1):

$$\Big|\int_{\{t\}\times Y}\langle\mathfrak{a},\mathcal{B}\rangle\Big| + \int_{[t_0,t]\times Y}(|B_A|^2+|E_A|^2+|\nabla_A^{\perp}\mathfrak{a}|^2+|\mathfrak{c}|^2) \leq c_0\Big((t-t_0)+\mathrm{E}_0+\int_{[t_0,t]}\Big(\int_{\{s\}\times Y}|\mathfrak{c}|\Big)^2 ds\Big). \tag{8.2}$$

Upper bounds for the integrals on the left hand side of the preceding inequality therefore follow from upper bounds for $\mathrm{E}_0$ and for the $\{s\}\times Y$ integrals of $|\mathfrak{c}|$ when $s \in [t_0, t]$.

The equation in (6.30) and the analogous equation that is obtained by integrating both sides of (6.31) lead (with the help of (8.2)) to the following:

$$\int_{\{t\}\times Y}(|c|+|\mathfrak{c}^-|) \leq c_0\Big((t-t_0)+\mathrm{E}_0+\int_{[t_0,t]}\Big(\int_{\{s\}\times Y}|\mathfrak{c}|\Big)^2 ds\Big) + \frac{t_0}{t}\int_{\{t_0\}\times Y}(|c|+|\mathfrak{c}^-|)\,. \tag{8.3}$$

The middle equation in (4.9) does not lead (at least directly) to a $|\mathfrak{c}^+|$ analog of (8.3) because of the appearance of a *negative* multiple of $\frac{1}{t}\mathfrak{c}^+$ on its left hand side.

However, the middle equation of (4.9) can be used to derive a bound for the change in the $\{s\} \times Y$ integrals of $|\mathfrak{c}^+|$ on any interval of the form $[t, 2t]$ for $t \in [t_0, \frac{1}{2} t_1]$. Here is the precise statement to that effect: If $t \in [t_0, \frac{1}{2} t_1]$ and $s \in [t, 2t]$, then

- $\int_{\{s\}\times Y} |\mathfrak{c}^+| \le 2 \int_{\{t\}\times Y} |\mathfrak{c}^+| + c_0\big((t - t_0) + E_0 + \int_{[t_0,t]} \big( \int_{\{s\}\times Y} |\mathfrak{c}|\big)^2 ds\big)$ .
- $\int_{\{s\}\times Y} |\mathfrak{c}^+| \ge \int_{\{t\}\times Y} |\mathfrak{c}^+| - c_0\big((t - t_0) + E_0 + \int_{[t_0,t]} \big( \int_{\{s\}\times Y} |\mathfrak{c}|\big)^2 ds\big)$ .

(8.4)

To prove this, take the inner product of both sides of the middle bullet of (4.9) with $\mathfrak{c}^+|\mathfrak{c}^+|^{-1}$ where $\mathfrak{c}^+ \ne 0$ to obtain inequalities that can be written as

- $\frac{\partial}{\partial t} \big( \frac{1}{t} |\mathfrak{c}^+|\big) \le c_0 \frac{1}{t} |B_A| + c_0 \frac{1}{t} |\mathfrak{c}|^2$ .
- $\frac{\partial}{\partial t} \big( \frac{1}{t} |\mathfrak{c}^+|\big) \ge - c_0 \frac{1}{t} |B_A| - c_0 \frac{1}{t} |\mathfrak{c}|^2$ .

(8.5)

Integrate both sides of these over $\{t\} \times Y$ and then integrate from the given value of t to a given $s \in [t, 2t]$. Doing this gives the following inequality

$$\Big| \frac{1}{s} \int_{\{s\}\times Y} |\mathfrak{c}^+| - \frac{1}{t} \int_{\{t\}\times Y} |\mathfrak{c}^+| \Big| \le c_0 \frac{1}{t} \Big(\sqrt{t}\, \Big( \int_{[t,2t]} |B_A|^2 \Big)^{1/2} + \int_{[t,2t]} |\mathfrak{c}|^2 \Big) ,$$

(8.6)

at which gives (8.4) when (8.2) is invoked.

**b) Initial assumptions for the $\{t\} \times Y$ integral of $|\mathfrak{c}|$**

Define the functions $m$ and $r$ on $[t_0, t_1]$ using the following rules:

$$m(t) = \sup_{s\in[t_0,t]} \int_{\{s\}\times Y} (|c| + |\mathfrak{c}^-|) \quad \textit{and} \quad r(t) = \sup_{s\in[t_0,t]} \int_{\{t\}\times Y} |\mathfrak{c}^+| \ .$$

(8.7)

Using this definition with (8.3) leads to this bound:

$$\int_{\{t\}\times Y} (|c| + |\mathfrak{c}^-|) \le c_0 (t - t_0)\, (1 + m^2 + r^2) + c_0 E_0 + \frac{t_0}{t}\, m(t_0).$$

(8.8)

Note in particular that if $t \in (t_0, t_1]$ is such that

$$m(t) < c_0^{-1} \frac{1}{t - t_0}$$

(8.9)

then the bound that is asserted by (8.8) leads to this bound:

$$m(\mathrm{t}) \leq c_0(\mathrm{t} - \mathrm{t}_0)\,(1 + r^2) + c_0 \mathrm{E}_0 + m(\mathrm{t}_0). \tag{8.10}$$

To see this, apply (8.8) at the point where the sup on the left in (8.7) is taken on. When doing this, keep in mind that $m$ and $r$ are non-decreasing. The essential point is this: If $\mathcal{X} \geq 1$ and if $m(\mathrm{t}) \leq \frac{1}{2\mathcal{X}} \frac{1}{\mathrm{t}-\mathrm{t}_0}$, then $\mathcal{X}(\mathrm{t}-\mathrm{t}_0)\, m(\mathrm{t})^2 \leq \frac{1}{2} m(\mathrm{t})$.

The bound in (8.10) says in effect that if there is no $r$ (which is to say no $\mathfrak{c}^+$), then an upper bound for $m$ is determined by $m(\mathrm{t}_0)$ and $\mathrm{E}_0$ if these aren't too big. This is because (8.10) without $r$ implies that (8.9) holds at $\mathrm{t} = \mathrm{t}_1$ if $m(\mathrm{t}_0)$ and $\mathrm{E}_0$ are at most $\mathcal{O}(1)$.

Under the circumstances $\mathfrak{c}^+$ is almost surely present. Therefore (8.10) bounds $m(\mathrm{t})$ given a bound for $r(\mathrm{t})$ and the initial conditions at $\mathrm{t}_0$, assuming that t is where (8.9) holds. But note in this regard that (8.9) holds by virtue of (8.10) if $\mathrm{E}_0$ and $m(\mathrm{t}_0)$ are not too large, and if the function $r(\mathrm{t})$ obeys $r(\mathrm{t}) \leq c_0^{-1} \frac{1}{\mathrm{t}-\mathrm{t}_0}$.

The preceding motivates an observation and then a definition and then a lemma. Here are the observation and definition: The function $\mathrm{t} \to (\mathrm{t}-\mathrm{t}_0)\, r(\mathrm{t})$ is an increasing function on $[\mathrm{t}_0, \mathrm{t}_1]$ because $r$ is non-decreasing and $(\mathrm{t}-\mathrm{t}_0)$ is increasing. It is also zero at $\mathrm{t}_0$. Therefore, given any $\mathrm{Z} > 1$, either $(\mathrm{t}_1 - \mathrm{t}_0)\, r(\mathrm{t}_1) < \frac{1}{\mathrm{Z}}$ or there is a unique t in $(\mathrm{t}_0, \mathrm{t}_1)$ where $(\mathrm{t}-\mathrm{t}_0)\, r(\mathrm{t})$ equala $\frac{1}{\mathrm{Z}}$. Let $\mathrm{t}_\mathrm{Z}$ denote either $\mathrm{t}_1$ or this unique t (as the case may be).

**Lemma 8.2**: *There exists* $\kappa > 1$ *with the following significance: Let* $(A, \mathfrak{a})$ *denote a Nahm pole solution and let* $(\mathrm{t}_0, \mathrm{t}_1)$ *be an open set in* $\Omega^-$ *as described at the beginning of Proposition 8.1. Assume that* $m(\mathrm{t}_0)$ *and* $\mathrm{E}_0$ *are less than* $\frac{1}{\kappa}$*; and that* $\mathrm{Z} > \kappa$. *If* $\mathrm{t} \in [\mathrm{t}_0, \mathrm{t}_\mathrm{Z}]$, *then* $m$ *obeys the bounds in (8.9) and (8.10).*

***Proof of Lemma 8.2*:** If $r(\mathrm{t}_\mathrm{Z}) \leq \frac{1}{\mathrm{Z}(\mathrm{t}_\mathrm{Z} - \mathrm{t}_0)}$, then $r(\mathrm{t}) \leq \frac{1}{\mathrm{Z}(\mathrm{t}_\mathrm{Z} - \mathrm{t}_0)}$ for $\mathrm{t} < \mathrm{t}_\mathrm{Z}$ because $r$ is non-decreasing. As a consequence, (8.8) leads to the following:

$$m(\mathrm{t}) \leq c_0(\mathrm{t}-\mathrm{t}_0)\, m^2 + c_0 \frac{1}{\mathrm{Z}(\mathrm{t}_\mathrm{Z} - \mathrm{t}_0)} + c_0\mathrm{E}_0 + m(\mathrm{t}_0). \tag{8.11}$$

Let $\mathrm{t}_{\mathrm{Z}*}$ denote the smallest $\mathrm{t} \in [\mathrm{t}_0, \mathrm{t}_\mathrm{Z}]$ where $m(\mathrm{t}) = \frac{1}{\mathrm{Z}^{1/2}(\mathrm{t} - \mathrm{t}_0)}$ if such t exists. Supposing it does exist, then it follows from (8.11) that

$$\frac{1}{\mathrm{Z}^{1/2}(\mathrm{t}_{\mathrm{Z}*} - \mathrm{t}_0)} \leq c_0 \frac{1}{\mathrm{Z}(\mathrm{t}_{\mathrm{Z}*} - \mathrm{t}_0)} + c_0\mathrm{E}_0 + m(\mathrm{t}_0)\,. \tag{8.12}$$

The preceding inequality is nonsense if $Z > c_0$ and if $E_0$ and $m(t_0)$ are less than $c_0^{-1}$. Therefore, $m(t) < \frac{1}{Z^{1/2}(t - t_0)}$ on $[t_0, t_Z]$ if $Z > c_0$ and if $E_0$ and $m(t_0)$ are less than $c_0^{-1}$.

**c) Bounding $r$ using the $\{t\} \times Y$ integrals of $|\mathfrak{t}|$**

With the observations from the preceding subsection in mind, the question now is this: Given Z, is $t_Z = t_1$ or not? Or, put differently: How large can $r(t)$ be? The following lemma takes the first step towards an answer.

**Lemma 8.3**: *There exists* $\kappa > 1$ *with the following significance: Fix a Nahm pole solution* $(A, \mathfrak{a})$ *and let* $(t_0, t_1)$ *denote an open set in* $\Omega^-$ *as described at the beginning of Proposition 8.1. Suppose that* $Z > \kappa$ *and* $m(t_0)$ *and* $E_0$ *are less than* $\frac{1}{\kappa}$. *If* $t \in [t_0, t_Z]$, *then*

$$r(t) \le \kappa \sup_{s \in [t_0, \frac{1}{2}t]} \Big( \int_{[s,2s]\times Y} |\mathfrak{t}| \Big) + \kappa\big((t - t_0) + E_0 + m(t_0)\big) .$$

***Proof of Lemma 8.3***: Start with (6.32) which, after integrating on $\{t\} \times Y$, implies this:

$$\frac{1}{t} \int_{\{t\}\times Y} |\mathfrak{c}^+| \le \int_{\{t\}\times Y} |\mathfrak{t}| + c_0 \int_{\{t\}\times Y} |\mathfrak{c}|^2 . \tag{8.13}$$

Now fix $t' \in [t_0, \frac{1}{2} t_1]$ and integrate (8.13) from $t'$ to $2t'$:

$$\frac{1}{t'} \int_{[t',2t']\times Y} |\mathfrak{c}^+| \le c_0 \int_{[t',2t']\times Y} |\mathfrak{t}| + c_0 \int_{[t',2t']\times Y} |\mathfrak{c}|^2 . \tag{8.14}$$

Meanwhile, in light of Lemma 8.2 and (8.10), the inequality in (8.2) implies the following when $t \in [t_0, t_Z]$:

$$\Big| \int_{\{t\}\times Y} \langle \mathfrak{a}, \mathcal{B} \rangle \Big| + \int_{[t_0,t]\times Y} \big(|B_A|^2 + |E_A|^2 + |\nabla_A^{\perp}\mathfrak{a}|^2 + |\mathfrak{c}|^2\big) \le c_0 (t - t_0)(1 + r(t)^2) + c_0 E_0 + c_0 m(t_0). \tag{8.15}$$

The right hand side of the $t = t'$ version of (8.15)'s inequality can be used to replace the $|\mathfrak{c}|^2$ integral on the right hand side of (8.14)'s inequality to obtain:

$$\frac{1}{t'} \int_{[t',2t']\times Y} |\mathfrak{c}^+| \le c_0 \int_{[t',2t']\times Y} |\mathfrak{t}| + c_0 (t - t_0)(1 + r(2t')^2) + c_0 E_0 + c_0 m(t_0). \tag{8.16}$$

Furthermore, in light of Lemma 8.2 and (8.10), this inequality and those in (8.4) with t replaced by $t'$ lead to the following: If $s \in [t', 2t']$, then

$$\int_{\{s\}\times Y} |\mathfrak{c}^+| \le c_0 \int_{[t', 2t']\times Y} |\mathfrak{t}| + c_0 (t - t_0)(1 + r(2t')^2) + c_0 E_0 + c_0 m(t_0).$$

(8.17)

Now take $t' \le \frac{1}{2}t$ and then s so that s is a point in $[t_0, t]$ where the $\{s\} \times Y$ integral of $|\mathfrak{c}^+|$ is $r(t)$. The inequality in (8.17) for such a choice for s leads to the bound

$$r(t) \le c_0 \int_{[t', 2t']\times Y} |\mathfrak{t}| + c_0 (t - t_0)(1 + r(t)^2) + c_0 E_0 + c_0 m(t_0).$$

(8.18)

Keeping in mind that if $r \le \frac{1}{Z(t_Z - t_0)}$ then $(t - t_0)\, r^2 < \frac{1}{Z}\, r$, this in turn implies the bound

$$(1 - c_0 \tfrac{1}{Z})\, r \le c_0 \int_{[t', 2t']\times Y} |\mathfrak{t}| + c_0 (t - t_0) + c_0 E_0 + c_0 m(t_0) .$$

(8.19)

The latter leads directly to the lemma's bound if $Z > c_0$.

**d) Bounding the $\{t\} \times Y$ integrals of $|\mathfrak{t}|$**

An appeal to Lemma 8.3 requires a reasonable bound on the integral of $|\mathfrak{t}|$ that appears on the right hand side of that lemma's inequality. The next lemma supplies this.

**Lemma 8.4**: *There exists* $\kappa > 32$ *with the following significance: Fix a Nahm pole solution* $(A, \mathfrak{a})$ *and an open set* $(t_0, t_1)$ *in* $\Omega^-$ *as described at the beginning of Proposition 8.1. Assume that* $m(t_0)$ *and* $E_0$ *are less than* $\frac{1}{\kappa}$. *If* $Z \ge \kappa$ *and if* $t \in [t_0, t_Z]$, *then*

$$\int_{\{t\}\times Y} |\mathfrak{t}| \le \kappa \tfrac{1}{t} .$$

*Therefore, if* $s \in [t_0, \frac{1}{2} t_Z]$, *then* $\int_{[s, 2s]\times Y} |\mathfrak{t}| \le \kappa \ln 2.$

The proof of this lemma requires the following auxilliary lemma.

**Lemma 8.5**: *There exists* $\kappa > c_*$ *with the following significance: Fix a Nahm pole solution* $(A, \mathfrak{a})$ *and an open set* $(t_0, t_1)$ *of* $\Omega^-$ *as described at the beginning of Proposition 8.1. Assume that* $m(t_0)$ *and* $E_0$ *are less than* $\frac{1}{\kappa}$. *If* $Z \ge \kappa$, *then* $K(s) \le \frac{\kappa}{s}$ *and* $N(s) < 2$ *for all* $s \in [t_0, \frac{1}{\kappa}\sqrt{t}]$.

Lemma 8.5 is proved momentarily.

***Proof of Lemma 8.4***: The proof assumes Lemma 8.5. Let $c_\Diamond$ denote the version of $\kappa$ that appears in Lemma 8.5. Fix $s \in [t, \frac{1}{c_\Diamond}\sqrt{t}\,]$ for the moment and invoke the second bullet of Lemma 5.7 to obtain the following inequality

$$\int_{\{t\}\times Y} |\mathfrak{t}| \le \int_{\{s\}\times Y} |\mathfrak{t}| \; + \; c_0\Big( \int_{[t,s]} \kappa^2 \Big)^{1/2} \Big( \int_{\{s\}\times Y} \langle \mathfrak{a}, \mathcal{B}\rangle - \int_{\{t\}\times Y} \langle \mathfrak{a}, \mathcal{B}\rangle + \kappa\Big(s + \int_{[t,s]} \kappa^2 \Big)^{1/2} \; . \tag{8.20}$$

Concerning the terms on the right hand side of this inequality: The $\{s\}\times Y$ integral of $|\mathfrak{t}|$ is at most $c_0\kappa^2(s)$ which is at most $c_0 \frac{1}{s^2}$ by virtue of Lemma 8.5. Likewise by virtue of Lemma 8.5, the integrals of $\kappa^2$ that appear are at most $c_0 \frac{1}{t}$. Meanwhile, if the $\{s\}\times Y$ integral of $\langle \mathfrak{a}, \mathcal{B}\rangle$ is negative, it can be ignored in (8.20); and if it is positive, then it is at most $c_0\kappa^2(s)$ which follows from the second bullet of Lemma 5.6 when t there is replaced by s and when s there is replaced by $\frac{1}{4} t_*$. Thus, by virtue of Lemma 8.5, this is at most $c_0 \frac{1}{s^2}$. Finally, the $\{t\} \times Y$ integral of $\langle \mathfrak{a}, \mathcal{B}\rangle$ is no smaller than $-c_0(\frac{1}{t} + 1)$ because of (8.15) (assuming $E_0$ and $m(t_0)$ are less than $c_0$) and Lemma 8.5 (which implies that $r(t) \le c_0 \frac{1}{t}$ in any event). Using all of these bounds in (8.20) leads to this inequality:

$$\int_{\{t\}\times Y} |\mathfrak{t}| \le c_0 \frac{1}{s^2} + c_0 \frac{1}{t} \; . \tag{8.21}$$

In particular, if $s = \frac{1}{2c_\Diamond} \sqrt{t}$, then the right hand side of (8.21) is at most $c_0 \frac{1}{t}$.

***Proof of Lemma 8.5***: The lemma follows from Lemma 5.9 if $N(t) \le 2$, and if the $\{t\}\times Y$ integral of $\langle \mathfrak{a}, \mathcal{B}\rangle$ is greater than $-\kappa^2(t)$. (The lower bound on $\kappa(t)$ in Item a) of Lemma 5.9 is satisfied here because of the asssumptions in Proposition 8.1 with regards to $t_1$ and, given these, by virtue of (5.17).) Two cases will be distinguished. The first is when $t \le 2t_0$ and the second is when $t \ge 2t_0$. The two parts of the proof consider these cases.

*Part 1*: Assume here that $t \le 2t_0$. Since $\kappa(\sqrt{t_1}) \ge 16 c_*(1 + \kappa(\frac{1}{c_*} t_*) + \kappa_*)$, what is said by (5.17) implies that the value of $\kappa$ at $\sqrt{t_0}$ is greater than $4c_*(1 + \kappa(\frac{1}{16c_*} t_*) + \kappa_*)$. Since $N(t_0)$ is less than 2 by assumption, the conditions for Lemma 5.9 are met for $t = t_0$. Indeed, the only other condition is Requirement c), which is automatic in the case that $t_0 > t^-$, and which follows from Lemma 3.6 and (3.1) in the case when $t_0 = t^-$. Lemma 5.9 with t replaced by $t_0$ (and with $\Theta = 0$ in the case $t_0 > t^-$ and with $\Theta = c_0 t^-$ in the case $t_0 = t^-$) leads directly to the assertion of Lemma 8.5 for $t \in [t_0, 2t_0]$.

*Part 2*: Now suppose that $t \ge 2t_0$. This requires $t_Z \ge 2t_0$ because the assumption is that $t \le t_Z$. That implies in turn that $r(t) \le \frac{2}{Zt_Z}$ for any $t \in [t_0, t_Z]$. Use this bound in (8.15) to conclude the following: If $t \in [t_0, t_Z]$, then

$$| \int_{\{t\}\times Y} \langle \mathfrak{a}, \mathcal{B} \rangle | + \int_{[t_0,t]\times Y} (|B_A|^2 + |E_A|^2 + |\nabla_A^{\perp} \mathfrak{a}|^2 + |\mathfrak{c}|^2) \le c_0 (1 + \mathrm{E}_0 + (t - t_0)\, m(t_0)^2 + \tfrac{1}{Z^2 t_Z}) \,. \tag{8.22}$$

Meanwhile, using the same $r(t) \le \frac{2}{Zt_Z}$ with (8.10) leads to the lower bound

$$\mathrm{K}(t) \ge \tfrac{\sqrt{3}}{2t} - c_0 (\mathrm{E}_0 + m(t_0) + \tfrac{1}{Z^2 t_Z}) \tag{8.23}$$

Thus, if $Z > c_0$ and if $\mathrm{E}_0$ and $m(t_0) \le c_0$ and if $t \le c_0^{-1}$, then $\mathrm{K}^2(t)$ will be much larger than the right hand side of (8.22). As a consequence, Lemma 5.9 can be invoked if $\mathrm{N}(t) \le 2$ to obtain the assertion of Lemma 8.5. The next paragraph explains why $\mathrm{N}(t) \le 2$.

To see about $\mathrm{N}(t)$, note that (8.22) and (8.23) (assuming $Z > c_0$ and that $\mathrm{E}$ and $m_0$ are less than $c_0$; and that $t \ge c_0^{-1}$) can be used in (5.9) to obtain an inequality for the derivative of $\mathrm{N}$ on $[t_0, t_Z]$ that has the form

$$\tfrac{d}{dt} \mathrm{N} \le \tfrac{\mathrm{N}(1 - \mathrm{N})}{t} + c_0 t + c_0 t \mathrm{N} \;\; . \tag{8.24}$$

This inequality implies that $\mathrm{N}$ is decreasing on $[t_0, t_Z]$ where $\mathrm{N}$ is greater than $1 + c_0 t^2$.

**e) Proof of Proposition 8.1**

The proof of the proposition is a matter of walking back through the preceding sections. To elaborate, suppose that $t \in (t_0, t_1)$. Then Lemma 8.4 can be invoked. Use the bound from Lemma 8.4 in Lemma 8.3 to bound $r(t)$ by $c_0$. (This is the first bullet of Proposition 8.1.) Now use the (8.10) bound to bound $m(t)$ by $c_0((t - t_0) + \mathrm{E}_0 + m(t_0))$ (which is the second bullet of Proposition 8.1). Use the $r(t) \le c_0$ bound in (8.15) to obtain the third bullet of Proposition 8.1.

## 9. A priori bounds on the whole of $(0, 1] \times Y$

With a Nahm pole solution $(A, \mathfrak{a})$ in hand, then Propositions 6.1 and 8.1 give bounds for various $[t_0, t_1] \times Y$ integrals and $\{t\} \times Y$ integrals (for $t \in [t_0, t_1]$) which required a priori assumptions about behavior at $t_0$ and/or $t_1$. This section refines those bounds; and (this is crucial), it replaces the a priori assumptions in those propositions with just a single constraint at a comparitively large value of t. In particular, the end result of the analysis in this section is (roughly speaking) this: The pair $(A, \mathfrak{a})$ is described by the results in Section 3 at times less than a positive time defined by $(A, \mathfrak{a})$ that depends only on an upper bound for the number $(1 + \kappa(\frac{1}{c_*} t_*) + \kappa_*)$. The upcoming Proposition 9.1 makes a precise statement to this effect. (As before, $c_*$ is the version of the number $\kappa$ from Lemma 5.8.)

**Proposition 9.1**: *Let* $(A, \mathfrak{a})$ *denote a given Nahm pole solution. Fix* $\varepsilon \in (0, 1]$ *and there exists* $t_\varepsilon \in (0, \frac{1}{8} t_*]$ *which depends on an upper bound for* $1 + \kappa(\frac{1}{c_*} t_*) + \kappa_*$ *but is otherwise independent of* $(A, \mathfrak{a})$*; and it is such that*

- $\int_{(0, t_\varepsilon] \times Y} (|B_A|^2 + |E_A|^2 + |\nabla_A^\perp \mathfrak{a}|^2) < \varepsilon$ .
- $|\langle \mathfrak{a} \otimes \mathfrak{a} \rangle - \frac{1}{4t^2} \mathfrak{g}| < \frac{\varepsilon}{t^2}$ *when* $t \in (0, t_\varepsilon] \times Y$.

Note in particular that the two bullets of this proposition makes quantitative assertions with regards to the Nahm pole requirements in Definition 1.1. A parenthetical remark: The proof shows that the top bullet holds with $|\ln t_\varepsilon|$ on the order of $\frac{1}{\varepsilon}$; and that the lower bullet holds with $|\ln t_\varepsilon|$ on the order of a high power of $\frac{1}{\varepsilon}$.

The top bullet of Proposition 9.1 is proved in Section 9c assuming the second bullet. The second bullet is proved in Section 9d. Sections 9a and 9b each supply a preliminary lemma that is used in the proposition's proof.

### a) Initial bounds for integrals

The upcoming lemma is a weaker, preliminary version of Proposition 9.1. This lemma introduce a critical time for any given Nahm pole solution, denoted by $T_*$ which is defined from the solution by the rule

$$T_* = \frac{1}{c_*^2 (1 + \kappa(\frac{1}{c_*} t_*) + \kappa_*)^2} \quad . \tag{9.1}$$

This critical time $T_*$ depends on $(A, \mathfrak{a})$ to the extent that $\kappa_*$ and $\kappa(\frac{1}{c_*} t_*)$ depend on $(A, \mathfrak{a})$.

**Lemma 9.2**: *There exist* $\kappa > 1$*, and given a Nahm pole solution* $(A, \mathfrak{a})$*, there exists* $\lambda > 1$ *which is less than* $\kappa$ *if* $\int_{\{t\} \times Y} \langle \mathfrak{a}, \mathcal{B} \rangle \le 0$ *at some point in* $[\frac{1}{\kappa} T_*, \frac{1}{\kappa} \sqrt{T_*}]$ *and which depends only on an upper bound for* $1 + \kappa(\frac{1}{c_*} t_*) + \kappa_*$ *otherwise; and these are such that*

- $\int_{(0, \frac{1}{\kappa} T_*] \times Y} (|B_A|^2 + |E_A|^2 + |\nabla_A^\perp \mathfrak{a}|^2 + |\mathfrak{c}|^2) \le \lambda$ ,

- $\displaystyle\int_{\{\frac{1}{\kappa}t\}\times Y} |\mathfrak{c}| \le \lambda$ *for all* $t \in (0, \frac{1}{\kappa}\mathrm{T}_*]$ ,
- $\displaystyle\int_{\{\frac{1}{\kappa}t\}\times Y} |\mathfrak{c}|^2 \le \frac{\lambda}{t}$ *for all* $t \in (0, \frac{1}{\kappa}\mathrm{T}_*]$ ,
- $\mathrm{N}(t) \le 1 + \kappa t^2$ *for all* $t \in (0, \frac{1}{\kappa}\mathrm{T}_*]$.

***Proof of Lemma 9.2***: Let $c_\lozenge$ denote the larger of the respective versions of $\kappa$ that appear in Lemma 5.8 and Proposition 6.1. Set $\mathrm{T}_\ddagger$ to be the largest time in $(t^-, \frac{1}{c_\lozenge} t_*]$ such that

$$\mathrm{K}(\sqrt{t}) \ge 16c_*((1+\mathrm{K}(\tfrac{1}{c_*} t_*) + \mathrm{K}_*) \quad \textit{for all } t \in (t^-, \mathrm{T}_\ddagger]. \tag{9.2}$$

What follows directly proves that the assertions of Lemma 9.2 hold if $\mathrm{T}_*$ is replaced by $\mathrm{T}_\ddagger$ in each of the four bullets. Indeed, the lemma with this replacement follows by alternating appeals to Proposition 6.1 and Proposition 8.1 starting from $t^-$. To elaborate: The time $t^-$ is either in $\Omega^+$ or $\Omega^-$. Suppose, for example, that $t^-$ is in $\Omega^+$ (one can show that this is the case for small enough $t^-$ if the traceless Ricci curvature of Y's metric is not zero). Apply Proposition 6.1 with $t_0 = t^-$ and with $t_1$ chosen so that it is either the lower boundary point of a component of $\Omega^-$ or it is described by CASE B of Proposition 6.1. Note with regards to this application of Proposition 6.1: The number $\mathrm{N}(t_0)$ obeys $\mathrm{N}(t_0) < 1 + c_0 t_0^2$ by virtue of the fourth bullet of Lemma 5.1. Meanwhile $\mathfrak{c}^-(t_0)$ is zero; and the $\{t_0\}\times Y$ integral of $|\mathfrak{c}|$ is at most $c_0 t_0$ and that of $|\mathfrak{c}^+|$ is at most $c_0 t_0 |\ln t_0|$. (These bounds follow from Lemmas 3.2 and 3.3). The number $\mathrm{E}_0$ in this case can be assumed to be less than 1 (by virtue of the choice for $t^-$.). If CASE B of Proposition 6.1 describes $t_1$, then stop because Lemma 9.1 with the replacement $\mathrm{T}_\ddagger$ for $\mathrm{T}_*$ follows directly. If $t_1$ is the lower boundary point of a component of $\Omega^-$, then apply Proposition 8.1 to this component. The $t_0$ requirements for Proposition 8.1 are met by virtue of what is said in Items A) and B) and D) of Proposition 6.1. If the upper boundary point for the $\Omega^-$ interval is less than $\mathrm{T}_\ddagger$, then Proposition 6.1 can be applied to the abutting interval in $\Omega^+$. The assumptions needed to apply Proposition 6.1 are supplied by Proposition 8.1. This alternating of appeals to Propositions 6.1 and 8.1 can be continued until such time t where $\mathrm{K}(\sqrt{t})$ is $16c_*(1+\mathrm{K}(\frac{1}{4c_*} t_*)+\mathrm{K}_*)$ or $t = \frac{1}{c_\lozenge} t_*$, which is to say until $t = \mathrm{T}_\ddagger$.

A confession before continuing: There is some sloppiness in the preceding paragraph if the set $\Omega^+$ has a countable or even uncountable collection of components. Nothing said yet has ruled this event out. (Taken literally, the sequential appeals to Propositions 6.1 and 8.1 in the manner just described might never end.) To deal with this event, note that the conclusions of Proposition 6.1 hold on some *open* interval containing the given interval $[t_0, t_1]$. The possibly uncountable collection of these open intervals and the components of $\Omega^-$ intersect $[t^-, \mathrm{T}_\ddagger]$ to give an open cover of $[t^-, \mathrm{T}_\ddagger]$. Since this interval is compact, there is a finite subcover; and a finite subcover can be used as described in the preceding paragraph to make the appeals to Propositions 6.1 and 8.1 . In particular, only finitely many appeals need be made.

The assertions of Lemma 9.2 follow if (9.1)'s time $T_*$ obeys $\frac{1}{c_0} T_* \le T_\ddagger$. This inequality is proved directly in four parts.

*Part 1*: Suppose first that $T_\ddagger$ is in $\Omega^+$. The assumptions of Lemma 5.9 hold and as a consequence so do its conclusion, in this case with $\Theta = 0$. In particular,

$$K(s) \le \frac{2c_*}{s} \quad \textit{and} \quad \int_{\{s\}\times Y} \langle \mathfrak{a}, \mathcal{B} \rangle \ge -c_0 \frac{1}{T_\ddagger} . \tag{9.3}$$

for all $s \in [T_\ddagger, \frac{1}{c_0}\sqrt{T_\ddagger}]$. If $T_\ddagger$ is in $\Omega^-$, then the assumptions of the $t = T_\ddagger$ version of Lemma 5.9 also hold by virtue of what is said in Proposition 8.1. In particular, the $t = T_\ddagger$ version of Lemma 5.9 with $\Theta \le c_0$ leads to the inequality in (9.3) for all $s \in [T_\ddagger, \frac{1}{c_0}\sqrt{T_\ddagger}]$ with a larger version of $c_0$.

The inequality in (9.3) (with a still larger $c_0$) also holds for $s \in [\frac{1}{c_0}\sqrt{T_\ddagger}, \sqrt{T_\ddagger}]$. That this is so follows from the top bullet of Lemma 5.5 with $t = \frac{1}{c_0}\sqrt{T_\ddagger}$ and $s = \sqrt{T_\ddagger}$ because $K(\frac{1}{c_0}\sqrt{T_\ddagger})$ is at most $c_0 \frac{1}{\sqrt{T_\ddagger}}$ (the top bullet of Lemma 5.9), and because $K(s)$ for s in $[\frac{1}{c_0}\sqrt{T_\ddagger}, \sqrt{T_\ddagger}]$ is at most $(1+c_0)\, K(\frac{1}{c_0}\sqrt{T_\ddagger})$ (by virtue of Lemma 5.1) and $K(\frac{1}{c_0}\sqrt{T_\ddagger})$ is at most $c_0 \frac{1}{\sqrt{T_\ddagger}}$ by virtue of Lemma 5.9.

*Part 2*: If $T_\ddagger \in \Omega^-$, then the $\{T_\ddagger\}\times Y$ integral of $\langle \mathfrak{a}, \mathcal{B} \rangle$ is greater than $-K^2(T_\ddagger)$ (this follows from Proposition 8.1). This is automatically true of $T_\ddagger \in \Omega^+$. Therefore, given a number $m \ge 1$, there is a maximal time $t' \in [T_\ddagger, \sqrt{T_\ddagger}]$ such that the $\{s\}\times Y$ integral of $\langle \mathfrak{a}, \mathcal{B} \rangle$ is greater than or equal to $-m K^2(s)$ for all $s \in [T_\ddagger, t']$. Let $s_m$ denote this maximal $t'$. By virtue of (5.9), the function $N(s)$ on $[T_\ddagger, s_m]$ obeys the differential equation

$$\frac{d}{dt} N \le \frac{N(1-N)}{s} + c_0 m s \ , \tag{9.4}$$

which leads to the bound $N(s) \le 1 + c_0(1+m)s$ for $s \in [T_\ddagger, s_m]$ given that $N(T_\ddagger) \le 1 + c_0 T_\ddagger$. (This bound for $N(T_\ddagger)$ follows from either Assertion D of Proposition 6.1 or the fourth bullet of Proposition 8.1 as the case may be).

*Part 2*: With $N(s) \le 1 + c_0 m s$ understood to hold for all $s \in [T_\ddagger, s_m]$, it then follows by virtue of (2.11) that $K(s)$ for such s obeys

$$K(s) \geq (1 - c_0 m s)\, \tfrac{T_‡}{s}\, K(T_‡)\ .$$

(9.5)

Meanwhile, by virtue of the fourth bullet of Lemma 9.1 (applied for times $t \leq T_‡$ where is is known to hold), $N(t) \leq 1 + c_0 t$ for all $t \in [t^-, T_‡]$. Since $K(t^-) \geq (1 - c_0 t^-)\frac{1}{t^-}$, it follows (by integrating (2.11) on the interval $[t^-, T_‡]$ that $K(T_‡) \geq (1 - c_0 T_‡)\frac{1}{T_‡}$ supposing that $T_‡ \leq c_0^{-1}$. Therefore, (9.5) leads to this lower bound:

$$K(s) \geq (1 - c_0 m s)\, \tfrac{1}{s}$$

(9.6)

when $s \in [T_‡, s_m]$.

*Part* 3*:* The part proves that $s_m = \sqrt{T_‡}$ if $m \geq m_0$ (which is a number less than $c_0$) and if $T_‡$ is less than $c_0^{-1} m_0^{-2}$. To prove this, note that if $s_m < \sqrt{T_‡}$, then the right hand inequality in (9.3) and (9.6) are consistent (remember that the $\{s_m\} \times Y$ integral of $\langle \mathfrak{a}, \mathcal{B} \rangle$ is *equal* to $-mK^2(s_m)$) only in the event that

$$m (1 - c_0 m s_m)^2 \frac{1}{{s_m}^2} \leq c_0 \frac{1}{T_‡}\ .$$

(9.7)

This says that if $s_m \leq c_0^{-1} m^{-1}$, then $s_m \geq c_0^{-1} \sqrt{m} \sqrt{T_‡}$. Therefore, taking $m = c_0$, this says in effect that $s_m \geq \sqrt{T_‡}$ if $T_‡ \leq c_0^{-1}$.

*Part 4*: If $m = c_0$ and $s_m = \sqrt{T_‡}$, then $K(\sqrt{T_‡}) \geq \frac{1}{4} \frac{1}{\sqrt{T_‡}}$ (assuming that $T_‡ \leq c_0^{-1}$). Given this lower bound, the fact that (9.2) is an equality at $T_‡$ implies that

$$\sqrt{T_‡} \geq \frac{1}{4 c_{\#} c_* (1 + K(\frac{1}{c_*} t_*) + K_*)}\ .$$

(9.8)

This implies in turn that $T_‡ \geq c_0^{-1} T_*$ which is what is needed to finish the lemma's proof.

Some of the preliminary conclusions in Parts 1-4 of Lemma 9.2's proof have implications for the time $T_*$ that are exploited in a subsequent section. These are stated by the next lemma. The lemma uses $c_\oplus$ to denote the version of $\kappa$ in Lemma 9.2

**Lemma 9.3**: *There exist* $\kappa \geq 1$ *with the following significance: Supposing that* $(A, \mathfrak{a})$ *is a Nahm pole solution, define its critical time* $T_*$ *as in (9.1). If* $t \in [\frac{1}{c_\oplus} T_*, (\frac{1}{2c_\oplus})^{1/2} \sqrt{T_*}]$, then

- $(1 - \kappa t)\frac{1}{t} \leq K(t) \leq \kappa \frac{1}{t}\ .$

- $\int_{\{t\}\times Y} \langle \mathfrak{a}, \mathcal{B}\rangle \geq -\kappa \min(\frac{1}{T_*}, \kappa^2(t))$ .
- $N(t) \leq 1 + \kappa t$ .

***Proof of Lemma 9.3***: Part 3 of Lemma 9.2's proof chooses a specific value of the number $m$ (which is less than $c_0$) so that the time $s_m$ is equal to the time $\sqrt{T_\ddagger}$ with $T_\ddagger$ defined in (9.2). Therefore, if $\frac{1}{c_\oplus} T_* < T_\ddagger$ (which it is), then $\frac{1}{\sqrt{c_\oplus}} \sqrt{T_*}$ is less than $s_m$ and thus so is any t between $\frac{1}{c_\oplus} T_*$ and $\frac{1}{\sqrt{c_\oplus}} \sqrt{T_*}$. This implies that the s = t version of (9.3) holds, that the $\{t\} \times Y$ integral of $\langle \mathfrak{a}, \mathcal{B}\rangle$ is greater than $-c_0 \kappa^2(t)$, and that $N(t) \leq 1 + c_0 t$. The latter implies the lower bound $\kappa(t) \geq (1 - c_0 t) \frac{1}{t}$ (see (9.6)).

**b) Preliminary pointwise bounds for $|\mathfrak{a}|$**

The lemma in this section asserts an apriori bound for $t|\mathfrak{a}|$. The lemma also uses $c_\oplus$ to denote Lemma 9.2's version of $\kappa$.

**Lemma 9.4**: *Let* $(A, \mathfrak{a})$ *denote a Nahm pole solution. There exists* $\lambda > 1$ *that depends on an upper bound for* $1 + \kappa(\frac{1}{c_*} t_*) + \kappa_*$ *but is otherwise independent of the Nahm pole solution; it is such that* $|\mathfrak{a}| < \frac{\lambda}{t}$ *when* $t \in (0, \frac{1}{2c_\oplus} T_*]$.

***Proof of Lemma 9.4*:** The proof has three parts. By way of notation: The proof uses $\lambda_*$ to denote Lemma 9.2's version of $\lambda$.

*Part 1*: Fix $y \in Y$ and $t \in (0, \frac{1}{2c_\oplus} T_*]$. Given $r \in [0, \frac{1}{2} t]$, let $B_r$ denote the ball in $(0,2t)\times Y$ centered at $(t,y)$ with radius r. Then, let $\hat{W}(r)$ denote $\frac{1}{r^4}$ times the integral of $|\mathfrak{a}|^2$ on $B_r$. Part 3 proves (using Lemma 9.2) that $\hat{W}$ can be written as

$$\hat{W} = \frac{3\omega_0(\alpha)}{4t^2} + \mathfrak{e}_0 + \mathfrak{e}_c \tag{9.9}$$

where the notation is as follows: First, $\alpha = \frac{r}{t}$ and $\omega_0(\cdot)$ is the function

$$\alpha \to \omega_0(\alpha) = \frac{4\pi}{3} \int_{-1}^{1} \frac{(1-s^2)^{3/2}}{(1-\alpha s)^2} ds \ . \tag{9.10}$$

Note in particular that $\omega_0(\alpha) = \frac{1}{2}\pi^2 + O(\alpha)$ for small $\alpha$. Meanwhile, the terms $\mathfrak{e}_0$ and $\mathfrak{e}_c$ are as follows: What is denoted by $\mathfrak{e}_0$ accounts for the fact that the metric on B is not the flat Euclidean metric. It's norm is bounded by $c_0 \frac{r^2}{t^2}$ . What is denoted by $\mathfrak{e}_c$ accounts for the $\mathfrak{c}$ part of $\mathfrak{a}$ (the left most term on the right hand side of (9.9), the term with the function $\omega_0(\cdot)$, accounts for the $-\frac{1}{2t}\sigma$ part of $\mathfrak{a}$). The norm of $\mathfrak{e}_c$ is bounded by

$$|\mathfrak{e}_c| \le c_0\left(\tfrac{1}{r^{3/2}t^{1/2}}\sqrt{M} + \tfrac{1}{r^2}M\right) \quad \textit{where} \quad M = \int_{[\frac{1}{2}t,\, 2t]\times Y} (|\nabla_A^{\perp}\mathfrak{c}|^2 + |\mathfrak{c}|^2)\ . \tag{9.11}$$

Note in particular that $M \le c_0\lambda_*$ by virtue of Lemmas 9.2 and Lemma 4.2.

*Part 2*: When $r \in (0, \frac{1}{2}t]$, let $\hat{K}(r)$ denote $r^{-3/2}$ times the $L^2$ norm of $|\mathfrak{a}|$ on $\partial B_r$. According to what is said in Section 3a of [T2], this function obeys $\hat{K}(r_1) \ge (1 - c_0 r_1{}^2)\hat{K}(r_0)$ when $r_1 > r_0$ (and both are less than $\frac{1}{4}t$), so it is essentially non-decreasing. Keeping this in mind, and noting that

$$\hat{W}(r) = \frac{1}{r^4}\int_0^r \hat{K}(s)^2 s^3\, ds \tag{9.12}$$

up to a correction bounded in norm by $c_0 r^2\hat{W}$ (which is due to the metric not being Euclidian), the following must hold:

- *If* $\delta \in (0, 1)$ *and if* $s \le \delta r$, *then* $\hat{K}^2(s) \le \frac{4}{1-\delta^4}(1 + c_0 r^2)\hat{W}(r)$ .
- *There exists* $s \in (\delta r, r)$ *such that* $\hat{K}^2(s) \ge 4(1 - c_0 r^2)(1-\delta^4)\hat{W}(r)$

(9.13)

The top bullet is of relevance here (both are used in the next subsection): Taking $r = \frac{1}{2}t$ and taking $s = \frac{1}{2}r$ (which is $\frac{1}{4}t$) leads to this: $\hat{K}(\frac{1}{4}t) \le c_0\lambda_*\frac{1}{t}$.

With the preceding understood, invoke the top bullet of Proposition 2.1 to see that $|\mathfrak{a}|$ on the radius $\frac{1}{8}t$ ball centered at $(t, y)$ is bounded by $c_0\lambda_*\frac{1}{t}$ also.

*Part 3*: The equations (9.9)–(9.11) are obtained by writing $\mathfrak{a}$ as $-\frac{1}{2t}\sigma + \mathfrak{c}$. The term with $\omega_0(\alpha)$ is the contribution from the $-\frac{1}{2t}\sigma$ part of $\mathfrak{a}$ when the metric on Y is a flat metric, and the $\mathfrak{e}_0$ term is the adjustment to the latter to account for the metric not being flat. The $\mathfrak{e}_c$ term is the contribution from $\mathfrak{c}$. By virtue of (4.4), its norm obeys

$$|\mathfrak{e}_c| \le \frac{1}{r^4}\left(\frac{1}{t}\int_{B_r}|c| + \int_{B_r}|\mathfrak{c}|^2\right) \tag{9.14}$$

To prove that (9.14) leads to (9.11), fix $s \in (0, 2t]$ and let $B_{r,s}$ denote the intersection of $B_r$ with the slice $\{s\}\times Y$. Since the volume of $B_{r,s}$ is at most $c_0 r^3$, the inequality in (9.14) leads to this bound:

$$|\mathfrak{e}_{\mathfrak{c}}| \le c_0 \frac{1}{r^{3/2}} \frac{1}{t} \int_{\frac{1}{2}t}^{2t} \Big( \int_{\{\cdot\}\times Y} |\mathfrak{c}|^6 \Big)^{1/6} + c_0 \frac{1}{r^2} \int_{\frac{1}{2}t}^{2t} \Big( \int_{\{\cdot\}\times Y} |\mathfrak{c}|^6 \Big)^{1/3} . \tag{9.15}$$

The preceding bound and a standard Sobolev inequality leads in turn to the following:

$$|\mathfrak{e}_{\mathfrak{c}}| \le c_0 \frac{1}{r^{3/2}} \frac{1}{t} \int_{\frac{1}{2}t}^{2t} \Big( \int_{\{\cdot\}\times Y} (|\nabla_A^{\perp}\mathfrak{c}|^2 + |\mathfrak{c}|^2) \Big)^{1/2} + c_0 \frac{1}{r^2} \int_{\frac{1}{2}t}^{2t} \Big( \int_{\{\cdot\}\times Y} (|\nabla_A^{\perp}\mathfrak{c}|^2 + |\mathfrak{c}|^2) \Big) , \tag{9.16}$$

and then to the next one (which changes just the left most term on the right hand side):

$$|\mathfrak{e}_{\mathfrak{c}}| \le c_0 \frac{1}{r^{3/2}} \frac{1}{t^{1/2}} \Big( \int_{(\frac{1}{2}t, 2t]\times Y} (|\nabla_A^{\perp}\mathfrak{c}|^2 + |\mathfrak{c}|^2) \Big)^{1/2} + c_0 \frac{1}{r^2} \int_{[\frac{1}{2}t, t]\times Y} (|\nabla_A^{\perp}\mathfrak{c}|^2 + |\mathfrak{c}|^2) \tag{9.17}$$

This last inequality is (9.11).

**c) Proof of Propostion 9.1's second bullet assuming the first bullet**

The second bullet of Proposition 9.1 is a consequence of its top bullet. The four parts of this subsection explain why. To set notation, let $c_{\oplus}$ again denote the version of $\kappa$ from Lemma 9.2. Let $\lambda_*$ now denote the larger of the $\lambda$'s from Lemmas 9.2 and 9.4.

*Part* 1: Fix $t \in (0, \frac{1}{2c_{\oplus}} \mathrm{T}_*]$ and return to the context of (9.9)–(9.11). Suppose that $\delta \in (0, \frac{1}{100}]$ and that the number $\mathrm{M}$ in (9.12) is less than $\delta^{26}$. If so, then $\mathfrak{e}_{\mathfrak{c}}$ in (9.9) obeys

$$|\mathfrak{e}_{\mathfrak{c}}| \le \frac{\delta}{t^2} \quad \textit{when } r \in [\delta^8 t, \tfrac{1}{2} t]. \tag{9.18}$$

Assuming $\delta < c_0^{-1}$, then by virtue of (9.10) and (9.11), the function $\hat{\mathrm{W}}(\cdot)$ on the interval $[\delta^2 t, \delta t]$ is almost constant in that it differs from $\frac{3\pi^2}{8t^2}$ by at most $\frac{\delta}{t^2}$. Therefore, it follows from the $r = \delta t$ version of the top bullet in (9.13) and the $r = \delta^2$ version of the lower bullet that there exists $r_\delta \in [\delta^3 t, \delta^2 t]$ such that $\hat{\kappa}^2(r_\delta)$ can be written as

$$\hat{\kappa}^2(r_\delta) = \Big( \frac{3\pi^2}{2} \frac{r_\delta^2}{t^2} \Big)(1 + e_\delta) \frac{1}{r_\delta^2} \tag{9.19}$$

with the norm of $e_\delta$ obyeing $|e_\delta| < c_0\delta$.

Because $\frac{r_\delta^2}{t^2}$ is less than $c_0\delta^4$, it follows as a consequence that the $r = r_\delta$ version of the second bullet in Proposition 2.1 is in play if $\delta < c_0^{-1}$. This bullet says (in part) that $|\nabla_A \mathfrak{a}| \le c_0 \frac{1}{r_\delta^2}$ on the radius $\frac{1}{2} r_\delta$ ball centered at $(t, y)$.

*Part 2*: Fix an orthonormal frame for $T^*Y$ at the point y and use the Levi-Civita connection's parallel transport along the geodesic segments from Y to extend this frame as an orthonormal frame for $T^*Y$ over a ball in Y centered at y of radius $c_0^{-1}$. Use this transported frame to view $\langle \mathfrak{a} \otimes \mathfrak{a}\rangle$ on this ball as a symmetric section of $\otimes_2 T^*Y|_y$.

Given $r \in (0, c_0^{-1}]$ define the symmetric element $\mathfrak{m}_r$ in $(\otimes_2 T^*Y)|_y$ by the rule

$$\mathfrak{m}_r = \frac{1}{r^4} \int_{B_r} \langle \mathfrak{a} \otimes \mathfrak{a} \rangle \ . \tag{9.20}$$

Up to an almost r independent factor, this is the average of the tensor $\langle \mathfrak{a} \otimes \mathfrak{a}\rangle$ over $B_r$. The argument used to prove (9.9)–(9.11) in the preceding subsection proves the following: If $r \in [\delta^5 r_\delta, \frac{1}{2} r_\delta]$, then $\mathfrak{m}_r$ can be written as

$$\left|\frac{\omega_0(\alpha)}{4t^2}\, \mathfrak{g} - \mathfrak{m}_r\right| \le \mathfrak{e}_0 + \mathfrak{e}_c \ , \tag{9.21}$$

where $\mathfrak{e}_0$ and $\mathfrak{e}_c$ are different than their namesakes in (9.9); but like their namesakes, $|\mathfrak{e}_0| \le \frac{r_\delta^2}{t^2}$ and $|\mathfrak{e}_c|$ obeys the bounds in (9.11). Therefore, by virtue of (9.10):

$$\left|\frac{1}{4t^2}\, \mathfrak{g} - \frac{2}{\pi^2}\, \mathfrak{m}_r\right| \le c_0 \frac{\delta}{t^2} \ . \tag{9.22}$$

Given that $|\nabla_A \mathfrak{a}| \le c_0 \frac{1}{r_\delta^2}$ on $B_r$ (because $r < \frac{1}{2} r_\delta$), the bound in (9.22) can hold only in the event that the pointwise bound

$$\left|\frac{1}{4t^2}\, \mathfrak{g} - \langle \mathfrak{a} \otimes \mathfrak{a}\rangle\right| \le c_0 \frac{\delta}{t^2} + c_0\, r \frac{1}{r_\delta^2} |\mathfrak{a}| \tag{9.23}$$

holds on the whole of $B_r$. Now $r_\delta \ge \delta^3 t$ and if $r \le \delta^4 r_\delta$ (remember that r can be as small as $\delta^5 r_\delta$), then (9.23) implies (among other things) that

$$\left|\frac{3}{4t^2} - |\mathfrak{a}|^2\right| \le c_0\, \delta\left(\frac{1}{t^2} + \frac{1}{t} |\mathfrak{a}|\right). \tag{9.24}$$

on the ball $B_r$. This last bound implies in turn that $\left|\frac{3}{4t^2} - |\mathfrak{a}|^2\right| \le c_0 \delta \frac{1}{t^2}$ on $B_r$. And, (9.24) with the latter bound implies in turn that

$$|\tfrac{1}{4t^2}\,\mathfrak{g} - \langle\mathfrak{a}\otimes\mathfrak{a}\rangle| \le c_0\,\delta\,\tfrac{1}{t^2} \tag{9.25}$$

on $B_r$. The bound in (9.25) on $B_r$ implies what is asserted by the second bullet of Proposition 9.1 because it holds with t constrained only to the extent that it be less than $\frac{1}{2c_\oplus}\mathrm{T}_*$ and that $\mathrm{M}$ at t obey $\mathrm{M} < \delta^{26}$ which is guaranteed by the first bullet of Proposition 9.1 if t is less than a time determined solely by the numbers $(1+\mathrm{K}(\frac{1}{c_*}t_*)+\mathrm{K}_*)$ and $\varepsilon$.

**d) Proof of the first bullet of Proposition 9.1**

The first bullet of Proposition 9.1 is proved momentarily. The proof invokes the lemma that follows directly. By way of notation, the lemma introduces

- $\mathrm{E}^- = \int_{(0,t^-]\times Y} (|B_A|^2 + |E_A|^2 + |\nabla_A^\perp\mathfrak{a}|^2)$ .
- $\Delta^- = \frac{|\ln t^-|}{t^-}\int_{\{t^-\}\times Y} |\mathfrak{b}^\perp|^2$ .

$$\tag{9.26}$$

Here, $\mathfrak{b}^\perp$ is defined to be the part of the $\otimes^2T^*Y$ valued 1-form on $(0,t^-]\times Y$ from (3.9) that annihilates the tangents to the $(0,\infty)$ factor of $(0,\infty)\times Y$. Other notation: The lemma uses $\lambda_*$ to denote the larger of the versions of $\lambda$ from Lemmas 9.2 and 9.4. The numbers $c_\oplus$ and $\mathrm{T}_*$ are as before ($c_\oplus$ is the version of $\kappa$ from Lemma 9.2 and $\mathrm{T}_*$ is defined in (9.1).)

**Lemma 9.5**: *Let* $(A,\mathfrak{a})$ *denote a Nahm pole solution. There exists* $\lambda > \lambda_*$ *that depends on an upper bound for* $1+\mathrm{K}(\frac{1}{c_*}t_*)+\mathrm{K}_*$, *but is otherwise independent of the Nahm pole solution; it is such that*

$$\int_{[t^-,\frac{1}{2c_\oplus}\mathrm{T}_*]\times Y} |\ln t|\,(|B_A|^2 + |E_A|^2 + |\nabla_A^\perp\mathfrak{a}|^2) \le \lambda(1+\mathrm{E}^- + \Delta^- + t^-|\ln t^-|)\ .$$

The next two paragraphs prove Proposition 9.1's top bullet given that Lemma 9.5 is true.

The crucial point for proving the top bullet is this: The integrand on the left hand side of Lemma 9.5's inequality does not depend on the choice of $t^-$. (This would not be the case if the integrand referred to the splitting of $\mathfrak{a}$ as $-\frac{1}{t}\sigma+\mathfrak{c}$ because different choices of $t^-$ can induce different splittings.) Therefore: With $t_\varepsilon$ as in Section 3 for fixed small $\varepsilon$, if there is a sequence $\{t_n^- \in (0,t_\varepsilon]\}_{n=1,2,\ldots}$ with limit zero such that corresponding sequences of $\mathrm{E}^-$ and $\Delta^-$ have limit zero, or are uniformly bounded by $c_0$, then the corresponding versions of Lemma 9.5 would lead to the inequality

$$\int_{(0,T_*]\times Y} |\ln t| \, (|B_A|^2 + |E_A|^2 + |\nabla_A^{\perp}\mathfrak{a}|^2) \le \lambda(1 + c_0) \, . \tag{9.27}$$

This last inequality implies in particular that if $t \in (0, \frac{1}{2c_{\oplus}} T_*]$, then

$$\int_{(0,t]\times Y} (|B_A|^2 + |E_A|^2 + |\nabla_A^{\perp}\mathfrak{a}|^2) \le \frac{1}{|\ln t|} \lambda(1 + c_0) \, . \tag{9.28}$$

which implies the first bullet of Proposition 9.1.

Now $E^-$ as a function of $t^-$ has limit zero as $t^-$ goes to zero by virtue of the Nahm pole condition in the third bullet of Definition 1.1 (see also the second bullet of (3.1)). What about $\Delta^-$? As explained directly, there are certainly sequence $\{t_n^-\}_{n=1,2,\ldots} \subset (0, t_\varepsilon]$ with the corresponding $\Delta^-$ sequence having limit zero also. To see why, fix $\delta > 0$ and suppose for the sake of argument that there exist $t_\delta \in (0, t_\varepsilon]$ such that

$$\frac{|\ln t|}{t} \int_{\{t\}\times Y} |\mathcal{b}^{\perp}|^2 \ge \delta \tag{9.29}$$

on $(0, t_\delta]$. If this is so, then the $(0, t_\varepsilon] \times Y$ integral of $\frac{1}{t^2} |\mathcal{b}^{\perp}|^2$ would not be finite (because the function $t \to \frac{1}{|\ln t| \, t}$ is not integrable near $t = 0$.) And if the integral of $\frac{1}{t^2} |\mathcal{b}^{\perp}|^2$ were infinite on $(0, t_\varepsilon] \times Y$, then the integral of $|\nabla_A^{\perp}\mathfrak{a}|^2$ would be too, by virtue of (3.13). But that integral is finite by assumption (the Nahm pole assumption).

***Proof of Lemma 9.5***: The proof has seven parts.

*Part 1*: Note first that by virtue of Lemma 9.2, it is sufficient to bound the integral in Lemma 9.5 with the integration domain being $[t^-, \frac{1}{2c_{\oplus}} T_*]$. With this understood, fix $t \in (t^-, \frac{1}{2} T_*]$ and consider the identity

$$\frac{d}{dt}\Big(-\ln t \int_{\{t\}\times Y} \langle \mathfrak{a}, \mathcal{B} \rangle\Big) = -\frac{1}{t} \int_{\{t\}\times Y} \langle \mathfrak{a}, \mathcal{B} \rangle + (-\ln t) \frac{d}{dt}\Big(\int_{\{t\}\times Y} \langle \mathfrak{a}, \mathcal{B} \rangle\Big) \, . \tag{9.30}$$

Integrating this identity over the interval $[t^-, \frac{1}{2} T_*]$ using (4.13)–(4.15) for the right most term leads to the following inequality

$$\int_{[t^-, \frac{1}{2c_{\oplus}} T_*]\times Y} |\ln t| (|B_A|^2 + |E_A|^2 + |\nabla_A^{\perp}\mathfrak{a}|^2) \le c_0 \int_{[t^-, \frac{1}{2c_{\oplus}} T_*]\times Y} |\ln t| \, |\mathfrak{c}|^2 + c_0 \int_{[t^-, \frac{1}{2c_{\oplus}} T_*]\times Y} \frac{1}{t} \langle \mathfrak{a}, \mathcal{B} \rangle + \mathcal{A}\big(\tfrac{1}{2c_{\oplus}} T_*\big) - \mathcal{A}(t^-) \tag{9.31}$$

where $\mathcal{A}(t) = |\ln t| \int_{\{t\}\times Y} \langle \mathfrak{a}, \mathcal{B}\rangle$ .

*Part 2*: With regards to the $\mathcal{A}(t^-)$ term: It need be considered only in the event that the $\{t^-\}\times Y$ integral of $\langle \mathfrak{a}, \mathcal{B}\rangle$ is negative. But if this integral is negative, it is no smaller than $-c_0 t^-$ (see the second bullet of Lemma 3.6) and so $-\mathcal{A}(t^-) \le c_0\, t^-|\ln t^-|$.

With regards to the $\mathcal{A}(\frac{1}{2c_\oplus} T_*)$ term: It need be considered only in the event that the $\{\frac{1}{2c_\oplus} T_*\}\times Y$ integral of $\langle \mathfrak{a}, \mathcal{B}\rangle$ is positive. Since the latter integral is no greater than $c_0 \frac{\lambda_*}{T_*}(1+\kappa(\frac{1}{c_*} t_*)+\kappa_*)^2$ (see the top bullet of Lemma 6.4), it follows as a consequence that $\mathcal{A}(\frac{1}{2c_\oplus} T_*) \le c_0|\ln(T_*)|\, \frac{\lambda_*}{T_*}(1+\kappa(\frac{1}{c_*} t_*)+\kappa_*)^2$.

The rest of the proof deals with the integrals on the right hand side of (9.31).

*Part 3*: This part and Part 4 consider the integral of $|\ln t|\,|\mathfrak{c}|^2$ that appears on the right hand side of (9.31). To this end, fix $\delta \in (0, \frac{1}{100})$ and then use (4.21) to see that

$$\int_{[t^-,\frac{1}{2c_\oplus}T_*]\times Y} |\ln t|\,|\nabla_A^{\perp}\mathfrak{a}|^2 \; + \; c_0\delta \int_{[t^-,\frac{1}{2c_\oplus}T_*]\times Y} \frac{|\ln t|}{t^2}\,|\mathscr{b}|^2 \;\ge$$

$$\frac{1}{2}\int_{[t^-,\frac{1}{2c_\oplus}T_*]\times Y} |\ln t|\,|\nabla_A^{\perp}\mathfrak{a}|^2 \; + \; \delta \int_{[t^-,\frac{1}{2c_\oplus}T_*]\times Y} |\ln t|\,|\nabla_A^{\perp}\mathfrak{c}|^2 \quad . \tag{9.32}$$

Now consider the left hand side integral that has $|\mathscr{b}|^2$. Write $\frac{1}{t^2}$ as $-\frac{d}{dt}(\frac{1}{t})$ and then integrate by parts (and remember that $E_A = \frac{\partial}{\partial t}\mathscr{b}$) to see that

$$\int_{[t^-,\frac{1}{2c_\oplus}T_*]\times Y} \frac{|\ln t|}{t^2}\,|\mathscr{b}|^2 \;\le c_0\, \frac{|\ln t^-|}{t^-} \int_{\{t^-\}\times Y} |\mathscr{b}|^2 \; + c_0 \int_{[t^-,\frac{1}{2c_\oplus}T_*]\times Y} \frac{1}{t^2}\,|\mathscr{b}|^2 + c_0 \int_{[t^-,\frac{1}{2c_\oplus}T_*]\times Y} |\ln t|\,|E_A|^2 \quad . \tag{9.33}$$

With (9.32) and (9.33) in hand, then (9.31) leads to the following when $\delta < c_0^{-1}$:

$$\int_{[t^-,\frac{1}{2c_\oplus}T_*]\times Y} |\ln t|\,(|B_A|^2 + |E_A|^2 + |\nabla_A^{\perp}\mathfrak{a}|^2 + \delta\,|\nabla_A^{\perp}\mathfrak{c}|^2) \le$$

$$c_0 \int_{[t^-,\frac{1}{2c_\oplus}T_*]\times Y} |\ln t|\,|\mathfrak{c}|^2 + c_0 \int_{[t^-,\frac{1}{2c_\oplus}T_*]\times Y} \frac{1}{t}\langle \mathfrak{a}, \mathcal{B}\rangle + c_0\delta \int_{[t^-,\frac{1}{2c_\oplus}T_*]\times Y} \frac{1}{t^2}\,|\mathscr{b}|^2 \; + c_0(\Delta^- + \Delta_{\mathcal{A}})\,. \tag{9.34}$$

where $\Delta^-$ is from (9.26) and where $\Delta_{\mathcal{A}}$ is any positive upper bound for $\mathcal{A}(\frac{1}{2c_\oplus} T_*) - \mathcal{A}(t^-)$.

With regards to that $\frac{1}{t^2}|\mathscr{b}|^2$ integral on the right hand side of (9.34): It follows from Lemma 4.1 that this integral is at most $c_0$ times the sum of $\mathrm{E}^-$ and the $[t^-, \frac{1}{2c_\oplus} T_*]\times Y$ integral of $|E_A|^2 + |B_A|^2 + |\nabla_A{}^{\perp}\mathfrak{a}|^2$. Therefore, Lemma 9.2 can be used to replace (9.18) by

$$\int_{[t^-,\frac{1}{2c_\oplus}T_*]\times Y} |\ln t|\,(|B_A|^2+|E_A|^2+|\nabla_A^{\perp}\mathfrak{a}|^2+\delta|\nabla_A^{\perp}\mathfrak{c}|^2) \le c_0 \int_{[t^-,\frac{1}{2c_\oplus}T_*]\times Y} |\ln t|\,|\mathfrak{c}|^2 + c_0 \int_{[t^-,\frac{1}{2c_\oplus}T_*]\times Y} \tfrac{1}{t}\langle \mathfrak{a},\mathcal{B}\rangle + c_0\,(\Delta^- + E^- + \lambda_{*} + \Delta_{\mathcal{A}}) . \tag{9.35}$$

The next part of the proof starts from (9.35).

*Part 4*: Invoke the $\mathfrak{q} = \mathfrak{c}$ version of Lemma 2.3 to see that the $\delta = c_0^{-1}$ version (9.35) leads to the following:

$$\int_{[t^-,\frac{1}{2c_\oplus}T_*]\times Y} |\ln t|\,(|B_A|^2+|E_A|^2+|\nabla_A^{\perp}\mathfrak{a}|^2+\delta|\nabla_A^{\perp}\mathfrak{c}|^2) \le \\ c_0 \int_{[t^-,\frac{1}{2c_\oplus}T_*]\times Y} |\ln t|\,\Big(\int_{\{\cdot\}\times Y} |\mathfrak{c}|\,\Big)^2 + c_0 \int_{[t^-,\frac{1}{2c_\oplus}T_*]\times Y} \tfrac{1}{t}\langle \mathfrak{a},\mathcal{B}\rangle + c_0\,(\Delta^- + E^- + \lambda_{*} + \Delta_{\mathcal{A}}) \tag{9.36}$$

And, what with the middle bullet of Lemma 9.2, this leads in turn to the following bound:

$$\int_{[t^-,\frac{1}{2c_\oplus}T_*]\times Y} |\ln t|\,(|B_A|^2+|E_A|^2+|\nabla_A^{\perp}\mathfrak{a}|^2+\delta|\nabla_A^{\perp}\mathfrak{c}|^2) \le c_0 \int_{[t^-,\frac{1}{2c_\oplus}T_*]\times Y} \tfrac{1}{t}\langle \mathfrak{a},\mathcal{B}\rangle + c_0\,(\Delta^- + E^- + \lambda_{*}^{2} + \Delta_{\mathcal{A}}). \tag{9.37}$$

There is just the one last integral on the right hand side of (9.37) to deal with.

*Part 5*: Turn the focus now to the integral of $\frac{1}{t}\langle \mathfrak{a},\mathcal{B}\rangle$ that appears on the right hand side of (9.37). Write $\mathcal{B}$ in terms of $\mathcal{b}$ as in (4.7) and having done so, then integrate by parts on any given $\{t\} \times Y$ slice to see that

$$\Big|\int_{\{t\}\times Y} \langle \mathfrak{a},\mathcal{B}\rangle\Big| \le c_0 \int_{\{t\}\times Y} |\nabla_A^{\perp}\mathfrak{a}|\,|\mathcal{b}| + c_0 \int_{\{t\}\times Y} |\mathfrak{a}|\,|\mathcal{b}|^2 . \tag{9.38}$$

Noting that $|\mathfrak{a}| \le \frac{\lambda_*}{t}$ (see Lemma 9.4), it follows from (9.38) that

$$\int_{[t^-,\frac{1}{2}T_*]\times Y} \tfrac{1}{t}\langle \mathfrak{a},\mathcal{B}\rangle \le c_0^{-1} \int_{[t^-,\frac{1}{2}T_*]\times Y} |\nabla_A^{\perp}\mathfrak{a}|^2 + c_0\,\lambda_* \int_{[t^-,\frac{1}{2}T_*]\times Y} \tfrac{1}{t^2}\,|\mathcal{b}|^2 . \tag{9.39}$$

Using this last bound in (9.37) with what was said previously about the $\frac{1}{t^2}|\mathcal{b}|^2$ integral and what is said in Lemma 9.2 directly to the bound in Lemma 9.5.

## 10. Proof of Theorems A and B

The proof of the first bullet of Theorem A is in Section 10.2. The proof of Theorem B is in Section 10e. The proof of the second bullet of Theorem A is in Section 10g. Intervening subsections supply material for subsequent arguments.

### a) Bounds for the whole of $(0, \infty) \times Y$

Proposition 9.1 implies the following: Fix a Nahm pole solution. Then, given $\varepsilon \in (0, c_0^{-1}]$, there exists $t_\varepsilon \in (0, \frac{1}{2} t_*]$ that depends only on $\varepsilon$ and an upper bound for the number $1 + K(\frac{1}{c_*} t_*) + K_*$ such that if $t < t_\varepsilon$, then $\mathfrak{a}$ can be written as in (3.4) with $\tau$ and $\mathfrak{c}$ obeying (3.5). A crucial fact is that this time $t_\varepsilon$ depends only on $\varepsilon$ and an upper bound for $1 + K(\frac{1}{c_*} t_*) + K_*$. Otherwise, it is independent of $(A, \mathfrak{a})$. Because of this fact, a time for which the bounds supplied by the lemmas in Section 3 are valid depend on $(A, \mathfrak{a})$ only via the upper bound for $1 + K(\frac{1}{c_*} t_*) + K_*$. The proposition that follows makes a precise assertion to this effect.

**Proposition 10.1**: *There exists* $\kappa > 1$ *with the following significance: Having fixed a Nahm pole solution,* $(A, \mathfrak{a})$*, there exists* $T \in (0, \frac{1}{8} t_*]$ *that depends only on an upper bound for* $1 + K(\frac{1}{c_*} t_*) + K_*$ *but is otherwise independent of* $(A, \mathfrak{a})$. *It is such that*

- *The bounds in* (3.1) *hold on* $(0, T] \times Y$ *for* $\varepsilon = \frac{1}{1000\kappa}$.
- *The* ad(P) *valued 1-form* $\mathfrak{a}$ *on* $(0, T] \times Y$ *can be written as* $-\frac{1}{2t}\tau + \mathfrak{c}$ *with* $\tau$ *and* $\mathfrak{c}$ *as (3.4) and as in the* $\varepsilon = \frac{1}{1000\kappa}$ *version of* (3.5). *Define* $\hat{\mathfrak{c}}$ *to be* $\langle \tau \otimes \mathfrak{c} \rangle$. *Then* $\mathfrak{c}$ *and the corresponding* $\hat{\mathfrak{c}}$ *obey*
  - a) $\int_{\{t\} \times Y} |\mathfrak{c}| \leq \kappa t |\ln t|$ *when* $t \in (0, T)$ ;
  - b) $\int_{\{t\} \times Y} |\mathfrak{c}|^2 \leq \kappa t$ *when* $t \in (0, T)$ ;
  - c) $\int_{[0,T] \times Y} |\nabla_A \mathfrak{c}|^2 \leq \kappa$ *and* $\int_{[0,T] \times Y} |\nabla \hat{\mathfrak{c}}|^2 \leq \kappa$ .
- *If* $\mathfrak{b}$ *is defined as in (3.9) on* $(0, T]$*, then* $\int_{[0,T] \times Y} (|\nabla \mathfrak{b}|^2 + \frac{1}{t^2}|\mathfrak{b}|^2) \leq \kappa$.
- *Fix a non-negative integer* k *and* $t \in (0, T]$. *There is a bound on the* $C^k$ *norm of the tensors* $\langle \mathfrak{c} \otimes \tau \rangle$ *and* $\mathfrak{b}$ *on* $[t, 2t] \times Y$ *that may depend on* k *and* t*, but not in* $(A, \mathfrak{a})$.

***Proof of Proposition 10.1***: The first bullet's bounds follow directly from Proposition 9.1. To prove the assertions in the second bullet concerning $\mathfrak{c}$: Repeat the arguments for Lemmas 3.2, 3.3 and 3.4 using Proposition 9.1 to obtain bounds (uniform, given $1 + K(\frac{1}{c_*} t_*) + K_*$) for the numbers that appear in Section 3's proofs of these lemmas. To prove the third bullet concerning $\mathcal{b}$: Repeat the arguments for Lemma 3.5 using the bounds from Proposition 9.1. The fourth bullet follows from Lemma 3.7.

Let $(A, \mathfrak{a})$ denote a given Nahm pole solution. Proposition 10.1 uses $(A, \mathfrak{a})$ to define the time $\mathrm{T}$. The second proposition in this subsection considers the behavior of $(A, \mathfrak{a})$ not just on $(0, \mathrm{T}) \times Y$, but on the whole of $(0, \infty) \times Y$, which is to say the behavior where t is not necessarily small.

To set the stage: No lemma or proposition in Sections 2–9 uses Definition 1.2's second bullet, nor does Proposition 10.1. (The first bullet is used only in Lemma 2.2 to see that $a_0 = 0$). The upcoming proposition does require the second bullet constraint in Definition 1.2. By way of a reminder, this bullet requires this:

$$\int_{[1,\infty)\times Y} \left(|d_A\mathfrak{a}|^2 + |\nabla_t\mathfrak{a}|^2 + |B_A - *(\mathfrak{a}\wedge\mathfrak{a})|^2 + |E_A|^2\right) < \infty. \tag{10.1}$$

Use $(A, \mathfrak{a})$ to define the function $\mathfrak{cs}$ on $(0, \infty)$ via the rule in (1.3). As noted in Section 1 (after (1.3)), this function is increasing on $(0, \infty)$ and, by virtue of the integral in (10.1) being finite, it has a $t \to \infty$ limit. This limit is denoted by $\mathfrak{cs}_\infty$.

The upcoming proposition and the subsequent one use $c_\oplus$ to denote Lemma 9.2's version of κ and it uses a given Nahm pole solution to define the time $\mathrm{T}_*$ via (9.1)

**Proposition 10.2**: *Fix a Nahm pole solution,* $(A, \mathfrak{a})$. *Given a time* $t > 0$, *and an integer* k, *there exists* $\xi > 1$ *that depends on* t *and* k. *It also depends on*

- *An upper bound for* $1 + \mathrm{K}(\frac{1}{c_*} t_*) + \mathrm{K}_*$;
- *If* $t > \frac{1}{2c_\oplus}\mathrm{T}_*$, *then also an upper bound for the number* $\mathfrak{cs}_\infty$.

*It is otherwise independent of* $(A, \mathfrak{a})$. *There is also an isomorphism (which is independent of* k*) from the product principle* SU(2)*-bundle to the bundle* P *over* $[t, 2t] \times Y$. *The number* $\xi$ *and the isomorphism have the following significance: Denote the isomorphism by* g *and write* $g^*A$ *as* $\theta_0 + \hat{a}_A$ *with* $\theta_0$ *denoting the product connection on* P. *Denote* $g^*\mathfrak{a}$ *by* $\hat{\mathfrak{a}}$. *Then, the norms of* $\hat{a}_A$ *and* $\hat{\mathfrak{a}}$ *and their* $\theta_0$*-covariant derivatives to order* k *are bounded by* $\xi t^{-k-1}$ *on* $[t, 2t] \times Y$.

This proposition is proved momentarily. It leads in a fairly standard way to the last proposition in this subsection:

**Proposition 10.3**: *Fix a Nahm pole solution,* $(A, \mathfrak{a})$. *Given times* $t_0 < t_1$ in $(0, \infty)$ *and integer* k, *there exists* $\xi > 1$ *that depends on* $t_0, t_1$ *and* k. *It also depends on*

- *An upper bound for* $1 + \mathrm{K}(\frac{1}{c_*} t_*) + \mathrm{K}_*$;

- *An upper bound for the number* $\mathfrak{cs}_\infty$ *it* $t_1 \geq \frac{1}{2c_\oplus} T_*$

*It is otherwise independent of* $(A, \mathfrak{a})$. *There is also an isomorphism (which is independent of* $k$*) from the product principle* SU(2)*-bundle to the bundle* $P$ *over* $[t_0,t_1] \times Y$. *The number* $\xi$ *and the isomorphism have the following significance: Denote the isomorphism by* $g$ *and write* $g^*A$ *as* $\theta_0 + \hat{a}_A$ *with* $\theta_0$ *denoting the product connection on* $P$. *Denote* $g^*\mathfrak{a}$ *by* $\hat{\mathfrak{a}}$. *Then, the norms of* $\hat{a}_A$ *and* $\hat{\mathfrak{a}}$ *and their* $\theta_0$*-covariant derivatives to order* $k$ *are bounded by* $\xi$ *on* $[t_0,t_1]\times Y$.

***Proof of Proposition 10.3***: Supposing that $n \in \{0, 1, \ldots,\}$, let $g_n$ denote the isomorphism supplied by Proposition 10.2 for the region $[2^n t_0, 2^{n+1} t_0] \times Y$. Let $h_n$ denote Proposition 10.2's isomorphism for the region $[2^{n+1/2} t_0, 2^{n+3/2} t_0] \times Y$. The domain $(2^{n+1/2} t_0, 2^{n+1} t_0) \times Y$ is contained in the overlap of the two regions. The 'transition' function $h_n g_n^{-1}$ on this overlap domain is smooth with bounds on the derivatives to any given order determined by the order and an upper bound for $1 + \kappa(\frac{1}{c_*} t_*) + \kappa_*$ and for $\mathfrak{cs}$. This is also the case for the transition function $g_n h_{n-1}^{-1}$ (for $n \geq 1$) on the overlap domain $(2^n t_0, 2^{n+1/2} t_0)$. Granted this, then it is a straightforward task to modify these isomorphisms sequentially, increasing in $n$, to obtain the required isomorphism on the whole of $[t_0, t_1]$. (Stop when $2^n t_0$ is greater than $t_1$.)

***Proof of Proposition 10.2***: Invoke Lemma 9.4 for $t < \frac{1}{2c_\oplus} T_*$ to get a uniform bound on $|\mathfrak{a}|$ by $\frac{\lambda}{t}$ where $\lambda$ depends only on an upper bound for $1 + \kappa(\frac{1}{c_*} t_*) + \kappa_*$ but is otherwise independent of $(A, \mathfrak{a})$. If $t < \frac{1}{4c_\oplus} T_*$, then Proposition 2.1 (in particular, its third bullet) can be invoked on balls of radius less than $\frac{1}{100 c_0 \lambda} t$ inside $[\frac{1}{2} t, 3t] \times Y$ to obtain an isomorphism from the product principle SU(2) bundle to the bundle P over the ball with the following properties: For any such ball, let $u$ denote the corresponding isomorphism from the product SU(2) bundle to P. Set $a^u = u^*A - \theta_0$ and set $\mathfrak{a}^u = u^*\mathfrak{a}$. The pair $(a^u, \mathfrak{a}^u)$ of SU(2)-Lie algebra valued 1-forms obey uniform $C^k$ bounds on each such ball with the covariant derivative defined by $\theta_0$. This means that the $C^k$ norm is at most $\xi_k t^{-k-1}$ with $\xi_k$ being independent of $(A, \mathfrak{a})$ and the ball in question, and independent of t.

These local bundle isomorphisms can then be modified on the overlaps of balls as just described (suitably chosen, and giving an open cover of $[t, 2t] \times Y$) to obtain an isomorphism from the product principle SU(2) bundle to P on the whole of $[t, 2t] \times Y$. This sort of modification on overlaps of balls from a cover is described in [Uh]. See also the Appendix of [Tan]. The required modifications on overlaps can be done because the transition function for any pair of overlapping balls will obey uniform $C^{k+1}$ bounds that come from the $C^k$ bounds on the corresponding $\hat{a}^u$. (The transition functions in this case are automorphisms of the product SU(2) bundle; if $u$ and $u'$ are the isomorphisms for two overlapping balls, then the transition function on their intersection is $u' u^{-1}$.)

Proposition 2.1 can also be invoked for suitable radius balls and with the same sort of modifications of local isomorphisms to obtain the assertion of Proposition 10.2 for $t > \frac{1}{4} T_*$ given suitably $(A, \mathfrak{a})$ independent upper bounds for $\kappa(\cdot)$ on $[\frac{1}{2} t, 3t] \times Y$. This means that the bounds can depend only on t and on upper bounds for $1 + \kappa(\frac{1}{c_*} t_*) + \kappa_*$ and $\mathfrak{cs}_\infty$. The derivation of the desired bound for $\kappa$ has two parts.

*Part 1*: Suppose that there is a bound of the form

$$\int_{[\frac{1}{2c_\oplus} T_*, \infty) \times Y} (|d_A \mathfrak{a}|^2 + |\nabla_t \mathfrak{a}|^2 + |B_A - *(\mathfrak{a} \wedge \mathfrak{a})|^2 + |E_A|^2) \leq \zeta + \mathfrak{cs}_\infty \tag{10.2}$$

with $\zeta$ depending only on an upper bound for $1 + \kappa(\frac{1}{c_*} t_*) + \kappa_*$. Of particular interest is the bound in (10.2) on $|\nabla_t \mathfrak{a}|^2$ because that bound (with the fundamental theorem of calculus) leads (where $t \geq \frac{1}{2c_\oplus} T_*$) to this one:

$$\kappa(t) \leq \kappa(\tfrac{1}{2c_\oplus} T_*) + (t - \tfrac{1}{2c_\oplus} T_*)^{1/2} (\zeta + \mathfrak{cs})^{1/2} \tag{10.3}$$

Meanwhile, by virtue of Lemma 9.3, the number $\kappa(\frac{1}{2c_\oplus} T_*)$ is bounded by $c_0 \frac{1}{T_*}$ Therefore, (10.3) exhibits an a priori bound for $\kappa(t)$ that depends only on t and upper bounds for $1 + \kappa(\frac{1}{c_*} t_*) + \kappa_*$.

*Part 2*: To derive a bound such as that in (10.2), fix for the moment times $t \in (0, \infty)$ and $s \in (0, t)$. Then there is the identity (see (2.23))

$$\mathfrak{cs}(s) - \mathfrak{cs}(t) = \tfrac{1}{2} \int_{[t, s] \times Y} (|\nabla_t \mathfrak{a}|^2 + |d_A \mathfrak{a}|^2 + |B_A - *(\mathfrak{a} \wedge \mathfrak{a})|^2 + |E_A|^2) \ . \tag{10.4}$$

And this implies in turn that

$$\int_{[t, s] \times Y} (|d_A \mathfrak{a}|^2 + |\nabla_t \mathfrak{a}|^2 + |B_A - *(\mathfrak{a} \wedge \mathfrak{a})|^2 + |E_A|^2) \leq \mathfrak{cs}_\infty - \mathfrak{cs}(t) \ . \tag{10.5}$$

This last identity puts the focus on $\mathfrak{cs}(t)$

What with (2.9), $\mathfrak{cs}(t)$ can be written as

$$\mathfrak{cs}(t) = -\tfrac{1}{3} \tfrac{NK^2}{t} + \tfrac{2}{3} \int_{\{t\} \times Y} \langle \mathfrak{a}, B_A \rangle \tag{10.6}$$

Now, for $t = \frac{1}{2c_\oplus} T_*$, the $\{t\} \times Y$ integral of $\langle \mathfrak{a}, B_A \rangle$ is no less than $-c_0 \frac{1}{T_*}$. This follows from the second bullet of Lemma 9.2. Meanwhile the value of $N$ at t for $t = \frac{1}{2c_\oplus} T_*$ is not much bigger than 1 (see the last bullet of Lemma 9.2); and the value of $K$ at $\frac{1}{2c_\oplus} T_*$ is bounded by $\frac{1}{T_*}$ (see Lemma 9.2 again).

It follows from what was is said in the preceding paragraph that $|\mathfrak{cs}(\frac{1}{4c_\oplus} T_*)|$ has an upper bound that depends only on an upper bound for $1 + K(\frac{1}{c_*} t_*) + K_*$. This fact with (10.5) leads directly to (10.2).

**b) Proof of the first bullet of Theorem A: Limits of $K(1)$ bounded sequences**

Suppose that $\{(A_n, \mathfrak{a}_n)\}_{n=1,2,\ldots}$ is a sequence of Nahm pole solutions with the following property: There exists $\mathcal{K} > 0$ and $cs \in \mathbb{R}$ such that each of the $(A_n, \mathfrak{a}_n)$ versions of $1 + K(\frac{1}{c_*} t_*) + K_*$ are bounded from above by $\mathcal{K}$, and such that each of the $(A_n, \mathfrak{a}_n)$ versions of $\mathfrak{cs}_\infty$ are bounded above by $cs$. The following proposition makes a formal assertion to the effect that this sequence has a subsequence that converges to a Nahm pole solution. (This proposition implies what is said by the first bullet of Theorem A.)

**Proposition 10.4**: *Let $\{(A_n, \mathfrak{a}_n)\}_{n=1,2,\ldots}$ denote a sequence of Nahm pole solutions as described above for the given values of $\mathcal{K}$ and $cs$. There exists*

a) *A subsequence $\Lambda \subset \{1, 2, \ldots\}$ ;*

b) *A time $\mathcal{T} \in (0, \frac{1}{8} t_*)$ (the latter depending only on $\mathcal{K}$) ;*

c) *A Nahm pole solution $(A, \mathfrak{a})$ with $1 + K(\frac{1}{c_*} t_*) + K_* \leq \mathcal{K}$, and with $\mathfrak{cs}_\infty \leq cs$. It is described by Propositions 9.1 and 10.1 where $t < \mathcal{T}$.*

d) *A sequence $\{g_n\}_{n\in\Lambda}$ of automorphisms of $P$ over $(0, \infty) \times Y$*

*This data has the properties listed below. The list writes $\mathfrak{a}$ where $t < \mathcal{T}$ as $-\frac{1}{2}\tau + \mathfrak{c}$ and it defines $\mathfrak{b}$ from $A$ where $t < \mathcal{T}$ as in (3.9). Define $\hat{\mathfrak{c}}$ to be $\langle \tau \otimes \mathfrak{c} \rangle$.*

- *If $n \in \Lambda$, then $(A_n, \mathfrak{a}_n)$ on $(0, \mathcal{T}] \times Y$ is also described by Propositions 9.1 and 10.1. This implies, in part, that, $\mathfrak{a}_n$ on $(0, \mathcal{T}] \times Y$ can be written as $-\frac{1}{2}\tau_n + \mathfrak{c}_n$ and $\mathfrak{b}_n$ on this same domain can be defined using the $A_n$ version of (3.9). Define $\hat{\mathfrak{c}}_n$ to be $\langle \tau_n \otimes \mathfrak{c}_n \rangle$.*
- $\lim_{n\in\Lambda} \int_{[0,\mathcal{T}]\times Y} (|\nabla(\mathfrak{b} - \mathfrak{b}_n)|^2 + \frac{1}{t^2}|(\mathfrak{b} - \mathfrak{b}_n)|^2) = 0$ .
- $\lim_{n\in\Lambda} \int_{[0,\mathcal{T}]\times Y} (|\nabla(\hat{\mathfrak{c}} - \hat{\mathfrak{c}}_n)|^2 + \frac{1}{t^2}|(\hat{\mathfrak{c}} - \hat{\mathfrak{c}}_n)|^2) = 0$ .
- $\lim_{n\in\Lambda} t^2|\langle \mathfrak{a}_n \otimes \mathfrak{a}_n \rangle - \langle \mathfrak{a} \otimes \mathfrak{a} \rangle| = 0$ *for $t \in (0, \mathcal{T})$ and this limit is uniform in* $t$.
- *The sequences $\{(g_n^* A_n - A\}_{n\in\Lambda}$ and $\{(g_n^* \mathfrak{a}_n - \mathfrak{a})\}_{n\in\Lambda}$ converge to $0$ in the $C^\infty$ topology on compact subsets in $(0, \infty) \times Y$.*

By way of a remark: It makes no sense to directly compare $\mathfrak{c}_n$ to $\mathfrak{c}$ and $\tau_n$ to $\tau$ (and likewise $A_n$ to A) because none of these are, by themselves, invariant under the action of Aut(P). But, $\langle \tau_n \otimes \mathfrak{c}_n \rangle$ and $\langle \tau \otimes \mathfrak{c} \rangle$ are Aut(P) invariant (being $\otimes_2 T^*Y$ valued 1-forms), as are $\mathscr{b}_n$ and $\mathscr{b}$.

***Proof of Proposition 10.4***: This is a straightforward consequence of what is said by Propositions 9.1, 10.1 and 10.3. The details of this are left to the reader.

**c) On $K(t)$ and $\mathfrak{cs}(t)$**

The function $\mathfrak{cs}(t)$ plays a crucial role in controlling the behavior of a Nahm pole solution at times that are $\mathcal{O}(1)$ or larger. More to the point, a priori bounds for $\mathfrak{cs}$ on sequences of Nahm pole solutions at $\mathcal{O}(1)$ times are needed to invoke the results in [T?] which describe the behavior of the sequence at $\mathcal{O}(1)$ times where the corresponding sequence of $K$ values at $\mathcal{O}(1)$ times diverges.

With the preceding as motivation, the following lemma provides a preliminary lower bound for $\mathfrak{cs}$. The lemma uses $T_*$ to denote the time defined in (9.1). The important point with regards to the size of $T_*$ is that $\frac{1}{\sqrt{T_*}} = c_*(1+K(\frac{1}{c_*}t_*)+K_*)$; which is to say that $T_*$ is on the order of the inverse of the square of $(1+K(\frac{1}{c_*}t_*)+K_*)$.

**Lemma 10.5**: *There exists* $\kappa > 1$ *with the following significance: Let* $(A, \mathfrak{a})$ *denote a given Nahm pole solution. Then* $\mathfrak{cs}(t) \geq -\kappa(1+K(\frac{1}{c_*}t_*)+K_*)^{3/2}$ *where* $t \geq \kappa T_*^{1/4}$ .

The proof of this lemma is given in the next subsection. The next lemma is a consequence of Lemma 10.5.

**Lemma 10.6**: *Let* $(A, \mathfrak{a})$ *denote a given Nahm pole solution. There exists* $\kappa > 1$ *which has the following significance: The function* $K$ *obeys*

$$K(t) \geq \frac{1}{\kappa}(1+K(\tfrac{1}{c_*}t_*)+K_*) \quad \textit{where} \quad t \leq \frac{1}{\kappa}\,\frac{1}{T_*^{1/4}}\,\frac{1}{(1+\mathfrak{cs}_\infty T_*^{3/4})} \ .$$

Lemma 10.5 says in effect that if $(1+K(\frac{1}{c_*}t_*)+K_*)$ is large, then size of $K$ is on the order of this or greater when $t \leq (1+K(\frac{1}{c_*}t_*)+K_*)^{1/2}$ (supposing that $\mathfrak{cs}_\infty$ is not also large).

***Proof of Lemma 10.6***: The proof assumes that Lemma 10.5 is true. Let $c_\ddagger$ denote the version of $\kappa$ in Lemma 10.5. For any $t \geq c_\ddagger T_*^{1/4}$, one has

$$|\mathrm{K}(t) - \mathrm{K}(c_\ddagger \mathrm{T}_*^{1/4})| \le c_0\sqrt{t}\,(\mathfrak{cs}_\infty - \mathfrak{cs}(c_\ddagger \mathrm{T}_*^{1/4}))^{1/2} \tag{10.7}$$

(Use (10.4) with the fundamental theorem of calculus.) By taking t first to be $\frac{1}{c_*} t_*$ and then to be the $\mathcal{O}(t_*)$ times used to define $\mathrm{K}_*$; and then to be any other time, one learns that

$$\mathrm{K}(t) \ge c_0^{-1}(1+\mathrm{K}(\tfrac{1}{c_*}t_*) + \mathrm{K}_*) \;\textit{ for all }\; t \le c_0^{-1}(1+\mathrm{K}(\tfrac{1}{c_*}t_*) + \mathrm{K}_*)^2(\mathfrak{cs}_\infty - \mathfrak{cs}(t_\ddagger))^{-1} \tag{10.8}$$

which implies the lemma's assertion.

The lemma is also true with $\kappa \le c_0$ for times $t \le c_0^{-1}\mathrm{T}_*^{1/2}$ by virtue of what is said by Lemma 9.3. The claim will also follow for $t \in [\, c_0^{-1}\mathrm{T}_*^{1/2},\ c_\ddagger \mathrm{T}_*^{1/4}]$ if the function $\mathrm{N}$ is non-negative in this interval. If $\mathrm{N}$ is zero at some time in this interval, let $t_\ddagger$ denote the smallest such time. Since $\mathrm{N}(t_\ddagger) = 0$, the $\{t_\ddagger\}\times Y$ integral of $\langle \mathfrak{a}\wedge\mathfrak{a}\wedge\mathfrak{a}\rangle$ is equal to that of $\langle \mathfrak{a}, B_A\rangle$. Therefore, if the $\{t_\ddagger\}\times Y$ integral of $\langle \mathfrak{a}, B_A\rangle$ is positive, then $\mathfrak{cs}(t_\ddagger) > 0$; in which case the argument that leads to (10.7) leads to the bound $|\mathrm{K}(t) - \mathrm{K}(t_\ddagger)| \le c_0\sqrt{t}\,\mathfrak{cs}_\infty^{1/2}$ for all times $t > t_\ddagger$. This implies the assertion of the lemma for $t \ge t_\ddagger$; and the positivity of $\mathrm{N}$ for $t \le t_\ddagger$ implies the statement otherwise. If the $\{t_\ddagger\}\times Y$ integral of $\langle \mathfrak{a}, B_A\rangle$ is negative, then $\mathfrak{cs}(t_\ddagger)$ is $\frac{2}{3}$ of that. But, by virtue of (4.13), the $\{t_\ddagger\}\times Y$ integral of $\langle \mathfrak{a}, B_A\rangle$ is no smaller than $- c_0 \frac{1}{\mathrm{T}_*}$ so $\mathfrak{cs}(t)$ for $t \ge t_\ddagger$ is likewise no less than this. Then, arguing as in (10.7):

$$|\mathrm{K}(t_\ddagger) - \mathrm{K}(c_\# \mathrm{T}_*^{1/4})| \le c_0 \mathrm{T}_*^{1/8}(\tfrac{1}{\mathrm{T}_*})^{1/2} \le c_0(\tfrac{1}{\mathrm{T}_*})^{3/8} \le c_0(1+\mathrm{K}(\tfrac{1}{c_*}t_*) + \mathrm{K}_*)^{3/4}\,. \tag{10.9}$$

This implies that the bound in the lemma also holds at $t \in [t_\ddagger, c_\ddagger \mathrm{T}_*^{1/4}]$ with $\kappa \le c_0$ (and then for $t < t_\ddagger$ because of $\mathrm{N}$ being positive there).

### d) Proof of Lemma 10.5

The proof appeals to the following auxilliary lemma. By way of a reminder about notation, $c_\oplus$ denotes the version of $\kappa$ in Lemma 9.2.

**Lemma 10.7**: *There exists* $\kappa > 1$ *with the following significance*: *Fix a Nahm pole solution* $(A, \mathfrak{a})$ *and use it to define the time* $\mathrm{T}_*$ *as in (9.1).*

- *If* $t \in [\frac{1}{c_\oplus}\mathrm{T}_*,\ \frac{1}{2c_\oplus}\sqrt{\mathrm{T}_*}]$, *then* $\displaystyle\int_{\{t\}\times Y} \langle \mathfrak{a}, \mathcal{B}\rangle \ge -\kappa\,\frac{1}{\sqrt{\mathrm{T}_*}}$ .
- *If* $t \in [\frac{1}{2c_\oplus}\sqrt{\mathrm{T}_*},\ \frac{1}{8}t_*]$, *then* $\displaystyle\int_{\{t\}\times Y} \langle \mathfrak{a}, \mathcal{B}\rangle \ge -\kappa\big(\frac{1}{\sqrt{\mathrm{T}_*}} + \frac{t}{\mathrm{T}_*}\big)$.

Note that the absolute value of Lemma 10.7's lower bound is on the order of the square root of the absolute value of the lower bound in the second bullet of Lemma 9.3.

This lemma is proved after the proof of Lemma 10.5.

***Proof of Lemma 10.5***: The proof has six parts.

*Part 1*: Let $t_\wedge$ denote the largest time $t \in (0, t_*]$ such that

$$\int_{\{s\}\times Y} \langle \mathfrak{a}\wedge\mathfrak{a}\wedge\mathfrak{a}\rangle \geq - \int_{\{s\}\times Y} \langle \mathfrak{a}, B_A\rangle \quad \textit{for all } s \leq t. \tag{10.10}$$

The point of (10.10) is this: For $t \leq t_\wedge$,

$$\mathfrak{cs}(t) \geq -\tfrac{4}{3} \int_{\{t\}\times Y} \langle \mathfrak{a}\wedge\mathfrak{a}\wedge\mathfrak{a}\rangle \ . \tag{10.11}$$

This implies that anything that is larger than $\frac{4}{3} \int_{\{t\}\times Y} \langle \mathfrak{a}\wedge\mathfrak{a}\wedge\mathfrak{a}\rangle$ is larger than $-\mathfrak{cs}$.

*Part 2*: The formula in (2.23) for the derivative of $\mathfrak{cs}$ can be written (after some integration by parts and commuting covariant derivatives) as

$$\tfrac{d}{dt}\mathfrak{cs} = \tfrac{1}{2} \int_{\{t\}\times Y} (|\nabla_t \mathfrak{a}|^2 + |\nabla_A^\perp \mathfrak{a}|^2 + |B_A|^2 + |\mathfrak{a}\wedge\mathfrak{a}|^2 + |E_A|^2 + \langle \mathrm{Ric}, \langle \mathfrak{a}\otimes\mathfrak{a}\rangle\rangle) \ . \tag{10.12}$$

Therefore, supposing that the inequality

$$\int_{\{t\}\times Y} (|\nabla_t \mathfrak{a}|^2 + |\nabla_A^\perp \mathfrak{a}|^2 + |B_A|^2 + |\mathfrak{a}\wedge\mathfrak{a}|^2 + |E_A|^2) \geq 10^4 (\sup_Y |\mathrm{Ric}|)\, \kappa^2 \tag{10.13}$$

holds at time t, then the derivative of $\mathfrak{cs}$ at t has the lower bound

$$\tfrac{d}{dt}\mathfrak{cs} \geq \tfrac{1}{2} \int_{\{t\}\times Y} |\mathfrak{a}\wedge\mathfrak{a}|^2 \ . \tag{10.14}$$

And, if the latter bound holds; and if (10.10) holds also, then

$$\tfrac{d}{dt}|\mathfrak{cs}| \leq -\tfrac{1}{100} |\mathfrak{cs}|^{4/3} \ . \tag{10.15}$$

This last inequality follows because $|\mathfrak{a}\wedge\mathfrak{a}|^2 \geq |\langle \mathfrak{a}\wedge\mathfrak{a}\wedge\mathfrak{a}\rangle|^{4/3}$.

The question of where (10.10) and (10.13 hold is considered momentarily. (But the certainly hold for times t less than Proposition 10.1's time $\tau$ by virtue of its third bullet.) Anyway, supposing that (10.14) holds on a maximal interval $(0, t_{\wedge\wedge}]$, then integrating it leads to the apriori bound

$$\mathfrak{cs} \geq -10^6 \tfrac{1}{t^3} \tag{10.16}$$

for all $t \in (0, t_{\wedge\wedge}]$.

*Part 3*: Fix $m > 10^6$ and introduce yet another time, this one denoted by $t_m$, which is the largest t in $(0, t_*]$ such that $s^3\mathfrak{cs}(s) \geq -m$ for all $s \leq t$. This new time $t_m$ must be greater than $t_{\wedge\wedge}$ because of (10.16). Because $|\mathfrak{a} \wedge \mathfrak{a}|^2 \geq |\langle \mathfrak{a} \wedge \mathfrak{a} \wedge \mathfrak{a}\rangle|^{4/3}$, it follows that the $\{t\}\times Y$ integral of $|\mathfrak{a} \wedge \mathfrak{a}|^2$ is greater than $m^{4/3} \frac{1}{t^4}$ when $t \leq t_m$. By way of comparison: There is the upper bound $\kappa^2 \leq c_0 \frac{1}{t^2}$ where $t \leq (\frac{1}{2c_\oplus})^{1/2}\sqrt{T_*}$ (see Lemma 9.3); and then the upper bound $\kappa \leq c_0\kappa((\frac{1}{2c_\oplus})^{1/2}\sqrt{T_*}$ for $t \in [(\frac{1}{2c_\oplus})^{1/2}\sqrt{T_*}, \frac{1}{8} t_*]$ (from Lemma 5.1). These bounds have the following implication: If $t_m \leq c_0^{-1} T_*^{1/4}$, then the failure of (10.16) for times $t \leq t_m$ is due to the failure of (10.10), not to the failure of (10.13).

*Part 4*: Supposing that $t_m$ is less than $c_0^{-1} T_*^{1/4}$, then (10.10) must fail at some time $t \leq t_m$. It must also fail at $t_m$ because if (10.10) is obeyed at $t_m$ then so is (10.15) (because (10.13) is obeyed). And, if (10.15) is obeyed at $t_m$, then the derivative of $t^3\mathfrak{cs}(t)$ at $t = t_m$ is positive (this follows if $m > 300$). The latter event requires that $t^3\mathfrak{cs}(t)$ be less then $-m$ on some open interval with upper endpoint $t_m$ which is nonsensical given the definition of $t_m$. Therefore if $t_m \leq c_0^{-1} T_*^{1/4}$, then (10.10) must fail at $t = t_m$. This is to say that

$$\int_{\{t_m\}\times Y} \langle \mathfrak{a}, B_A\rangle < -\int_{\{t_m\}\times Y} \langle \mathfrak{a}\wedge\mathfrak{a}\wedge\mathfrak{a}\rangle . \tag{10.17}$$

This requires that $\int_{\{t_m\}\times Y} \langle \mathfrak{a}, B_A\rangle \leq -c_0 m \frac{1}{t_m^3}$ because of (2.22) and the definition of $t_m$.

*Part 5*: Choose $m = 2\cdot 10^6$. The first point to note is that if $t_m \leq c_0^{-1}T_*^{1/6}$, then $t_m$ must be in the set $\Omega^-$ (which is to say that the $\{t_m\}\times Y$ integral of $\langle \mathfrak{a}, \mathcal{B}\rangle$ is negative). Indeed, if $t_m$ is not in $\Omega^-$, then the $\{t_m\} \times Y$ integral of $\langle \mathfrak{a}, B_A\rangle$ is no smaller than $-c_0 \frac{1}{t_m}$ if $t_m \leq \frac{1}{2c_\oplus} \sqrt{T_*}$, and it is no smaller than $-c_0\sqrt{T_m}$ if $t_m \in [\frac{1}{2c_\oplus} \sqrt{T_*}, \frac{1}{8} t_*]$. If these are to equal to $-m \frac{1}{t_m^3}$, then $t_m$ must be greater than $c_0^{-1}T_*^{1/6}$.

*Part 6*: Suppose that the $\{t_m\}\times Y$ integral of $\langle \mathfrak{a}, \mathcal{B}\rangle$ is negative. As explained directly, the number $t_m$ in this case is not less than $\frac{1}{c_\oplus} T_*$. Indeed, supposing that $t < \frac{1}{c_\oplus} T_*$ and that t is a point in $\Omega^-$, then the same argument that led to Lemma 9.2 (repeat the argument but stop at t) proves that both the $[0, t] \times Y$ integral of $|\mathfrak{c}|^2$ and that the $\{t\} \times Y$ integral of $|\mathfrak{c}|$ are bounded by $c_0$. (The argument for Lemma 9.3 but stopping at t proves that Lemma 9.3 holds at t with $\lambda < c_0$.) These $|\mathfrak{c}|^2$ and $|\mathfrak{c}|$ bounds with Lemma 3.6 and the top bullet of Lemma 5.5 (take t there to be $t^-$ and take s there to be the time t here) imply that the $\{t\}\times Y$ integral of $\langle \mathfrak{a}, B_A\rangle$ is greater than $-c_0 \frac{1}{t}$. This forbids the case $t_m \leq \frac{1}{c_\oplus} T_*$ unless $t_m$ and thus $T_*$ are greater than $c_0^{-1}$.

Now suppose that $t_m \geq \frac{1}{c_\oplus} T_*$. Then, by virtue of Lemma 10.7 the time $t_m$ obeys

- $\frac{1}{t_m^3} \leq c_0 \frac{1}{t_m}$ *if* $t_m \leq \frac{1}{2c_\oplus} \sqrt{T_*}$ .
- $\frac{1}{t_m^3} \leq c_0 (\frac{1}{\sqrt{T_*}} + \frac{t_m}{T_*})$ *if* $t_m \in [\frac{1}{2c_\oplus} \sqrt{T_*}, \frac{1}{8} t_*]$ .

(10.18)

And, these bounds can't be satisfied unless $t_m \geq c_0^{-1} T_*^{1/4}$ unless both are greater than $c_0^{-1}$.

***Proof of Lemma 10.7***: The second bullet follows from the first using the second bullet of Lemma 5.5. As for the top bullet: This bullet is true if the $\{t\}\times Y$ integral of $\langle \mathfrak{a}, \mathcal{B}\rangle$ is positive on the whole of $[\frac{1}{c_\oplus} T_*, \frac{1}{2c_\oplus} \sqrt{T_*}]$, so there is no generality lost to assume that there exists a time t in this interval where $\{t\}\times Y$ integral of $\langle \mathfrak{a}, \mathcal{B}\rangle$ is negative. Under these circumstances, the number $\lambda$ in Lemma 9.2 is bounded by $c_0$. It then follows by virtue of the top bullet of Lemma 9.2 and Lemma 3.6 and the top bullet of Lemma 5.5 (with $t = t^-$ and $s = \frac{1}{c_\oplus} T_*$) that the $\{\frac{1}{c_\oplus} T_*\} \times Y$ integral of $\langle \mathfrak{a}, \mathcal{B}\rangle$ is no less than $-c_0$. With this understood, then the top bullet of Lemma 10.7 follows from the top bullet of Lemma 10.5 (with $t = \frac{1}{c_\oplus} T_*$ and $s \in [\frac{1}{c_\oplus} T_*, \frac{1}{2c_\oplus} \sqrt{T_*}]$) provided that

$$\int_{[\frac{1}{c_\oplus} T_*, t] \times Y} (|\mathfrak{c}^+|^2 + |\mathfrak{c}^-|^2) \leq c_0 \frac{t}{T_*} \quad \textit{for } t \in [\frac{1}{c_\oplus} T_*, c_0^{-1} \sqrt{T_*}] \,. \tag{10.19}$$

The proof that (10.19) holds has four parts.

*Part 1*: The plan is to redo arguments for Lemma 6.7. Those arguments require bounds for $\kappa^2 - \frac{3}{4t^2}$ for $t \in [\frac{1}{c_\oplus} T_*, \frac{1}{2c_\oplus} \sqrt{T_*}]$ which are obtained by reworking the proof Lemma 6.6. Those arguments require, in turn, an a priori bound for the $\{t\}\times Y$ integral of $|\mathfrak{t}|$ when t is in $[\frac{1}{c_\oplus} T_*, \frac{1}{2c_\oplus} \sqrt{T_*}]$. A suitable bound is this one:

$$\int_{\{t\}\times Y} |t| \le c_0 \tfrac{1}{T_*} \quad \textit{for} \quad t \in [\tfrac{1}{c_\oplus} T_*, (\tfrac{1}{2c_\oplus})^{1/2}\sqrt{T_*}] \ . \tag{10.20}$$

It is obtained from the second bullet of Lemma 5.7 using bounds from Lemmas 9.2 and 9.3. (Use the second bullet of Lemma 5.7 with $s = (\frac{1}{2c_\oplus} T_*)^{1/2}$.)

*Part 2*: With the bound from (10.20) in hand, now return to the proof of Lemma 6.6. Reworking the derivation of (6.5) leads to this replacement for (6.5):

$$K^2 \le \tfrac{3}{4^{1/3}} \Big( \int_{\{t\}\times Y} \mathcal{Q} + \Delta\Big)^{2/3} + c_0 \tfrac{1}{T_*} \tag{10.21}$$

with $\Delta$ again defined by (6.6). The inequality in (6.7) is replaced in this case by

$$\int_{\{t\}\times Y} \mathcal{Q} \le \tfrac{N K^2}{3t} + c_0 \tfrac{1}{T_*} \tag{10.22}$$

The term $c_0 \frac{1}{T_*}$ here is an upper bound for the $\{t\}\times Y$ integral of $\langle \mathfrak{a}, \mathcal{B}\rangle$. It comes about by using the bound from the second bullet of Lemma 5.5 for the $\{t\}\times Y$ integral of $\langle \mathfrak{a}, \mathcal{B}\rangle$ in the event that the latter integral is positive with the time s taken to be $(\frac{1}{2c_\oplus})^{1/2}\sqrt{T_*}$. Note that the $\{s\}\times Y$ integral of $\langle \mathfrak{a}, \mathcal{B}\rangle$ is at most $c_0 K^2(s)$ which less than $c_0 \frac{1}{T_*}$ in this case (see Lemma 9.3).

What with (10.20) and (10.22), then (6.8) is replaced by the inequality

$$K^2 \le \tfrac{3}{4^{1/3}} \Big(\tfrac{N K^2}{3t} + \Delta\Big)^{2/3} + c_0 \tfrac{1}{T_*} \ . \tag{10.23}$$

With regards to $N$ here: Lemma 9.3's third bullet bounds it by $(1 + c_0 t)$. With regards to $\Delta$: Use (6.14) and (6.15) with $t = \frac{1}{c_\oplus} T_*$ with the definition of $T_*$ from (9.1) to see that

$$\Delta(t) \le c_0 \tfrac{1}{T_*} \tag{10.24}$$

when $t \in [\frac{1}{c_\oplus} T_*, \frac{1}{2c_\oplus}\sqrt{T_*}]$.

*Part 3*: To continue reworking Lemma 6.6, consider now (10.23) at a time t from the interval $t \in [\frac{1}{c_\oplus} T_*, \frac{1}{2c_\oplus}\sqrt{T_*}]$ where $K^2(t) \ge \frac{3}{4t^2}$. At such a time, the inequality in (10.23) implies that

$$\kappa^2 \leq (\tfrac{3}{4t^2})^{1/3} \kappa^{4/3} (1 + c_0 t + c_0 t^2 \tfrac{1}{T_*})^{2/3} .$$

(10.25)

This last inequality implies that $\kappa^2 - \frac{3}{4t^2} \leq c_0 \frac{1}{T_*}$ for $t \in [\frac{1}{c_\oplus} T_*, c_0^{-1}\sqrt{T_*}]$.

*Part 4*: With the bound $\kappa^2 - \frac{3}{4t^2} \leq c_0 \frac{1}{T_*}$ in hand, the proof of Lemma 6.7 can be repeated with only notational changes to conclude (with Lemma 9.2 to bound the $\{\frac{1}{c_\oplus} T_*\} \times Y$ integral of $|\mathcal{c}|$ by $c_0$) that

$$\int_{[\frac{1}{c_\oplus} T_*, t] \times Y} (|\mathfrak{c}^+|^2 + |\mathfrak{c}^-|^2) \leq c_0 \frac{t}{T_*} \quad \textit{for } t \in [\tfrac{1}{c_\oplus} T_*, (\tfrac{1}{2c_\oplus})^{1/2}\sqrt{T_*}] .$$

(10.26)

(With regards to the $\{\frac{1}{c_\oplus} T_*\} \times Y$ integral of $|\mathcal{c}|$: Keep in mind that $\lambda$ in Lemma 9.3 is bounded by $c_0$ because there is, by assumption, a time in $[\frac{1}{c_\oplus} T_*, (\frac{1}{2c_\oplus})^{1/2}\sqrt{T_*}]$ where the $\{t\} \times Y$ integral of $\langle \mathfrak{a}, \mathcal{B} \rangle$ is negative.)

**e) Proof of Theorem B: Positive Ricci curvature**

Theorem B says in effect that if the metric on Y has positive Ricci curvature, then any sequence of Nahm pole solutions has a subsequence that converges to a Nahm pole solution after applying a suitable automorphism of P on $(0, \infty) \times Y$ to each term. This theorem is an immediate consequence of Proposition 10.4 given what is said by the following lemma.

**Lemma 10.8**: *Supposing that the metric on* Y *has positive Ricci curvature, then there exists* $\kappa > 1$ *with the following significance: Let* $(A, \mathfrak{a})$ *denote a given Nahm pole solution. Then the corresponding value of* $1 + \kappa(\frac{1}{c_*} t_*) + \kappa_*$ *is no greater than* $\kappa$.

***Proof of Lemma 10.8***: When the Ricci curvature is positive, the identity in (2.10) leads to the differential inequality

$$-\frac{d^2}{dt^2} \kappa + r^2 \kappa \leq 0$$

(10.27)

with r being positive. This has the following consequence: Given positive times t and s, with $s > t$, and given $x \in (t, s)$, then the value of $\kappa$ at x obeys

$$K(x) \leq K(t)\, e^{-r(x-t)} + K(s)\, e^{-r(s-x)} .$$

(10.28)

Now fix x and use (10.7) with t replaced by s to see that $\lim_{s\to\infty} K(s)e^{-rs} = 0$. Take this limit in (10.28) to see that

$$K(x) \leq K(t)\; e^{-r(x-t)} .$$

(10.29)

Now let $c_!$ denote the version of $\kappa$ from Lemma 10.6. Take x in (10.29) to be the time $\frac{1}{c_!} \frac{1}{T_*^{1/4}} \frac{1}{(1+\mathfrak{cs}_\infty T_*^{3/4})}$ . Take t to be the time $(\frac{1}{2c_\oplus})^{1/2}\sqrt{T_*}$. Then $K(t) \geq c_* (1+K(\frac{1}{c_*} t_*)+K_*)$ (see Lemma 9.3) and $K(x) \geq c_0^{-1}(1+K(\frac{1}{c_*} t_*)+K_*)$ (see Lemma 10.6). These bounds are not compatible with (10.29) unless $T_* \geq c_0^{-1}$ which is to say that $(1+K(\frac{1}{c_*} t_*)+K_*) \leq c_0$.

With regards to positive Ricci curvature: The proof just given for Theorem B only uses the positive Ricci curvature assumption at the very end, in Lemma 10.8. But, note that the proofs of many of the preliminary results in Sections 4-9 that lead to Proposition 10.4 and then to Theorem B after Lemma 10.8 can be simplified (or the result is irrelevant) in the event that the Ricci curvature is positive.

### f) Renormalization

The proof of the second bullet of Theorem A in the next subsection uses the observations in this section about Nahm pole solutions with $(1+K(\frac{1}{c_*} t_*)+K_*)$ large (which is to say that $T_*$ is small). By virtue of what is said in Lemmas 10.3 and 10.6, the function $K$ is nearly constant on an interval that stretches from $\sqrt{T_*}$ to nearly $\frac{1}{T_*^{1/4}}$ . (The henceforth implicit assumption in what follows is that $\mathfrak{cs}_\infty$ is fixed and that $T_*^{3/4}\mathfrak{cs}_\infty \leq 1$.) To be precise, the two lemmas imply that $c_0^{-1} \leq K(t)/K(s) \leq c_0$ on this interval. Keeping this in mind, introduce by way of notation $\hat{\mathfrak{a}} = \frac{1}{K(1)} \mathfrak{a}$. The following lemma summarizes some basic facts about $\hat{\mathfrak{a}}$.

**Lemma 10.9**: *There exists* $\kappa > 100$ *with the following significance: Suppose that* $(A, \mathfrak{a})$ *is a Nahm pole solution with* $T_* \leq \frac{1}{100} t_*$ *and* $T_*^{3/4}\mathfrak{cs}_\infty \leq 1$. *If* $t \in [\sqrt{T_*}, \frac{1}{\kappa T_*^{1/4}}]$, *then*

- $\kappa^{-1} \leq \int_{\{t\}\times Y} |\hat{\mathfrak{a}}|^2 \leq \kappa$ ,
- $\int_{[t,t+1]\times Y} |d_A\hat{\mathfrak{a}}|^2 \leq \kappa T_*^{1/4}$ ,
- $\int_{[t,t+1]\times Y} |\hat{\mathfrak{a}}\wedge\hat{\mathfrak{a}}|^2 \leq \kappa T_*$ ,

- $\int_{[t,t+1]\times Y} |\nabla_A^{\perp}\hat{\mathfrak{a}}|^2 \le \kappa$ *(this is* $\kappa T_*^{1/4}$ *if the Ricci curvature of* Y *is non-negative),*
- $\int_{[t,t+1]\times Y} |\nabla_t\hat{\mathfrak{a}}|^2 \le \kappa T_*^{1/4}$ ,

*In addition* $|\hat{\mathfrak{a}}| \le \kappa(1 + \frac{1}{t^{3/2}})$ *if* $t \in [2\sqrt{T_*}, \frac{1}{2\kappa T_*^{1/4}}]$.

Lemma 10.9 is proved momentarily. What follows directly sets the stage for the next lemma which is little more than a corollary to Lemma 10.9.

Introduce $c_\#$ to denote Lemma 10.9's version of $\kappa$. Now define

$$\mathcal{T} = \langle \hat{\mathfrak{a}} \otimes \hat{\mathfrak{a}} \rangle \tag{10.30}$$

which is viewed as a symmetric section of $\otimes_2 T^*Y$ over $[2\sqrt{T_*}, \frac{1}{2c_\# T_*^{1/4}}]$. The upcoming lemma concerns the tensor $\mathcal{T}$.

This upcoming lemma reintroduces (from Section 1b) two differential operators on Y that act on symmetric tensors. The first is denoted by div. It is the formal, $L^2$ adjoint of the map from 1-forms to 2-tensors that sends any given 1-form $u$ to its covariant derivative (which is a priori a section of $\otimes_2 T^*Y$. When written using a local orthonormal frame $\{e^i\}_{i=1,2,3}$ for $T^*Y$, the 1-form $\mathrm{div}(\mathcal{T})$ is $\nabla_k \mathcal{T}_{ki}\, e^i$ where $\nabla_k$ denotes the directional covariant derivative along the vector field dual to $e^k$. The second differential operator is denoted by curl. This sends symmetric 2-forms to sections of $\otimes_2 T^*Y$. When written using the orthonormal frame, $\mathrm{curl}(\mathcal{T})$ has $e^i \otimes e^j$ component $\varepsilon^{imn}\nabla_m \mathcal{T}_{nj}$ with $\varepsilon$ denoting the volume 3-form. A final piece of notation: If $\mathcal{V}$ is a section of $\otimes_2 T^*Y$, then $\mathcal{T}\bullet\mathcal{V}$ is the section with $e^i \otimes e^j$ component given by $\mathcal{T}_{im}\mathcal{V}_{mj}$.

**Lemma 10.10**: *There exists* $\kappa > 1$ *with the following significance: Suppose that* $(A, \mathfrak{a})$ *is a Nahm pole solution with* $T_* \le \frac{1}{100} t_*$ *and* $T_*^{3/4}\mathfrak{cs}_\infty \le 1$. *If* $t \in [4\sqrt{T_*}, \frac{1}{4\kappa T_*^{1/4}}]$, *then the tensor* $\mathcal{T}$ *defined in (10.30) has the properties listed below.*

- $|\mathcal{T}| \le \kappa(1 + \frac{1}{t^3})$ .
- $\int_{[t,t+1]\times Y} |\nabla^{\perp}\mathcal{T}|^2 \le \kappa(1 + \frac{1}{t^3})$ *(this is* $\kappa T_*^{1/4}$ *if the Ricci curvature of* Y *is non-negative),*
- $\int_{[t,t+1]\times Y} |\nabla_t \mathcal{T}|^2 \le \kappa(1 + \frac{1}{t^3})\, T_*^{1/4}$.
- $\int_{[t,t+1]\times Y} |\mathrm{tr}(\mathcal{T}\bullet\mathcal{T}) - (\mathrm{tr}\mathcal{T})^2| \le \kappa T_*$
- $\int_{[t,t+1]\times Y} |\mathrm{div}\mathcal{T} - \frac{1}{2} d\, \mathrm{tr}(\mathcal{T})|^2 \le \kappa\, (1 + \frac{1}{t^3})\, T_*^{1/4}$.

- $\int_{[t,t+1]\times Y} |\mathcal{T} \bullet \operatorname{curl}(\mathcal{T})| \le \kappa\, (1 + \frac{1}{t^3})\, T_*^{1/8}$ .

***Proof of Lemma 10.10***: The first bullet follows from the $|\hat{\mathfrak{a}}|$ bound in Lemma 10.9, and the second and third follows from the $|\hat{\mathfrak{a}}|$ bound and the bound in the lemma's fourth and fifth bullets. The fourth bullet of this lemma directly restates the third bullet of Lemma 10.9. The fifth bullet follows directly from the bound in Lemma 10.9's second bullet given that $d_A{*}\hat{\mathfrak{a}} = 0$ (which is so because $d_A{*}\mathfrak{a} = 0$). The sixth bullet follows from Lemma 10.9's second and third bullets because the $e^i \otimes e^j$ component of $\mathcal{T} \bullet \operatorname{curl}(\mathcal{T})$ is

$$\langle \hat{\mathfrak{a}}_i \hat{\mathfrak{a}}_m \rangle \langle ({*}d_A \hat{\mathfrak{a}})_m\, \hat{\mathfrak{a}}_j \rangle + \varepsilon^{mnk} \langle \hat{\mathfrak{a}}_i \hat{\mathfrak{a}}_m \rangle \langle \hat{\mathfrak{a}}_k \nabla_{An} \hat{\mathfrak{a}}_j \rangle\, ; \tag{10.31}$$

and $\varepsilon^{mnk} \langle \hat{\mathfrak{a}}_i \hat{\mathfrak{a}}_m \rangle \langle \hat{\mathfrak{a}}_k \nabla_{An} \hat{\mathfrak{a}}_j \rangle$ can be written as $\frac{1}{4} \langle \hat{\mathfrak{a}}_i [\hat{\mathfrak{a}}_k, [\hat{\mathfrak{a}}_m, \nabla_{An} \hat{\mathfrak{a}}_j]] \rangle$ which is the same as $\frac{1}{4} \langle [\hat{\mathfrak{a}}_i, \hat{\mathfrak{a}}_k][\hat{\mathfrak{a}}_m, \nabla_{An} \hat{\mathfrak{a}}_j] \rangle$ whose norm is bounded by $c_0 |\hat{\mathfrak{a}} \wedge \hat{\mathfrak{a}}|\, |\hat{\mathfrak{a}}|\, |\nabla_A^{\perp} \hat{\mathfrak{a}}|$.

***Proof of Lemma 10.9***: The first bullet follows from Lemmas 10.3 and 10.6. The second and the fifth bullets follow from (10.5) and Lemma 10.5. The third and fourth bullets follow from Lemma 10.5 and (10.12).

The only assertion remaining is that of the pointwise bound. To establish that, fix $c_+$ so that the five bullets of the lemma hold with $\kappa = c_+$. Now fix $t \in [2\sqrt{T_*}, \frac{1}{2c_+ T_*^{1/4}}]$ and a point $p \in Y$. Having done so, let $G_{(t,p)}$ denote the Dirichelet Green's function for the operator $-\nabla_t^2 + d^\dagger d$ on $[\frac{1}{2} t, \frac{3}{2} t] \times Y$ with pole at the point $(t, p)$. This Green's function obeys the bounds

$$|G_{t,p}| \le c_0 \frac{1}{\operatorname{dist}(\cdot, (t,p))^2} \quad \textit{and} \quad |\nabla_t G_{(t,p)}| \le c_0 \frac{1}{\operatorname{dist}(\cdot, (t,p))^3} \tag{10.32}$$

on $[\frac{1}{2} t, \frac{3}{2} t] \times Y$. It is also positive on $(\frac{1}{2} t, \frac{3}{2} t) \times Y$. Multiply both sides of (2.4) by $G_{(t,p)}$ and then integrate over $[\frac{1}{2} t, \frac{3}{2} t] \times Y$. Two applications of integration by parts and an appeal to the top bullet of Lemma 10.9 leads to this:

$$\frac{1}{2} |\hat{\mathfrak{a}}(t,p)|^2 + \int_{[\frac{1}{2}t, \frac{3}{2}t]\times Y} G_{(t,p)} (|\nabla_t \hat{\mathfrak{a}}|^2 + |\nabla_A^{\perp} \hat{\mathfrak{a}}|^2 + \tfrac{1}{c_0 T_*} |\hat{\mathfrak{a}} \wedge \hat{\mathfrak{a}}|^2) \le c_0 \int_{[\frac{1}{2}t, \frac{3}{2}t]\times Y} G_{(t,p)} |\hat{\mathfrak{a}}|^2 + c_0 \frac{1}{t^3}\ . \tag{10.33}$$

The $c_0 \frac{1}{t^3}$ term on the right hand side invokes the bound in (10.30) for $\nabla_t G_{(t,p)}$ on the two boundary components of integration domain. As explained in the next paragraph, the integral on the right is no greater than $c_0 \frac{1}{t^2}$ . Granted this, then (10.33) implies that

$$\tfrac{1}{2}\,|\hat{\mathfrak{a}}\,(t,p)|^2 + \int_{[\frac{1}{2}t,\,\frac{3}{2}t]\times Y} G_{(t,p)}(|\nabla_t\hat{\mathfrak{a}}|^2 + \;|\nabla^{\perp}_{A}\hat{\mathfrak{a}}|^2 + \tfrac{1}{c_0 T_*}\,|\hat{\mathfrak{a}}\wedge\hat{\mathfrak{a}}|^2) \;\le c_0\,\tfrac{1}{t^3}\;\;,$$

(10.34)

which gives the lemma's sixth bullet because $G_{(t,p)} \ge 0$.

To see about the size of the integral on the right hand side of (10.33), note first that the integrals of both $|\hat{\mathfrak{a}}|^2$ and $|\nabla|\hat{\mathfrak{a}}\,||^2$ restricted to any $[s, s+1]\times Y$ for $s \in [\sqrt{T_*},\, \frac{1}{c_+ T_*^{1/4}}\,]$ are bounded by $c_0$. This is by virtue of the first, fourth and fifth bullets of Lemma 10.9. As a consequence, an appeal to what is a version of Hardy's inequality proves this: Let p denote any given point in Y. Let $B \subset [s, s+1]\times Y$ denote the ball centered at the point $o = (\frac{3}{2}\,s, y)$ with radius $\frac{1}{4}\,s$. Then

$$\int_B \tfrac{1}{\mathrm{dist}(\cdot\,, o)^2}\,|\hat{\mathfrak{a}}|^2 \;\;\le c_0\,\tfrac{1}{t^2}\;.$$

(10.35)

(The proof of this version of Hardy's inequality amounts to little more than an integration by parts with respect to the radial coordinate function on a Gaussian coordinate chart for $[s, s+1]\times Y$ centered at $(\frac{3}{2}\,s, p)$.)

**g) Proof of the second bullet of Theorem A: Limits of $\kappa(1)$ diverging sequences**

Returning to the context and notation of Theorem A, suppose that $\{(A_n, \mathfrak{a}_n)\}_{n=1,2,\ldots}$ is a sequence of Nahm pole solutions whose corresponding $\kappa(1)$ sequence has no bounded subsequences. (The latter is the numerical sequence whose n'th term is the $\mathfrak{a}_n$ version of the value of $\kappa$ at $t = 1$.) Define the corresponding sequence $\{\hat{\mathfrak{a}}_n\}_{n=1,2,\ldots}$ and then the corresponding sequence $\{\mathcal{T}_n = \langle\hat{\mathfrak{a}}_n \otimes \hat{\mathfrak{a}}_n\rangle\}_{n=1,2,\ldots}$. (The subscript n here is not a 1-form index.) By virtue of Lemma 10.10's first, second and third bullets, there is a subsequence $\Lambda \subset \{1, 2, \ldots\}$ such that the corresponding sequence $\{\mathcal{T}_n\}_{n\in\Lambda}$ converges weakly in the $L^2_1$-Sobolev topology on compact subsets of $(0, \infty)\times Y$ to a t-independent, Sobolev class $L^2_1 \cap L^\infty$ symmetric section of $\otimes_2 T^*Y$. Let $\mathcal{T}$ denote the limit. This $\mathcal{T}$ has almost everywhere rank 1 because of the Lemma 10.10's fourth bullet. Furthermore, $\mathrm{div}(\mathcal{T}) = \frac{1}{2}\,d\,\mathrm{tr}(\mathcal{T})$ by virtue of Lemma 10.10's fifth bullet, and $\mathcal{T}\bullet\mathrm{curl}(\mathcal{T}) = 0$. Thus, $\mathcal{T}$ defines a weak $\mathbb{Z}/2$ harmonic 1-form on Y. (This limit $\mathcal{T}$ is covariantly constant in the event that Y has negative Ricci curvature because of the second bullet of Lemma 10.10.)

**Appendix**: **The Mazzeo-Witten definition of the Nahm pole**

The Mazzeo and Witten [MW] definition of the Nahm pole $t \to 0$ asymptotics is not the one given by Definition 1.1 (which is to say (3.1)). To say more about the [MW] definition, let (A, $\mathfrak{a}$) denote a given solution to (2.7). The definition in [MW] of the Nahm pole solution requires a writing of the ad(P) valued 1-form $\mathfrak{a}$ as

$$\mathfrak{a} = -\tfrac{1}{2t}\upsilon + \mathfrak{e} \tag{A.1}$$

with $\upsilon$ being an isomorphism from TY to ad(P) having the properties that are described momentarily. To set the stage, introduce

$$\hat{\mathfrak{e}} = \langle \mathfrak{e} \otimes \upsilon \rangle \,, \tag{A.2}$$

a section of $\otimes_2 T^*Y$. It is the analog of $\langle \mathfrak{c} \otimes \tau \rangle$ from (3.6) (this will be denoted here by $\hat{\mathfrak{c}}$). Now introduce the $(\otimes_2 T^*Y)$–valued 1-form

$$\mathfrak{B} = \langle \upsilon \otimes \nabla_A \upsilon \rangle \tag{A.3}$$

It is the analog of $\mathscr{b}$ from (3.9).

Now for the Nahm pole requirement in [MW]: There exist $\lambda > 0$ and a writing of $\mathfrak{a}$ as in (A.1) on $(0, \varepsilon) \times Y$ with $\hat{\mathfrak{e}}$ and $\mathfrak{B}$ obeying

$$\lim_{t \to 0} t^{\lambda}(|\mathfrak{B}| + |\hat{\mathfrak{e}}|) = 0 \textit{ for some } \lambda < 1. \tag{A.4}$$

This is Condition (1) in Section 2.4.1 of [MW]. (There is also an implicit slice condition near $t = 0$ that is assumed in [MW] which is their (2.23). But the techniques in [He] can be used to prove that this slice condition can be imposed by acting via an automorphism of P defined near $t = 0$ if (A.4) holds.) The results in [MW] (see, e.g. Proposition 5.9 of [MW]) imply the following: Suppose that $(A, \mathfrak{a})$ obeys (2.7) and is such that $(A, \mathfrak{a})$ is defined where t is small by (A.1) and (A.3) using an isomorphism $\upsilon$ and tensors $\hat{\mathfrak{e}}$ and $\mathfrak{B}$ that obey (A.4). Then (3.1) holds for any given $\varepsilon$.

With the preceding understood, the substantive question is this: Given a Nahm pole solution as defined here (it is assumed only to obey (3.1)), is there an isomorphism $\upsilon$ from T*Y to ad(P), defined where t is small on $(0, \infty) \times Y$, such that the corresponding $\hat{\mathfrak{e}}$ and $\mathfrak{B}$ obey (A.4)? As explained in the seven parts that follow, the (A.4) condition is obeyed for sufficiently small t when $\upsilon$ is the isomorphism $\tau$ in (3.4), in which case, $\hat{\mathfrak{e}}$ is $\hat{\mathfrak{c}} = \langle \tau \otimes \mathfrak{c} \rangle$ and $\mathfrak{B}$ is $\mathscr{b}$.

*Part 1*: With regards to the argument for $|\mathscr{b}|$: An important point is that $|\nabla \mathscr{b}|$ is bounded by $c_0 \frac{1}{t^2}$. This can be proved using (3.42) and (3.46) and what is said by Proposition 2.1: Rescale around any given point where t is small so that the radius $c_0^{-1} t$ ball about the point is rescaled to have radius 1. Also, rescale $\mathfrak{a}$ to $\hat{\mathfrak{a}} = t\,\mathfrak{a}$. The rescaled versions of (3.42) and (3.46) with the derivative bounds from Proposition 2.1 (and the $|\mathscr{b}| \le c_0 \frac{1}{t}$ bound from Lemma 3.1) lead to a $c_0$ bound for the rescaled version of $|\nabla \mathscr{b}|$. (This uses standard elliptic regularity arguments.) Undoing the rescaling gives the desired $c_0 \frac{1}{t^2}$ bound for $|\nabla \mathscr{b}|$.

*Part 2*: Define $\mathscr{c}$ and $\mathfrak{c}^+$ as in (3.6). Let $\mathfrak{h}$ denote either $\mathscr{c}$ or $\mathfrak{c}^+$ or $\mathscr{b}^\perp$ or $\mathscr{b}_t$; so each is an incarnation of $\mathfrak{h}$. Suppose that there exists $t_0 > 0$ and $\mu > 0$ such that

$$\int_{(0,t]\times Y} \frac{1}{s^2} |\mathfrak{h}|^2 \le t^\mu \tag{A.5}$$

for all t sufficiently small. Given (A.5), suppose in addition that $|\mathfrak{h}| \ge t^{-1+\mu/8}$ at a point in $[\frac{1}{4}t, \frac{1}{2}t] \times Y$ to generate some nonsense. If this is the case, then (because $|\nabla\mathfrak{h}| \le c_0 \frac{1}{t^2}$ for each incarnation of $\mathfrak{h}$), there will be a ball in $[\frac{1}{8}t, t] \times Y$ of radius greater than $c_0^{-1} t^{1+\mu/8}$ where $|\mathfrak{h}|$ is greater than $\frac{1}{2}\, t^{-1+\mu/8}$. The contribution to the left hand side of (A.5) from this ball will be greater than $c_0^{-1} t^{3\mu/4}$ which is much greater than $t^\mu$ when t is small. This catastrophe is inevitable unless $|\mathfrak{h}|$ is less than $t^{-1+\mu/8}$ for all t sufficiently small. The latter bound for all of the incarnations of $\mathfrak{h}$ is what is required by (A.4) (with $\lambda$ being any positive number less than $\mu/8$).

*Part 3*: Certain properties of the $\{t\}\times Y$ integrals of $\langle \mathfrak{a}, \mathcal{B}\rangle$ for t near zero will be used momentarily to find a $t_0 > 0$ and a $\mu > 0$ for (A.5) in each incarnation of $\mathfrak{h}$. To set the stage for this, take t very small and write the $\{t\}\times Y$ integral of $\langle \mathfrak{a}, \mathcal{B}\rangle$ as

$$\int_{\{t\}\times Y} \langle \mathfrak{a}, \mathcal{B}\rangle = \int_{\{t\}\times Y} \langle \mathfrak{a} \wedge \tfrac{1}{4}\varepsilon_{ijk} (d_\Gamma \mathscr{b}^\perp + \mathscr{b}^\perp \wedge \mathscr{b}^\perp)_{ij} \tau_k \rangle \ . \tag{A.6}$$

Now integrate by parts to take the derivative of $\mathscr{b}^\perp$ so as to write (A.6) as

$$\int_{\{t\}\times Y} \langle \mathfrak{a}, \mathcal{B}\rangle = \int_{\{t\}\times Y} \langle d_A \mathfrak{a} \wedge \tfrac{1}{4}\varepsilon_{ijk} (\mathscr{b}^\perp)_{ij} \tau_k \rangle - \int_{\{t\}\times Y} \langle \mathfrak{a} \wedge \tfrac{1}{4}\varepsilon_{ijk} (\mathscr{b}^\perp \wedge \mathscr{b}^\perp)_{ij} \tau_k \rangle \tag{A.7}$$

Since $|\mathfrak{a}| \le c_0 \frac{1}{t}$ and since $|d_A \mathfrak{a}| \le c_0 |\nabla_A{}^\perp \mathfrak{a}|$, the latter bound leads to this one:

$$\Big| \int_{\{t\}\times Y} \langle \mathfrak{a}, \mathcal{B}\rangle \Big| \le c_0 \int_{\{t\}\times Y} |\nabla_A^\perp \mathfrak{a}| |\mathscr{b}^\perp| + c_0 \tfrac{1}{t} \int_{\{t\}\times Y} |\mathscr{b}^\perp|^2 \ ; \tag{A.8}$$

which leads in turn to

$$| \int_{\{t\}\times Y} \langle \mathfrak{a}, \mathcal{B} \rangle \,| \leq c_0 t \int_{\{t\}\times Y} |\nabla_A^\perp \mathfrak{a}|^2 \tag{A.9}$$

because of the second bullet in (3.13).

*Part 4*: Fix $t_1 \in (0, t_\varepsilon)$ is small enough so that (A.9) holds on $(0, t_1]$. Introduce by way of notation $f$ to denote the function on $(0, t_1]$ that is defined by the rule

$$f(t) = \int_{(0,t]\times Y} (|B_A|^2 + |E_A|^2 + |\nabla_A^\perp \mathfrak{a}|^2)\,. \tag{A.10}$$

Now, by virtue of (A.9) and Lemma 3.6, this function and its derivative obey

$$t \tfrac{d}{dt} f \geq c_0^{-1} f - c_0 t\,. \tag{A.11}$$

Write the version of $c_0$ in this inequality as $\frac{1}{2\mu}$. Integrating (A.11) leads to the bound

$$f(t) \leq (\tfrac{t}{t_1})^{2\mu} (f(t_1) + c_0 t_1^{\,2})\,. \tag{A.12}$$

And, (A.12) implies this: There is a time $t_0 \in (0, t_1)$ such that if $t < t_0$, then

$$\int_{(0,t]\times Y} (|B_A|^2 + |E_A|^2 + |\nabla_A^\perp \mathfrak{a}|^2) \leq t^\mu\,. \tag{A.15}$$

(Note that the number $\mu$ does not depend the Nahm pole solution $(A, \mathfrak{a})$.)

*Part 5*: The bound in (A.15) leads directly to the bound in (A.5) when $\mathfrak{h} = \mathcal{b}^\perp$ because of the second bullet in (3.12). To obtain the $\mathfrak{h} = \mathfrak{c}^+$ version of (A.5), it is sufficient to bound the $\{t\}\times Y$ integral of $|\mathfrak{c}^+|^2$ by $z t^{1+\mu}$ when t is small with $z$ being independent of t. Granted this, and granted (3.8), it is sufficient to bound the $\{t\}\times Y$ integral of $|\mathfrak{t}|^2$ by $z t^{-1+\mu}$ when t is small (with $z$ here being independent of t also). To do the latter task, return to (3.27). Let $t_{1\varepsilon}$ now denote a positive time (but less than $t_\varepsilon$) where (A.15) holds. Supposing that $t \in (0, t_{1\varepsilon})$, then integrating (3.27) from t to $t_{1\varepsilon}$ using now the bound from (A.15) leads to this bound:

$$( \int_{\{t\}\times Y} |\mathfrak{t}|^2 )^{1/2} \le ( \int_{\{t_{\varepsilon 1}\}\times Y} |\mathfrak{t}|^2 )^{1/2} + c_0\, t^{-(1-\mu)/2}$$

(A.16)

to replace the bound in (3.28). The latter gives the required bound $c_0\, t^{-1+\mu}$ for the $\{t\}\times Y$ integral of $|\mathfrak{t}|^2$ because the $\{t_{\varepsilon 1}\}\times Y$ of $|\mathfrak{t}|^2$ on the left hand side of (A.16) is t-independent.

*Part 6*: The $\mathfrak{h} = c$ version of (A.5) will hold if the $\{t\}\times Y$ integral of $c^2$ is bounded by $z\,t^{1+\mu}$ when t is small with $z$ being independent of t. To obtain this bound, use the $c_0 t^{1+\mu}$ bound for the $\{t\}\times Y$ integral of $|\mathfrak{c}^+|^2$ from the preceding paragraph to go from the inequality in (3.33) to the following inequality (it replaces (3.34)):

$$\frac{\partial}{\partial t}(t^3 \int_{\{t\}\times Y} c^2 ) \le c_0 t^4 \int_{\{t\}\times Y} |B_A|^2 + c_0\, \varepsilon\, t^{3+\mu} \ .$$

(A.17)

Integrating the latter on (0, t) leads to the desired bound.

*Part* 7: This part considers (A.5) for the case $\mathfrak{h} = \mathscr{b}_t$. Just the $|B_A|^2$ bound in (A.15) with (3.41) and the bounds from Parts 5 and 6 imply that

$$\int_{(0,t]\times Y} |\nabla_{A_t}\mathfrak{a} - \tfrac{1}{2t^2}\mathfrak{t}|^2 \le c_0\, t^{\mu}$$

(A.18)

for t sufficiently small. This last bound with (3.13) then gives

$$\int_{(0,t]\times Y} (\tfrac{1}{s^2}|\mathscr{b}_t|^2 + |\tfrac{\partial}{\partial t} c|^2 + |\nabla_t \mathfrak{c}^+|^2) \le c_0\, t^{\mu} \ .$$

(A.19)

The $|\mathscr{b}_t|^2$ part of this integral gives the required bound.